\newif\ifHAL
\HALtrue

\ifHAL
\documentclass[10pt,a4paper]{article}
\else
\documentclass[10pt]{m2an}
\fi

\usepackage[colorlinks,citecolor=cyan,linkcolor=magenta]{hyperref}
\usepackage[latin1]{inputenc}
\usepackage{amsmath}
\usepackage{amsfonts}
\usepackage{amssymb}
\usepackage{subcaption}
\usepackage{graphicx}
\usepackage{xspace}
\usepackage{tikz}
\usepackage{fullpage}

\usepackage[colorinlistoftodos]{todonotes}

\usepackage{pict2e}
\usepackage{picture}

\newcommand{\RR}{\mathbb{R}}      

\newcommand{\vertiii}[1]{{\|\kern-0.25ex | #1
    | \kern-0.25ex \|}}

\newcommand{\mean}[1]{\{\kern-1.1mm\{#1\}\kern-1.1mm\}}                  

\newcommand{\ndg}[1]{| \kern -.25mm \|{#1}| \kern -.25mm \|}
\newcommand{\nsdg}[1]{| \kern -.25mm \|{#1}| \kern -.25mm \|_{\rm s}}

\newcommand{\su}{\sum_{K\in \mesh}}

\newcommand{\vx}{{\tiny\textbullet } }

\newcommand{\fes}{\hat{V}_{h}^k}

\newcommand{\fesE}{\hat{V}_{K}^k}

\newcommand{\n}{{\bf n}}

\newcommand{\Fall}{\mathcal{F}_h}

\newcommand{\upi}{^{\mathrm{i}}}
\newcommand{\upb}{^{\mathrm{b}}}
\newcommand{\Fb}{\mathcal{F}_h\upb}
\newcommand{\Fint}{\mathcal{F}_h\upi}
\newcommand{\dK}{\partial K}
\newcommand{\dKi}{\partial K\upi}
\newcommand{\dKb}{\partial K\upb}
\newcommand{\FK}{\mathcal{F}_{\dK}}
\newcommand{\Ihk}{\mathcal{\hat{I}}_h^k}
\newcommand{\FKi}{\mathcal{F}_{\dK\upi}}
\newcommand{\FKb}{\mathcal{F}_{\dK\upb}}
\newcommand{\meshi}{\mesh\upi}
\newcommand{\meshb}{\mesh\upb}

\newcommand{\mesh}{\mathcal{T}_h}

\newcommand{\ncdg}[1]{| \kern -.25mm \|{#1}| \kern -.25mm \|_{\rm DG}}

\renewcommand{\tilde}[1]{\widetilde{#1}}
\renewcommand{\hat}[1]{\widehat{#1}}
\makeatletter
\newcommand*{\rom}[1]{\text{\expandafter\@slowromancap\romannumeral #1@}}
\makeatother

\ifHAL
\usepackage{theorem}
\newtheorem{corollary}{Corollary}[section]
\newtheorem{lemma}[corollary]{Lemma}
\newtheorem{theorem}[corollary]{Theorem}
\newtheorem{proposition}[corollary]{Proposition}
\newtheorem{definition}[corollary]{Definition}

\newtheorem{remark}[corollary]{Remark}

\newtheorem{assumption}[corollary]{Assumption}

\newcommand{\qed}{ \vspace{-0.5cm} \hfill $\Box$ }

\else
\theoremstyle{plain}
\newtheorem{theorem}{Theorem}[section]
\newtheorem{remark}[theorem]{Remark}
\newtheorem{lemma}[theorem]{Lemma}
\newtheorem{corollary}[theorem]{Corollary}

\newtheorem{assumption}[theorem]{Assumption}

\fi

\begin{document}

\title{Hybrid high-order method for singularly perturbed fourth-order problems on curved domains}

\ifHAL
\author{
    Zhaonan Dong\thanks{
        Inria, 2 rue Simone Iff, 75589 Paris, France,
        and CERMICS, Ecole des Ponts, 77455 Marne-la-Vall\'{e}e 2, France.
        {\tt{zhaonan.dong@inria.fr}}.
    }
    \and Alexandre Ern\thanks{
        CERMICS, Ecole des Ponts, 77455 Marne-la-Vall\'{e}e 2, France,
        and  Inria, 2 rue Simone Iff, 75589 Paris, France.
        {\tt{alexandre.ern@enpc.fr}}.
}}
\else
\author{Zhaonan Dong}
\address{Inria, 2 rue Simone Iff, 75589 Paris, France;
    \email{zhaonan.dong@inria.fr}}
    \secondaddress{CERMICS, Ecole des Ponts, 77455 Marne-la-Vall\'{e}e 2, France. \email{alexandre.ern@enpc.fr} }
\author{Alexandre Ern}
\sameaddress{2, 1}
\fi

\date{\today}

\ifHAL
\else

\begin{abstract}
    We propose a novel hybrid high-order method (HHO) to approximate singularly perturbed
    fourth-order PDEs on domains with a possibly curved boundary. The two key ideas in
    devising the method are the use of a Nitsche-type boundary penalty technique to weakly
    enforce the boundary conditions and a scaling of the weighting parameter in the
    stabilization operator that compares the singular perturbation parameter to the
    square of the local mesh size. With these ideas in hand, we derive stability and
    optimal error estimates over the whole range of values for the singular perturbation
    parameter, including the zero value for which a second-order elliptic problem is recovered.
    Numerical experiments illustrate the theoretical analysis.
\end{abstract}


\fi

\maketitle

\ifHAL
\begin{abstract}
    We propose a novel hybrid high-order method (HHO) to approximate singularly perturbed
    fourth-order PDEs on domains with a possibly curved boundary. The two key ideas in
    devising the method are the use of a Nitsche-type boundary penalty technique to weakly
    enforce the boundary conditions and a scaling of the weighting parameter in the
    stabilization operator that compares the singular perturbation parameter to the
    square of the local mesh size. With these ideas in hand, we derive stability and
    optimal error estimates over the whole range of values for the singular perturbation
    parameter, including the zero value for which a second-order elliptic problem is recovered.
    Numerical experiments illustrate the theoretical analysis.
\end{abstract}
\else
\fi

\section{Introduction} \label{Introduction}

Fourth-order singular perturbed PDEs are used in the modeling of various physical phenomena,
such as thin plate elasticity, micro-electromechanical systems, and phase separation to mention a few examples. In the present work, we consider the following model problem: Find $u:\Omega\to\mathbb{R}$ such that
\begin{equation}\label{pde}
\left\{\begin{alignedat}{2}
\varepsilon \Delta^2 u -\Delta u &=f &\qquad &\text{in $\Omega$}, \\
u & = g_{\rm D}  &\qquad &\text{on $\partial \Omega$}, \\
\varepsilon \n_\Omega{\cdot}\nabla u  &= \varepsilon g_{\rm N}  &\qquad &\text{on $\partial \Omega$},
\end{alignedat}\right.
\end{equation}
where $\Omega$ is a open bounded Lipschitz domain in $\mathbb{R}^d$, $d=2,3$,
with boundary $\partial\Omega$ and unit outward normal $\n_\Omega$. The problem
data are the forcing term $f:\Omega\to\mathbb{R}$ and the Dirichlet and Neumann data
$g_{\rm D},g_{\rm N}:\partial\Omega\to\mathbb{R}$. The assumptions on the data are
specified below. {The use of other boundary conditions in \eqref{pde} is
currently under study.} Moreover, the perturbation parameter $\varepsilon$ is a nonnegative real number,
i.e., we only assume that $\varepsilon\ge0$, and we are especially interested in
the singularly perturbed regime where $\varepsilon\ll \ell_\Omega^2$, where $\ell_\Omega$ is some suitable length scale associated with $\Omega$, e.g., its diameter ($\ell_\Omega=1$ if the problem is written in nondimensional form).
Notice that the Neumann boundary condition is scaled by $\varepsilon$
so that the model problem~\eqref{pde}
becomes the Poisson problem with Dirichlet boundary condition when $\varepsilon=0$.
Another feature of interest here
is that the domain $\Omega$ can have a curved boundary.

The purpose of this work is to design and analyze a hybrid high-order (HHO) method to
approximate the model problem~\eqref{pde}. The key feature of the proposed method is its
ability to handle in a robust way the whole scale for the singular perturbation parameter
$\varepsilon\in [0,\ell_\Omega^2]$ (notice that the value $\varepsilon=0$ is allowed).
HHO methods were introduced in \cite{DiPEL:14} for linear diffusion and in \cite{DiPEr:15} for locking-free linear elasticity. In such methods, discrete unknowns are attached to the mesh cells and to the mesh faces.
The two key ingredients to devise HHO methods are a local reconstruction operator
and a local stabilization operator in each mesh cell.
HHO methods offer various attractive features, such as the support of polytopal meshes,
optimal error estimates, local conservation properties, and computational efficiency
due to compact stencils and local elimination of the cell unknowns by static condensation.
As a result, these methods have been developed extensively over the past few years and now cover a broad range of applications; we refer the reader to the two recent monographs \cite{di2020hybrid,cicuttin2021hybrid} for an overview. As shown in \cite{CoDPE:16}, HHO methods can be embedded into the broad framework of hybridizable discontinuous Galerkin (HDG) methods, and they can be bridged to nonconforming virtual element methods (ncVEM).
Moreover, HHO methods are closely related to weak Galerkin (WG) methods. Indeed,
the reconstruction operator in the HHO method corresponds to the weak gradient
(or any other differential operator) in WG methods, so that the only relevant
difference between HHO and WG methods lies in the choice of the discrete
unknowns and the design of the stabilization operator.

Various HHO methods for the biharmonic operator were devised and analyzed recently in \cite{DongErn2021biharmonic}, including a comparison with existing WG methods for the biharmonic operator. We also refer the reader to \cite{BoDPGK:18} for the first HHO method for the biharmonic operator in primal form. In \cite{DongErn2021biharmonic}, two HHO methods were proposed (called HHO-A and HHO-B). Both methods use cell unknowns to approximate the solution in each mesh cell, face unknowns to approximate its trace on the mesh faces, and face unknowns to approximate its normal derivatives on the mesh faces. In both methods, the cell unknowns are polynomials of degree $(k+2)$ and the face unknowns for the normal derivative
are polynomials of degree $k$, with $k\ge0$. HHO-A is restricted to two space dimensions
and uses polynomials of degree $(k+1)$ for the face unknowns related to the trace, whereas
HHO-B supports any space dimension but uses polynomials of degree $(k+2)$ for these face unknowns. Moreover, the HHO-A method was combined in \cite{DongErn2021biharmonic} with a Nitsche-type boundary penalty technique, originally introduced in \cite{CaChE:20} to weakly enforce Dirichlet conditions in HHO methods for second-order PDEs and further developed in \cite{BurEr:18,BurErCiDe:21} to handle unfitted meshes in problems with a curved interface or boundary. In particular, one of the advances in \cite{BurErCiDe:21} is that the weighting parameter in the boundary penalty term does not need to be large enough, but only positive.

In the present work, our starting point is the HHO-B method from \cite{DongErn2021biharmonic}. Consistently with the paradigm considered for singularly perturbed second-order elliptic PDEs, the boundary conditions in \eqref{pde} are weakly enforced by means of a Nitsche-type boundary penalty technique. This is the first key idea to capture possible boundary layers and to achieve robustness for the singularly perturbed fourth-order elliptic problem.
The second key idea to achieve robustness is to revisit the weighting of the
stabilization operator in the HHO-B method by including a scaling factor that
compares the singular
perturbation parameter $\varepsilon$ with the square of the (local) mesh size.
With these two ideas in hand, we can devise a novel HHO method that remains
uniformly stable over the full range $\varepsilon \in [0,\ell_\Omega^2]$ and that delivers
optimally decaying error estimates, both in the case $\varepsilon\approx\ell_\Omega^2$
(representative of a fourth-order PDE) and in the case $\varepsilon\ll\ell_\Omega^2$
and even $\varepsilon=0$ (representative of a second-order PDE). In a nutshell
(see Theorem~\ref{Theorem: main} for a more precise statement and the remarks below
for a discussion), the error estimate
takes the general form
$\varepsilon^{\frac12}\|\nabla^2 e\|_\Omega + \|\nabla e\|_\Omega
\le C(\varepsilon^{\frac12} h^{k+1}+ h^{k+2})$,
where $e$ represents the approximation error, $h$ the mesh size, and the constant
$C$ depends on the regularity of the exact solution, the shape-regularity of the
underlying meshes, and the polynomial degree $k\ge0$.
An additional benefit of using the Nitsche-type boundary penalty technique is the
seamless support of domains with a curved boundary, in the wake of the ideas
developed in \cite{BurEr:18,BurErCiDe:21} for second-order PDEs.

Let us briefly put our contribution in perspective with the literature on
singularly perturbed fourth-order PDEs. Consistently with the present approach,
we focus on discretization methods that hinge on the primal form of the PDE and thus lead,
at the algebraic level, to a symmetric positive definite linear system.
To the best of our knowledge, the present method appears to be the first in the literature
that, at the same time, supports polytopal meshes and offers a robust behavior over the
full range of values for the singular perturbation parameter $\varepsilon$.
On the one hand, robust approximation methods developed on specific meshes
(composed, e.g., of simplices or cuboids) include {$C^0$-interior penalty
discontinuous Galerkin (IPDG) methods \cite{brenner2011c}}
and methods based on the modified Morley element
\cite{nilssen2001robust,wang2006modified,wang2007robust,guzman2012family,wang2013uniformly},
for which a weak enforcement of the boundary conditions using Nitsche-type techniques
was considered more recently in \cite{wang2018morley,huang2021morley}.
On the other hand, discretization methods for singularly perturbed fourth-order operators
on polytopal meshes include the WG method from \cite{cui2020uniform}
and the $C^0$-ncVEM from \cite{zhang2020nonconforming}. Both methods, however,
do not support the limit with $\varepsilon=0$ and in this case lead, at the
algebraic level, to a singular linear system. {More precisely, in this limit,
the sub-blocks coupling the face unknowns discretizing the gradient (or the normal gradient)
on the faces either to the other unknowns or to themselves all vanish. Thus, to recover a nonsingular linear system when $\varepsilon=0$, one needs to manually remove these gradient face unknowns, but unfortunately this fix cannot be applied when $\varepsilon\ll1$, leading to serious conditioning issues in this case. This situation is instead avoided by the present method: by including a mesh-dependent cutoff in the stabilization coefficient coupling the gradient face unknowns to the other unknowns, the linear system at the limit $\varepsilon=0$ remains nonsingular without any need to remove manually some unknowns.}

The rest of this work is organized as follows.  We present the weak formulation of the model problem together with the discrete setting in Section \ref{sec: weak form and discrete}.  In Section \ref{sec: HHO-N}, we introduce the present HHO method.  In Section \ref{sec:Stability and error analysis}, we present our main results on the stability and error analysis of the HHO method. Numerical results are discussed in Section \ref{sec: numerical examples}. Finally, the proofs of our main results are collected in Section \ref{sec: proof of main resutls}.

\section{Weak formulation and discrete setting}\label{sec: weak form and discrete}

In this section, we present the weak formulation of the model problem~\eqref{pde} together with the discrete setting.

\subsection{Weak formulation} \label{Basic notation and weak formulation}

We use standard notation for the Lebesgue and Sobolev spaces. In particular, when considering fractional-order Sobolev spaces, we use the Sobolev--Slobodeckij seminorm based on the double integral. For an open, bounded, Lipschitz set $S$ in $\mathbb{R}^d$, $d\in\{1,2,3\}$, with a piecewise smooth boundary, we denote by $(v,w)_S$ the $L^2(S)$-inner product, and we employ the same notation when $v$ and $w$ are vector- or matrix-valued fields. We denote by $\nabla w$ the (weak) gradient of $w$ and by $\nabla^2 w$ its (weak) Hessian. It is convenient to consider the following inner product and corresponding seminorm on $\varepsilon H^2(S)+H^1(S)$ (this is just a shortcut notation for $H^2(S)$ if $\varepsilon>0$ and $H^1(S)$ if $\varepsilon=0$):
\begin{equation} \label{eq:def_inner_prod}
(\nabla v, \nabla w)_{S,\varepsilon} : = \varepsilon(\nabla^2 v, \nabla^2 w)_{S} + (\nabla v, \nabla w)_{S}, \qquad  \|\nabla v\|_{S,\varepsilon}^2 : = (\nabla v, \nabla v)_{S,\varepsilon}.
\end{equation}

Let $\n_S$ be the unit outward normal vector on the boundary $\partial S$ of $S$. Assuming that the functions $v$ and $w$ are smooth enough, we have the following integration by parts formula:
\begin{align}
\varepsilon(\Delta^2 v,w)_S - (\Delta v,w)_S={}& (\nabla v, \nabla w)_{S,\varepsilon} + \varepsilon \big((\nabla\Delta v,\n_S w)_{\partial S}-(\nabla^2v\n_S,\nabla w)_{\partial S}\big)
-(\nabla v,\n_S w)_{\partial S}. \label{eq:ipp1}
\end{align}
To alleviate the notation, it is implicitly understood that within integrals over $\partial S$, $\partial_n$ denotes the normal derivative on $\partial S$ along $\n_S$. Moreover, $\partial_t$ denotes the ($\RR^{d-1}$-valued) tangential derivative on $\partial S$.
We also denote by $\partial_{nn} v$ the (scalar-valued) normal-normal second-order derivative and by $\partial_{nt} v$ the ($\RR^{d-1}$-valued) normal-tangential second-order derivative. The integration by parts formula~\eqref{eq:ipp1} can then be rewritten as
\begin{align}
\varepsilon(\Delta^2 v,w)_S - (\Delta v,w)_S ={}& (\nabla v, \nabla w)_{S,\varepsilon} + \varepsilon \big((\partial_n\Delta v,w)_{\partial S}-(\partial_{nn}v,\partial_n w)_{\partial S}-(\partial_{nt}v,\partial_t w)_{\partial S}\big) -(\partial_n v, w)_{\partial S}. \label{eq:ipp2}
\end{align}

Let us consider the Hilbert spaces $V:= \varepsilon H^2(\Omega) + H^1(\Omega)$ and $V_0:= \varepsilon H^2_0(\Omega) + H^1_0(\Omega)$ equipped with the inner product $(\nabla v, \nabla w)_{\Omega,\varepsilon}$. Assume that the source term in \eqref{pde} satisfies $f\in L^2(\Omega)$ and that the boundary data are such that there is $u_g\in V$ such that
$u_g = g_{\rm D}$ and $\varepsilon \n{\cdot}\nabla u_g = \varepsilon g_{\rm N}$
on $\partial \Omega$.
Using the above integration by parts formula, the following weak formulation for \eqref{pde} is derived: Find  $u\in u_g+V_0$ such that
\begin{equation}\label{weak form}
(\nabla u, \nabla w)_{\Omega,\varepsilon} = (f,w)_\Omega, \qquad  \forall w \in V_0.
\end{equation}
The Lax--Milgram lemma readily shows that this problem is well-posed.

\subsection{Polytopal and curved meshes}\label{Polytopic  meshes}

In this work, we assume that
$\partial\Omega$ can be covered by a finite number of closed $C^2$ manifolds with
nonoverlapping interior, and we write $\partial\Omega=\bigcup_{m\in\{1,\ldots,N_\partial\}}S_m$.
Let $\{\mesh\}_{h>0}$ be a mesh family such that each mesh $\mesh$ covers $\Omega$ exactly.
A generic mesh cell is denoted by $K\in\mesh$, its diameter by $h_K$, and
its unit outward normal by $\n_K$. We define the following mesh-dependent
parameter for measuring locally the dominant operator in the PDE: For all $K\in\mesh$,
\begin{equation}\label{Def:sigma}
\sigma_K := \max\{1,\varepsilon h^{-2}_K\}.
\end{equation}

We partition the boundary $\dK$ of any mesh cell $K\in\mesh$
by means of the two subsets $\dKi :=  \overline{\dK\cap \Omega}$ and
$\dKb:=\dK\cap \partial\Omega$. Similarly, we partition the mesh as
$\mesh = \meshi \cup \meshb$, where $\meshb$ is the collection of all the mesh
cells $K$ such that $\dKb$ has positive measure.
The mesh faces are collected in the set $\Fall$, which is split as
$\Fall=\Fint\cup\Fb$, where $\Fint$ is the collection of the interior faces
(shared by two distinct mesh cells) and
$\Fb$ the collection of the boundary faces. For all $F\in \Fall$, we orient $F$
by means of the fixed unit normal vector $\n_F$ whose direction is arbitrary
for all $F\in\Fint$ and $\n_F:=\n_\Omega$ for all $F\in\Fb$.
For any mesh cell $K\in\mesh$, the mesh faces composing its boundary
$\partial K$ are collected in the set $\FK$, which is partitioned as
$\FK=\FKi\cup \FKb$ with obvious notation.
To avoid distracting technicalities, we assume that each mesh $\mesh$ is
compatible with the decomposition $\partial\Omega=\bigcup_{m\in\{1,\ldots,N_\partial\}}S_m$,
so that for all $K\in\meshb$, each face $F\in\FKb$ is a closed $C^2$ manifold.

In this work, we consider mesh sequences satisfying the following mesh shape-regularity assumption.

\begin{figure}[h]
\centering
	\setlength{\unitlength}{3cm}
    \includegraphics[scale=0.4]{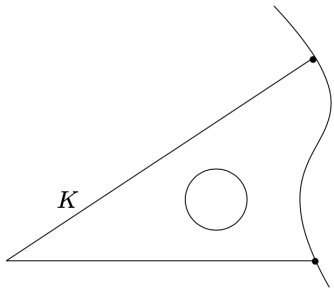}
%
%
%
%
	\caption{Example of a boundary mesh cell with a curved boundary face and which is star-shaped with respect an interior ball; \vx denotes a mesh vertex at the boundary.}
	\label{fig:curved_star-shaped}
\end{figure}

\begin{assumption}[Mesh shape-regularity]\label{mesh assumption}
(i) Any interior mesh cell $K\in\meshi$ is a polytope with planar faces, and the sequence of interior meshes $\{\meshi\}_{h>0}$ is shape-regular in the sense of {\cite[Definition 1]{DiPEr:15}}.
(ii) For any boundary mesh cell $K\in\meshb$, all the faces in $\FKi$ are planar with diameter uniformly equivalent to $h_K$, and all the faces in $\FKb$ are subsets of $\partial\Omega$ which are closed $C^2$ manifolds. Moreover, for each $F\in\FKb$, $K$ can be decomposed into a finite union of nonoverlapping subsets, $\{K_{F,m}\}_{m\in \{1,\ldots,n_{K,F}\}}$, so that each $K_{F,m}$ is star-shaped with respect to an interior ball of radius uniformly equivalent to $h_K$; see Figure~\ref{fig:curved_star-shaped} for an illustration with $n_{K,F}=1$. (Notice that the star-shapedness assumption implies that $n_{K,F}$ is uniformly bounded.)
\end{assumption}

\begin{remark}[Mesh assumptions]
The above assumptions on the mesh sequence are fairly general. Let us
briefly discuss some of the most significant ones.
(i) The assumption that each mesh $\mesh$ covers $\Omega$ exactly is reasonable
in the present context where boundary conditions are enforced by means of a
Nitsche-like penalty method. In particular, all the integrals in the mesh cells
and their faces are performed in the physical space without invoking a geometric mapping
that can introduce some error due to the approximation of the geometry. To avoid
distracting technicalities, we do not consider quadrature errors in our error analysis.
(ii) The assumption that the interior faces of the mesh are planar is important in
the present setting of HHO methods which use polynomial functions as discrete
unknowns attached to these interior faces. Notice though that the use of
a Nitsche-like penalty method allows us to avoid introducing discrete unknowns on
the boundary faces; thus, such faces
do not need to be planar.
(iii) The star-shapedness assumption on the mesh boundary cells is introduced to
invoke a Poincar\'e-type inequality in such cells (and, more generally, polynomial
approximation properties; see Lemma~\ref{lemma: Polynomial approximation} below).
An alternative is to invoke an extension operator when asserting polynomial
approximation properties, as, for instance, in \cite{BurEr:18}.
In this case, only the multiplicative trace inequality
(see Lemma~\ref{lemma: trace inequality} below) requires a star-shapedness assumption
for the mesh boundary cells, but for this latter result to hold, star-shapedness with respect to an interior point for each $K_{F,m}$ is a sufficient assumption \cite[Lemma 32]{DGpolybook}
(and in this setting, the number $n_{K,F}$ does no longer need to
be uniformly bounded).
\end{remark}

\subsection{Analysis tools}\label{Inverse inequality and polynomial approximation}

Let us briefly review the main analysis tools used in this work. We simply state the results and refer the reader to Remark~\ref{rem:proofs} for some comments on the proofs.
In what follows,
we always consider a shape-regular mesh sequence satisfying Assumption~\ref{mesh assumption}.
Moreover, in various bounds, we use
the symbol $C$ to denote any {positive} generic constant (its value can change at each occurrence)
that is independent of $h>0$, the considered mesh cell $K\in\mesh$, and the considered
function in the inequality. The value of $C$ can depend on the parameters quantifying
the shape-regularity of the mesh sequence and the polynomial degree (whenever relevant).

\begin{lemma}[Discrete inverse inequalities]\label{lemma: Inverse inequality}
Let $l\ge0$ be the polynomial degree.
There is $C$ (depending on $l$)
such that for all $h>0$, all $K\in\mesh$, and all $v_h\in\mathbb{P}^l(K)$,
\begin{align}\label{trace inverse inequality}
\|{v}_h\|_{\partial K} &\leq C h_K^{-\frac12}  \|{v}_h\|_{ K}, \\
\label{H1 inverse inequality}
\|\nabla {v}_h\|_{K} &\leq C h_K^{-1}  \|{v}_h\|_{ K}, \\
\label{inverse inequality tangential}
\|\partial_t {v}_h\|_{F} &\leq C h_K^{-1}  \|{v}_h\|_{ F}, \quad \forall F\in\FKi.
\end{align}
\end{lemma}

\begin{lemma}[Multiplicative trace inequality]\label{lemma: trace inequality}
There is $C$ such that for all $h>0$,
all $K\in\mesh$, and all $v\in H^1(K)$,
\begin{equation}\label{trace inequality}
\|{v} \|_{\dK} \leq C \big(
h_K^{-\frac12}  \|{v}\|_{K} + h_K^{\frac12}  |v|_{H^1(K)}\big).
\end{equation}
\end{lemma}

\begin{lemma}[Polynomial approximation]\label{lemma: Polynomial approximation}
Let $l\geq 0$ be the polynomial degree.
There is $C$ (depending on $l$) such that for all $r\in [0,l+1]$,
all $m\in \{0,\ldots,\lfloor r\rfloor\}$, all $h>0$, all $K\in\mesh$,
and all $v\in H^r(K)$,
\begin{equation}\label{Polynomial approximation}
|{v} - \Pi^{l}_{K} (v)|_{H^m(K)}
\leq C h^{r-m}_K |{v} |_{H^r(K)},
\end{equation}
where $\Pi_K^l$ denotes the $L^2$-orthogonal projection onto $\mathbb{P}^{l}(K)$.
\end{lemma}

Let us briefly highlight some useful consequences of the above results. First,
\eqref{Polynomial approximation} includes the following Poincar\'e-like inequalities
(take, respectively, $l=0$, $r=1$, and $l=1$, $r=2$, $m\in\{0,1\}$
in \eqref{Polynomial approximation}): For all $K\in\mesh$,
\begin{alignat}{2}
&\|v-\Pi_K^0(v)\|_K \le Ch_K\|\nabla v\|_K,&\qquad &\forall v\in H^1(K), \label{eq:PK1}
\\
&\|v-\Pi_K^1(v)\|_K+h_K\|\nabla(v-\Pi_K^1(v))\|_K\le Ch_K^2\|\nabla^2v\|_K,&\qquad &\forall
v\in H^2(K).\label{eq:PK2}
\end{alignat}
The following consequences of \eqref{eq:PK1}-\eqref{eq:PK2} combined with the multiplicative trace inequality \eqref{trace inequality} will be useful in our analysis: Letting $k\ge0$ be the polynomial degree, we have for all $v\in H^2(K)$,
\begin{align}
&h_K^{-\frac12} \|v-\Pi_K^{k+2}(v)\|_{\dK} \le C\|\nabla(v-\Pi_K^{k+2}(v))\|_K, \label{eq:tr_PK1}\\
&{h_K^{-\frac32}} \|v-\Pi_K^{k+2}(v)\|_{\dK} + h_K^{-\frac12}\|\nabla(v-\Pi_K^{k+2}(v))\|_{\dK}
\le C\|\nabla^2(v-\Pi_K^{k+2}(v))\|_K. \label{eq:tr_PK2}
\end{align}
Recalling that $\sigma_K$ is defined in~\eqref{Def:sigma} and the norm $\|{\cdot}\|_{K,\varepsilon}$ in~\eqref{eq:def_inner_prod} (with $S:=K$), we have for all $v\in H^2(K)$,
\begin{equation} \label{eq:tr_sigma}
\sigma_K^{\frac12} h_K^{-\frac12} \|v-\Pi_K^{k+2}(v)\|_{\dK} \le C\|\nabla(v-\Pi_K^{k+2}(v))\|_{K,\varepsilon}.
\end{equation}

\begin{remark}[Proofs] \label{rem:proofs}
Let us briefly comment on the proofs of the above lemmas. Concerning Lemma~\ref{lemma: Inverse inequality} and Lemma~\ref{lemma: trace inequality}, the proof on mesh cells having flat faces can be found, e.g., in \cite[Sec. 1.4.3]{DiPietroErnbook}.
On mesh cells having a curved face, these results are established, e.g., in \cite{BurEr:18,Wu_unfitted} assuming that the curved face is a $C^2$ manifold.
More recently, these results were extended  in \cite{cangiani2019hp} with fully explicit constants to $C^1$ manifolds (and sometimes even Lipschitz) and some mild additional geometric assumptions.

Concerning Lemma~\ref{lemma: Polynomial approximation}, the key step is to establish the
Poincar\'e inequality \eqref{eq:PK1} since \eqref{Polynomial approximation} can then be derived by using recursively the Poincar\'e inequality.
On the interior mesh cells, which can be decomposed as
a finite union of (convex) subsimplices, this latter inequality is established by proceeding as in \cite{VeeserPoicare,ErnGuermond:17}. On the boundary mesh cells, which can have a curved face,
one invokes the star-shapedness assumption with respect to a ball. We refer the reader to \cite{zheng2005friedrichs} for the derivation of this inequality with an explicitly determined constant under such an assumption.
\end{remark}

\section{HHO discretization} \label{sec: HHO-N}

In this section, we first introduce the local ingredients to formulate the HHO discretization in each mesh cell and then we derive the global discrete problem.

\subsection{Local unknowns, reconstruction, and stabilization}\label{sec:local}

Let $k\geq 0$ be the polynomial degree. For all $K\in \mesh$, the local HHO space is
\begin{equation} \label{eq:def_fesE}
\fesE: =\mathbb{P}^{k+2}(K) \times \mathbb{P}^{k+2}(\FKi) \times
\mathbb{P}^{k}(\FKi),
\end{equation}
where $\mathbb{P}^{l}(\FKi) := \times_{F\in\FKi} \mathbb{P}^l(F)$ for all $l\ge0$. Notice that we do not introduce any discrete unknowns on the faces of $K$ that lie on the boundary.
A generic element in  $\fesE$ is denoted $\hat{v}_K = (v_K, v_{\partial K}, \gamma_{\partial K})$ with $v_K \in \mathbb{P}^{k+2}(K)$, $v_{\partial K} \in \mathbb{P}^{k+2}(\dKi)$, and  $\gamma_{\partial K} \in \mathbb{P}^{k}(\dKi)$. The first component of $\hat{v}_K$ aims at representing the solution inside the
mesh cell, the second its trace on the interior part of the cell boundary, and the third its normal derivative on the interior part of the cell boundary (along the direction of the outward normal $\n_K$).

We define the local reconstruction operator $R_K\upi: \fesE \rightarrow \mathbb{P}^{k+2}(K)$ such that, for all $\hat{v}_K \in \fesE$ with $\hat{v}_K:= (v_K,v_{\partial K}, \gamma_{\partial K})$, the polynomial
$R_K\upi(\hat{v}_K)\in \mathbb{P}^{k+2}(K)$ is uniquely defined by solving the following problem with test functions $w \in \mathbb{P}^{k+2}(K)^\perp:=\{w\in \mathbb{P}^{k+2}(K)\;|\; (w,1)_K=0\}$:
\begin{align}
(\nabla R_K\upi (\hat{v}_K), \nabla w)_{K, \varepsilon}
= {}& (v_K, \varepsilon\Delta^2 w)_K-(v_K,\Delta w)_K + (v_{\partial K} , \partial_n  w)_{\dKi} \nonumber \\
& - \varepsilon \Big\{ (v_{\dK} , \partial_n \Delta  w)_{\dKi}
- (\gamma_{\dK},  \partial_{nn}  w)_{\dKi}
- (\partial_t v_{\dK},  \partial_{nt}  w)_{\dKi}\Big\}
, \label{eq:rec_ipp HHO-N}
\end{align}
together with the condition $(R_K\upi (\hat{v}_K),1)_K = (v_K,1)_K$.
Integration by parts shows that~\eqref{eq:rec_ipp HHO-N} is equivalent to
\begin{align}
(\nabla R_K\upi(\hat{v}_K), \nabla w)_{K, \varepsilon}
={}&	(\nabla  {v}_K, \nabla w)_{K,\varepsilon}
- (v_K -v_{\partial K} , \partial_n  w)_{\dKi}
- (v_K  , \partial_n  w)_{\dKb}\nonumber \\
&+ \varepsilon\Big\{ (v_K -v_{\partial K} , \partial_n \Delta  w)_{\dKi}
- (\partial_n v_K - \gamma_{\partial K},  \partial_{nn}  w)_{\dKi}
- (\partial_t (v_K - v_{\partial K}),  \partial_{nt}  w)_{\dKi} \nonumber \\
& + (v_K , \partial_n \Delta  w)_{\dKb}
- (\nabla  v_K, \nabla  \partial_{n}  w)_{\dKb}  \Big\}.\label{eq: reconstruction}
\end{align}
Notice that $(\nabla  v_K, \nabla \partial_{n} w)_{\dKb}=(\partial_nv_K,\partial_{nn}w)_{\dKb} + (\partial_tv_K,\partial_{nt}w)_{\dKb}$.

The local stabilization bilinear form is composed of a contribution on $\dKi$ and one on $\dKb$. These two contributions are defined such that, for all $(\hat{v}_K, \hat{w}_K)\in \fesE \times \fesE$, with  $\hat{v}_K:= (v_K,v_{\partial K}, \gamma_{\partial K})$ and $\hat{w}_K:= (w_K,w_{\partial K}, \chi_{\partial K})$,
\begin{equation}\label{def: stabilisation}
\begin{aligned}
S\upi_{\dK}(\hat{v}_K,\hat{w}_K)
:={}& \sigma_Kh_K^{-1} \big( v_{\partial K}- v_K,w_{\partial K}- {w}_K
\big)_{\dKi} +  \sigma_K h_K
\big( \Pi^{k}_{\dKi}(\gamma_{\partial K}- \partial_n {v}_K),
\chi_{\partial K}- \partial_n {w}_K\big)_{\dKi},
\end{aligned}
\end{equation}
 and
\begin{equation}\label{def: stabilisation Boundary}
\begin{aligned}
S\upb_{\dK}(v_K,w_K)
:= {} & \sigma_Kh_K^{-1}\big( v_K, {w}_K\big)_{\dKb}
+ \varepsilon h_K^{-1} \big( \nabla  {v}_K,\nabla {w}_K\big)_{\dKb},
\end{aligned}
\end{equation}
where $\Pi^k_{\dKi}$ denotes the $L^2$-orthogonal projection onto the broken polynomial space $\mathbb{P}^k(\FKi)$.

Finally, we define the local bilinear form $\hat{a}_K$ on $\fesE \times \fesE$ such that
\begin{equation}\label{local bilinear form}
\hat{a}_K(\hat{v}_K,\hat{w}_K)
:=
(\nabla R_K\upi (\hat{v}_K),\nabla R_K\upi (\hat{w}_K))_{K,\varepsilon}
+S\upi_{\dK}(\hat{v}_K,\hat{w}_K)
+S\upb_{\dK}(v_K,w_K).
\end{equation}

\begin{remark}[Reconstruction]
There are two differences with the reconstruction operator introduced in
\cite{DongErn2021biharmonic} for the biharmonic problem. First, as expected,
the terms related to the second-order operator are added, whereas the terms
related to the fourth-order operator are scaled by $\varepsilon$.
The second difference is more subtle and is inspired from
the ideas in \cite{BurErCiDe:21} for the second-order operator
and extended here to the fourth-order operator as well. It consists in
discarding the integrals over $\dKb$ and only keeping the integrals over $\dKi$
for all the boundary terms on the right-hand side of~\eqref{eq:rec_ipp HHO-N}.
Following the ideas in \cite{CaChE:20}, it is also possible to keep the
boundary terms and to use the trace of the cell
unknown $v_K$ and its normal derivative on $\dKb$ to evaluate them. The
disadvantage of this latter approach is that the Nitsche-type boundary penalty terms
need then to be weighted by a coefficient that is large enough, whereas
the weighting coefficient needs only to be positive in the present setting.
\end{remark}

\begin{remark}[Stabilization]
The interior stabilization bilinear form $S\upi_{\dK}$ is inspired from
\cite{DongErn2021biharmonic} and is weighted here by the local coefficient
$\sigma_K$ defined in~\eqref{Def:sigma} to cover both regimes of interest (dominant Laplacian and dominant bi-Laplacian). Moreover, the boundary stabilization bilinear form $S\upb_{\dK}$ is associated with the Nitsche-type boundary penalty. We emphasize that this latter bilinear form does not need to be weighted by a coefficient which is large enough.
\end{remark}

\subsection{The global discrete problem}

Recall that $k\ge0$ is the polynomial degree and that the local HHO space $\fesE$ is defined in~\eqref{eq:def_fesE} for all $K\in\mesh$. The global HHO space is defined as
\begin{equation} \label{eq:def_fes}
\fes : = \mathbb{P}^{k+2}(\mesh) \times \mathbb{P}^{k+2}(\Fint)\times \mathbb{P}^k(\Fint).
\end{equation}
A generic element in $\fes$ is denoted $\hat{v}_h:=(v_{\mesh},v_{\Fint},\gamma_{\Fint})$
with $v_{\mesh}:=(v_K)_{K\in\mesh}$, $v_{\Fint}:=(v_F)_{F\in\Fint}$, and $\gamma_{\Fint}
:=(\gamma_F)_{F\in\Fint}$, where $\gamma_F$ is meant to approximate the normal derivative in the direction of the unit normal vector $\n_F$ orienting $F$.
For all $K\in \mesh$, the local components of $\hat{v}_h$ are collected in the triple
$\hat{v}_K:=(v_K,v_{\dK},\gamma_{\dK})\in\fesE$ with $v_{\dK}|_F: = v_F$ and  $\gamma_{\dK}|_F: = (\n_F {\cdot}\n_K)\gamma_{F}$ for all $F \in \FKi$.

For all $\hat{v}_h, \hat{w}_h \in \fes$, the global bilinear form $\hat{a}_h$ is assembled cellwise as follows:
\begin{equation}\label{bilinear form}
\hat{a}_h(\hat{v}_h, \hat{w}_h):= \su \hat{a}_K(\hat{v}_K, \hat{w}_K),
\end{equation}
with $\hat{a}_K$ defined in~\eqref{local bilinear form}.
To assemble the right-hand side of the discrete problem, we define the discrete linear form
\begin{align}\label{linear form}
{\ell}_h(\hat{w}_h):= &\su \Big\{
(f,{w}_K)_K
+\big(g_{\rm D},
\sigma_K h_K^{-1} w_{K}
+ \varepsilon\partial_n \Delta (R_K\upi(\hat{w}_{K}))
- \partial_n  (R_K\upi(\hat{w}_{K})) \big)_{\dKb} \nonumber \\
&+\varepsilon \big( g_{\rm N} \n + (\partial_t g_{\rm D}){\bf t},
h_K^{-1}\nabla w_{K} -  \nabla \partial_n (R_K\upi(\hat{w}_{K})) \big)_{\dKb}\Big\}.
\end{align}
The devising of $\ell_h$ is motivated by the consistency error analysis (see the proof of
Lemma \ref{lemma: consistency} in Section~\ref{sec:proof_consistency}).
Finally, the discrete problem consists in finding ${\hat{u}_h}\in \fes$ such that
\begin{equation}\label{discrete problem}
\hat{a}_h(\hat{u}_h,\hat{w}_h) = {\ell}_h (w_h), \qquad \forall w_h\in \fes.
\end{equation}
In the next section, we establish stability and consistency properties for~\eqref{discrete problem}, leading to robust and optimal error estimates. Let us also mention that, at the
algebraic level, the discrete problem~\eqref{discrete problem} is amenable to static condensation: all the cell unknowns can be eliminated locally, leading to a global problem coupling only the face unknowns approximating the trace and the normal derivative of the solution at the mesh interfaces.

\begin{remark}[Limit regime $\varepsilon=0$]
We emphasize that the discrete problem~\eqref{discrete problem} remains well-posed even
in the limit regime where $\varepsilon=0$. The resulting HHO discretization though differs
from the usual HHO discretizations for second-order PDEs. Indeed, taking $\varepsilon=0$
in \eqref{discrete problem}, one still has triples of local unknowns. In other words,
discrete unknowns approximating the normal derivative at the mesh interfaces are
still present and coupled to the other discrete unknowns.
\end{remark}

\begin{remark}[Other HHO method]
In the two-dimensional case, one can also think of using the HHO-A method developed in \cite[section 3]{DongErn2021biharmonic} for the biharmonic operator. The advantage is that the global HHO space can be reduced to $\mathbb{P}^{k+2}(\mesh) \times \mathbb{P}^{k+1}(\Fint)\times \mathbb{P}^k(\Fint)$. However, it is not yet clear how to design  the stabilization bilinear form so as to derive stability and error estimates that remain robust in the singularly perturbed regime $\varepsilon\ll1$.
\end{remark}

\section{Main results}\label{sec:Stability and error analysis}

In this section, we state our main results concerning the analysis of the above HHO method: stability and well-posedness, polynomial approximation and bound on consistency error, and, finally, the main error estimate leading to robust and optimally decaying convergence rates. The proofs of these results are contained in Section~\ref{sec: proof of main resutls}. Recall that in this work, we use
the symbol $C$ in bounds
to denote any {positive} generic constant (its value can change at each occurrence)
that is independent of $h>0$, the considered mesh cell $K\in\mesh$, and the considered
function in the bound. The value of $C$ can depend on the parameters quantifying
the shape-regularity of the mesh sequence and the polynomial degree.
In addition, the value of $C$ is independent of the singular perturbation parameter
$\varepsilon\ge0$.

\subsection{Stability and well-posedness}

We define the local energy seminorm defined such that, for all $K\in\mesh$ and all $\hat{v}_K:=(v_K, v_{\partial K}, \gamma_{\partial K}) \in \fesE$,
\begin{align}
|\hat{v}_K|^2_{\fesE}: ={} & \|\nabla v_K\|_{K,\varepsilon}^2
+ \sigma_K h_K^{-1}\|v_{\partial K} - v_K\|_{\dKi}^2
+ {\sigma_K}{h_K} \|\gamma_{\partial K} - \partial_n v_K\|_{\dKi}^2 \nonumber \\
&+  \sigma_K h_K^{-1} \|v_{K}\|^2_{\dKb}
+ \varepsilon h_K^{-1} \|\nabla v_{K}\|^2_{\dKb}. \label{H2_seminorm_elem}
\end{align}
The proof of the following result is postponed to Section~\ref{sec:proof_stability}.

\begin{lemma}[Local stability and boundedness] \label{lem: stability and boundedness}
There is a real number $\alpha>0$, depending only on the mesh shape-regularity and the polynomial degree $k$, such that, for all $h>0$, all $K \in \mesh$, and all $ \hat{v}_K\in \fesE$,
\begin{equation}\label{local equivalent}
\alpha|\hat{v}_K|^2_{\fesE}
\leq  \|\nabla R_K\upi (\hat{v}_K) \|_{K,\varepsilon}^2
+ S\upi_{\dK}(\hat{v}_K,\hat{v}_K)+ S\upb_{\dK}(v_K,v_K)
\leq 	\alpha^{-1} |\hat{v}_K|^2_{\fesE}.
\end{equation}
\end{lemma}

We equip the space $\fes$ with the norm  $\|\hat{v}_h\|_{\fes}^{2}:= \su |\hat{v}_K|^2_{\fesE}$.
It is readily verified that $\hat{v}_h\mapsto \|\hat{v}_h\|_{\fes}$ indeed defines
a norm on $\fes$. An immediate consequence of Lemma~\ref{lem: stability and boundedness} is the following bound establishing that the discrete bilinear form $\hat{a}_h$ is coercive on $\fes$:
\begin{equation}\label{coercivity}
\hat{a}_h(\hat{v}_h,\hat{v}_h)  \geq \alpha \|\hat{v}_h\|_{\fes}^2, \qquad \forall \hat{v}_h\in \fes.
\end{equation}
Invoking the Lax--Milgram lemma readily yields the following result.

\begin{corollary}[Well-posedness]\label{lemma: coercivit}
The discrete problem \eqref{discrete problem} is well-posed.
\end{corollary}

\subsection{Approximation and consistency} \label{Consistency and error analysis}

For all $K\in \mesh$, we define the local reduction operator $\mathcal{\hat{I}}^k_K : H^2(K) \rightarrow \fesE$ such that, for all $v\in H^2(K)$,
\begin{equation}\label{interpoltation}
	\mathcal{\hat{I}}^k_K(v):= (\Pi_{K}^{k+2}(v),\Pi_{\dKi}^{k+2} (v), \Pi_{\dKi}^k (\n_K{\cdot}\nabla v)) \in \fesE.
\end{equation}
In addition, we define the operator $\mathcal{E}_K\upi:= R_K\upi \circ \mathcal{\hat{I}}^k_K: H^2(K) \rightarrow \mathbb{P}^{k+2} (K)$. This operator does not have approximation properties if $K\in\mesh\upb$ because some boundary terms have been removed in the definition of the reconstruction operator. This leads us to define the lifting operator $\mathcal{L}_K:  H^2(K) \rightarrow  \mathbb{P}^{k+2}(K)$ for all $K\in \mesh$ such that, for all $v \in H^2(K)$ and all $w \in \mathbb{P}^{k+2}(K)^\perp$,
\begin{equation}\label{def: lifting}
(\nabla \mathcal{L}_K(v), \nabla w)_{K,\varepsilon}
:= -\varepsilon\Big((v,  \partial_{n} \Delta w))_{\dKb}
-(\partial_n v, \partial_{nn}  w )_{\dKb}-(\partial_t v,\partial_{nt}  w )_{\dKb}\Big)
+ (v, \partial_{n} w)_{\dKb},
\end{equation}
together with the condition $(\mathcal{L}_K(v),1)_K =  0$.
Notice that $\mathcal{L}_K(v) = 0$ for all $K\in \meshi$.
We then define the operator $\mathcal{E}_K: H^2(K) \rightarrow \mathbb{P}^{k+2} (K)$ such that
\begin{equation}
\mathcal{E}_K(v): =   \mathcal{E}_K\upi(v) + \mathcal{L}_K(v).
\end{equation}
The definition of $R_K\upi$ implies that
\begin{align*}
(\nabla R_K\upi (\mathcal{\hat{I}}^k_K (v)), \nabla w)_{K,\varepsilon}
={}& 	\varepsilon\Big(( \Pi_{K}^{k+2}(v) , \Delta^2 w)_K
- ( \Pi_{\dKi}^{k+2} (v) , \partial_n \Delta  w)_{\dKi}
+ ( \Pi^k_{\dKi} (\partial_n v),  \partial_{nn}  w)_{\dKi} \\
&
+  (\partial_t( \Pi^{k+2}_{\dKi} (v)),  \partial_{nt}  w)_{\dKi}\Big) - ( \Pi_{K}^{k+2}(v) , \Delta w)_K + {(\Pi^{k+2}_{\dKi}(v) , \partial_n w)_{\dKi}}.
\end{align*}
Since $w \in  \mathbb{P}^{k+2}(K)$, we infer that
\begin{align*}
(\nabla R_K\upi (\mathcal{\hat{I}}^k_K (v)), \nabla w)_{K,\varepsilon}
={}& 	\varepsilon\Big(( v , \Delta^2 w)_K - ( v , \partial_n \Delta  w)_{\dKi}
+ ( \partial_n v,  \partial_{nn}  w)_{\dKi} +(\partial_t v,  \partial_{nt}  w)_{\dKi}\\
& + (\partial_t( \Pi^{k+2}_{\dKi} (v)-v),  \partial_{nt}  w)_{\dKi}\Big) - (v , \Delta w)_K + {(v , \partial_n w)_{\dKi}}.
\end{align*}
Using the definitions of $\mathcal{L}_K$ and $\mathcal{E}_K$ shows that
\begin{align}
(\nabla \mathcal{E}_K(v), \nabla w)_{K,\varepsilon}
={}& 	\varepsilon\Big(( v , \Delta^2 w)_K - ( v , \partial_n \Delta  w)_{\dK}
+ ( \partial_n v,  \partial_{nn}  w)_{\dK} +(\partial_t v,  \partial_{nt}  w)_{\dK} \nonumber \\
& +  (\partial_t( \Pi^{k+2}_{\dKi} (v)-v),  \partial_{nt}  w)_{\dKi}\Big) - (v , \Delta w)_K + {(v , \partial_n   w)_{\dK}}. \label{eq:id_EK}
\end{align}
Integration by parts then implies that
\begin{equation*}
( \nabla  (\mathcal{E}_K (v) - v), \nabla w)_{K,\varepsilon}
= \varepsilon (\partial_t( \Pi^{k+2}_{\dKi} (v)-v),  \partial_{nt}  w)_{\dKi}.
\end{equation*}
This shows that the operator $\mathcal{E}_K$ is a projection and that it coincides with the $H^1$-elliptic projection if $\varepsilon=0$. The following result establishes the approximation properties of the projection operator $\mathcal{E}_K$ in the general case, as well as the approximation properties for the stabilization operators. To state the result, we consider the following norm for all $K\in\mesh$ and all $v\in H^{2+s}(K)$, $s>\frac{3}{2}$:
\begin{equation}\label{Def: spectical norm}
\begin{aligned}
\|v\|^2_{\sharp, K}&:= \|\nabla v\|_{K,\varepsilon}^2 +
\varepsilon \Big(h^{3}_K \|\partial_{n} \Delta v \|_{\partial K}^2
+h_K \|\partial_{nn}  v \|_{\partial K}^2
+h_K \|\partial_{nt}  v \|_{\partial K}^2 \Big) +h^{}_K \|\partial_{n} v \|_{\partial K}^2.
\end{aligned}
\end{equation}
The proof of the following lemma is postponed to Section~\ref{sec:proof_approximation}.

\begin{lemma}[Approximation] \label{lem: approxmation}
The following holds for all $K\in \mesh$ and all $v\in H^{2+s}(K)$, $s>\frac32$:
\begin{equation}\label{stabilisation bound}
\|v - \mathcal{E}_K (v)\|_{\sharp, K}^2 + S\upi_{\dK}(\mathcal{\hat{I}}^k_K(v),\mathcal{\hat{I}}^k_K(v)) +  S\upb_{\dK}(v - \Pi_K^{k+2}(v),v - \Pi_K^{k+2}(v))
\leq
C \| v - \Pi_K^{k+2}(v)\|_{\sharp,K}^2.
\end{equation}
\end{lemma}

The global interpolation operator $\mathcal{\hat{I}}_h^k:H^2(\Omega)\to \fes$ is defined such that, for all $v\in H^2(\Omega)$,
\begin{equation} \label{def:Ihk_bis}
\Ihk(v):=\big( (\Pi_K^{k+2}(v))_{K\in\mesh},(\Pi^{k+2}_F(v))_{F\in\Fint},(\Pi_F^{k}(\n_F{\cdot}\nabla v))_{F\in\Fint} \big)\in \fes,
\end{equation}
so that the local components of $\Ihk(v)$ are $\mathcal{\hat{I}}^k_K(v|_K)$ for all $K\in\mesh$. We define the consistency error $\delta_h\in (\fes)^\prime$ such that, for all $\hat{w}_h\in \fes$,
\begin{equation}
	\langle \delta_h,\hat{w}_h \rangle:= {\ell}(\hat{w}_h) - \hat{a}_h(\Ihk(u), \hat{w}_h),
\end{equation}
where the brackets refer to the duality pairing between $(\fes)^\prime$ and $\fes$.
The proof of the following bound on the consistency error is postponed to Section~\ref{sec:proof_consistency}.

\begin{lemma}[Consistency]\label{lemma: consistency}
Assume that $u\in H^{2+s}(\Omega)$, $s>\frac{3}{2}$. The following holds:
\begin{equation} \label{consistency}
\langle \delta_h,\hat{w}_h \rangle
\leq C
\left( \su \|u-\Pi_K^{k+2} (u)\|^2_{\sharp,K}\right)^{\frac12} \|\hat{w}_h\|_{\fes},
\qquad \forall \hat{w}_h\in \fes.
\end{equation}
\end{lemma}

\subsection{Error estimate}

The above results lead to the following error bound. The proof is postponed to Section~\ref{sec:proof_main}. {Let $\hat u_h\in\fes$ be the discrete HHO solution and recall from Section~\ref{sec:local} that for all $K\in\mesh$, $\hat{u}_K\in\fesE$ denotes the
local components of $\hat{u}_h$
associated with the mesh cell $K$ and its faces in $\FK\upi$.}
To simplify the notation, we set $R_K(\hat{u}_K):=R^i_K(\hat{u}_K)+\mathcal{L}_K(u)$ for all $K\in\mesh$, and recall that $\mathcal{L}_K(u)$ is nonzero only on boundary cells where it is fully computable from the boundary data.

\begin{theorem}[Error estimate]\label{Theorem: main}
Assume that $u\in H^{2+s}(\Omega)$ with $s>\frac{3}{2}$. The following holds true:
\begin{equation}\label{abstract error bound}
\su \|\nabla (u- R_K(\hat{u}_K))\|_{K,\varepsilon}^2 \leq C \su \| u-\Pi_K^{k+2}(u)\|_{\sharp,K}^2.
\end{equation}
Consequently, if $k\ge1$, assuming $u|_K\in H^{k+3}(K)$ for all $K\in\mesh$, we have
\begin{equation}\label{eq:err2}
\su   \|\nabla (u- R_K(\hat{u}_K))\|_{K,\varepsilon}^2 \leq C \su \big( \sigma^{\frac12}_K h_K^{k+2}|u|_{H^{k+3}(K)}\big)^2,
\end{equation}
and if $k=0$, assuming $u|_K\in H^{4}(K)$ for all $K\in\mesh$, we have
\begin{equation} \label{eq:err3}
\su \|\nabla (u- R_K(\hat{u}_K))\|_{K,\varepsilon}^2 \leq C \su \big(  \sigma^{\frac12}_K h_K^2(|u|_{H^{3}(K)}+h_K|u|_{H^{4}(K)})\big)^2.
\end{equation}
\end{theorem}

\begin{remark}[Error estimate \eqref{eq:err2}] \label{rem:err2}
In the case where $\varepsilon \approx \ell_\Omega^2$, i.e., the fourth-order operator is dominant, we have $\sigma_K \approx \mathcal{O}(h_K^{-2})$, so that the error estimate \eqref{eq:err2} implies that
$$
\su \|\nabla^2 (u- R_K(\hat{u}_K))\|_{K}^2
\leq C 	\su \big(  h_K^{k+1}|u|_{H^{k+3}(K)}\big)^2,
$$
which corresponds to the error estimate obtained in \cite{DongErn2021biharmonic}
for the biharmonic problem.
Instead, in the case where $\varepsilon\ll1$, one has in practice $\sigma_K=1$ (unless extremely fine meshes are used), and the error estimate \eqref{eq:err2} implies that
$$
\su \|\nabla (u-R_K(\hat{u}_K))\|_{K}^2
\leq C 	\su \big(  h_K^{k+2}|u|_{H^{k+3}(K)}\big)^2.
$$
Similar comments can be made for \eqref{eq:err3}.
\end{remark}

\begin{remark}[Regularity assumption]
The present error analysis requires that the exact solution has the minimal regularity
$u\in H^{2+s}(\Omega)$, with $s>\frac32$. This assumption
is consistent with the rather classical paradigm encountered in the literature when analyzing nonconforming approximation methods. Notice that this regularity requirement is less stringent
than the one needed to achieve optimal decay rates as soon as $k\ge1$. Moreover,
this requirement
can be lowered to $s>1$ by using the techniques developed in \cite{ErnGuer2021}
and \cite[Chap.~40-41]{Ern_Guermond_FEs_II_2021} in the context of second-order
elliptic PDEs.
\end{remark}

\begin{remark}[$k=0$]
The regularity assumption $u|_{K}\in H^4(K)$ on the exact solution is slightly suboptimal in the case $k=0$ whenever $s<1$. This assumption can be avoided if a multiplicative
trace inequality in fractional Sobolev spaces is available on cells with a curved boundary.
Specifically, we need to assert that for $s\in(\frac12,1)$,
there is $C$ such that for all $h>0$,
all $K\in\mesh$, and all $v\in H^1(K)$,
we have $\|{v} \|_{\dK} \leq C \big(
h_K^{-\frac12}  \|{v}\|_{K} + h_K^{s-\frac12}  |v|_{H^s(K)}\big)$.
This inequality can be established on cells with a flat boundary by invoking affine
geometric mappings, see \cite[Lem.~7.2]{ErnGuermond:17}.
\end{remark}

\section{Numerical examples} \label{sec: numerical examples}

In this section, we present numerical examples to illustrate the theoretical results on the present HHO method. We first study convergence rates and robustness for smooth solutions
in domains with a polygonal (Section~\ref{sec:res_poly}) and a curved (Section~\ref{sec:res_curved}) boundary. Then we consider a more challenging test case with unknown analytical solution and a boundary layer forming as $\varepsilon\to0$. All the computations were
run with \texttt{Matlab R2018a} on the NEF
platform at INRIA Sophia Antipolis M\'editerran\'ee using 12 cores,
and all the linear systems after static condensation are solved using the
\verb*|backslash| function. {The quadratures in polygonal cells are performed by sub-triangulating the polygon into triangles. For every curved element, the sub-triangulation is constructed by considering a sufficiently fine decomposition of its curved edge into smaller straight sub-edges. In our implementation, we consider $30$ sub-edges; this number was verified to be sufficient on the finest meshes and highest polynomial degrees reported in Section~\ref{sec:res_poly}. We emphasize that the sub-triangulation is only used to generate the quadrature rules and that these calculations are fully parallelizable.}

\subsection{Convergence rates and robustness in a polygonal domain}
\label{sec:res_poly}

\begin{figure}[t]
    \centering
    \includegraphics[scale=0.33]{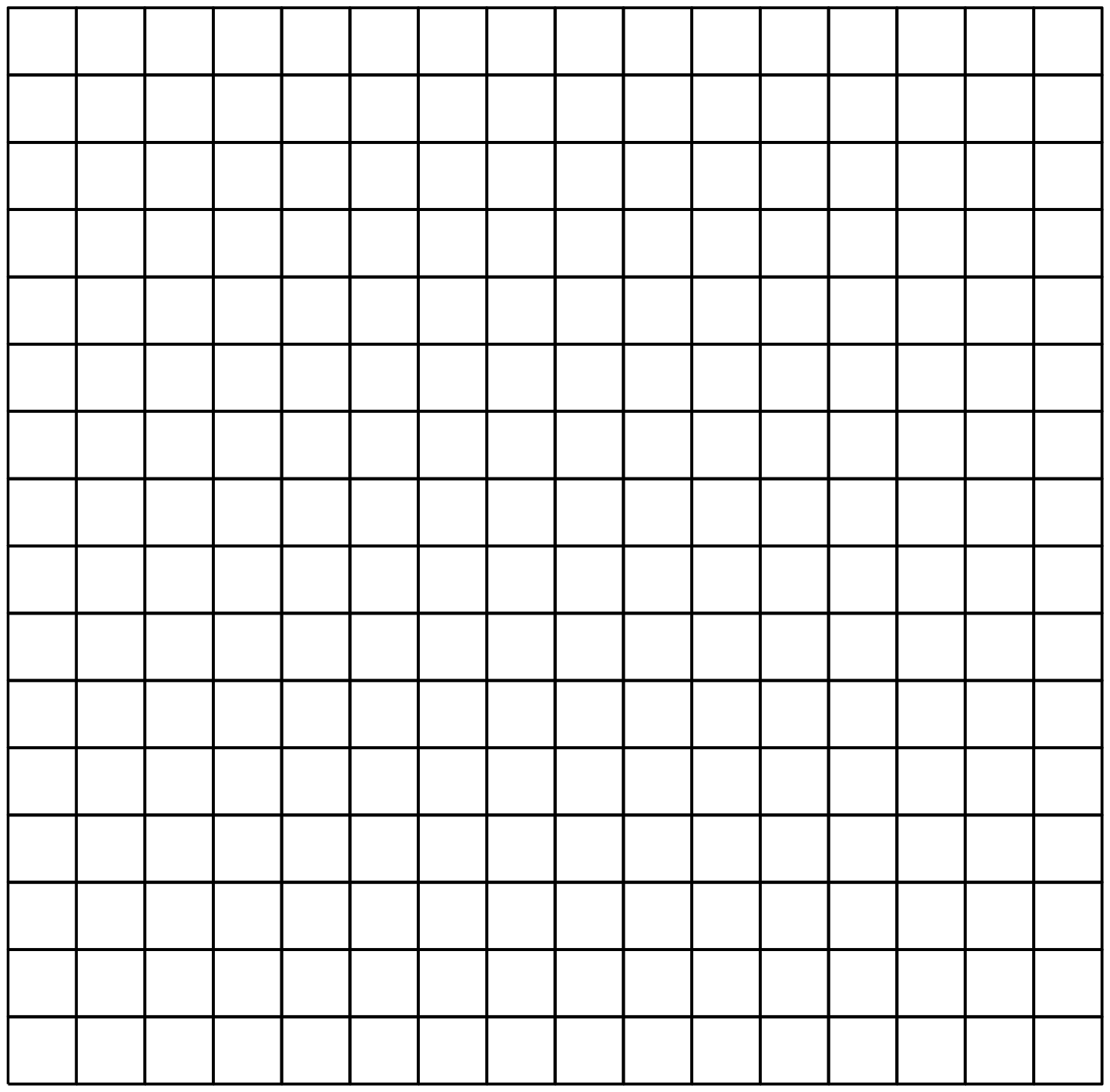}
    \hspace{1cm}
    \includegraphics[scale=0.33]{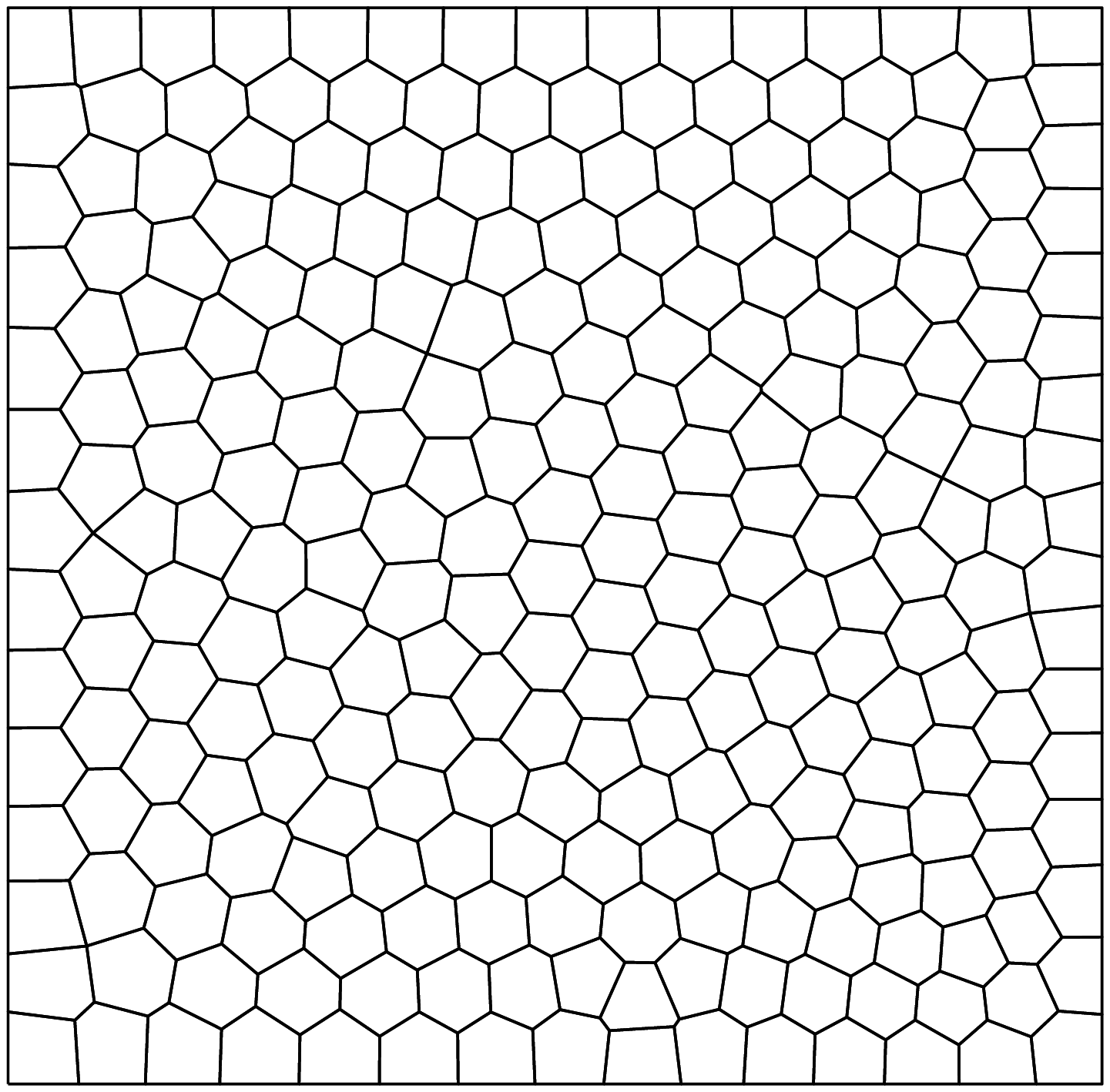}
    \caption{Examples of rectangular (left) and polygonal (Voronoi-like) (right) meshes composed of 256 cells.}\label{ex1:mesh_figure}
\end{figure}
\begin{figure}[!htb]
    \begin{center}
        \begin{tabular}{cc}
            \includegraphics[scale=0.4]{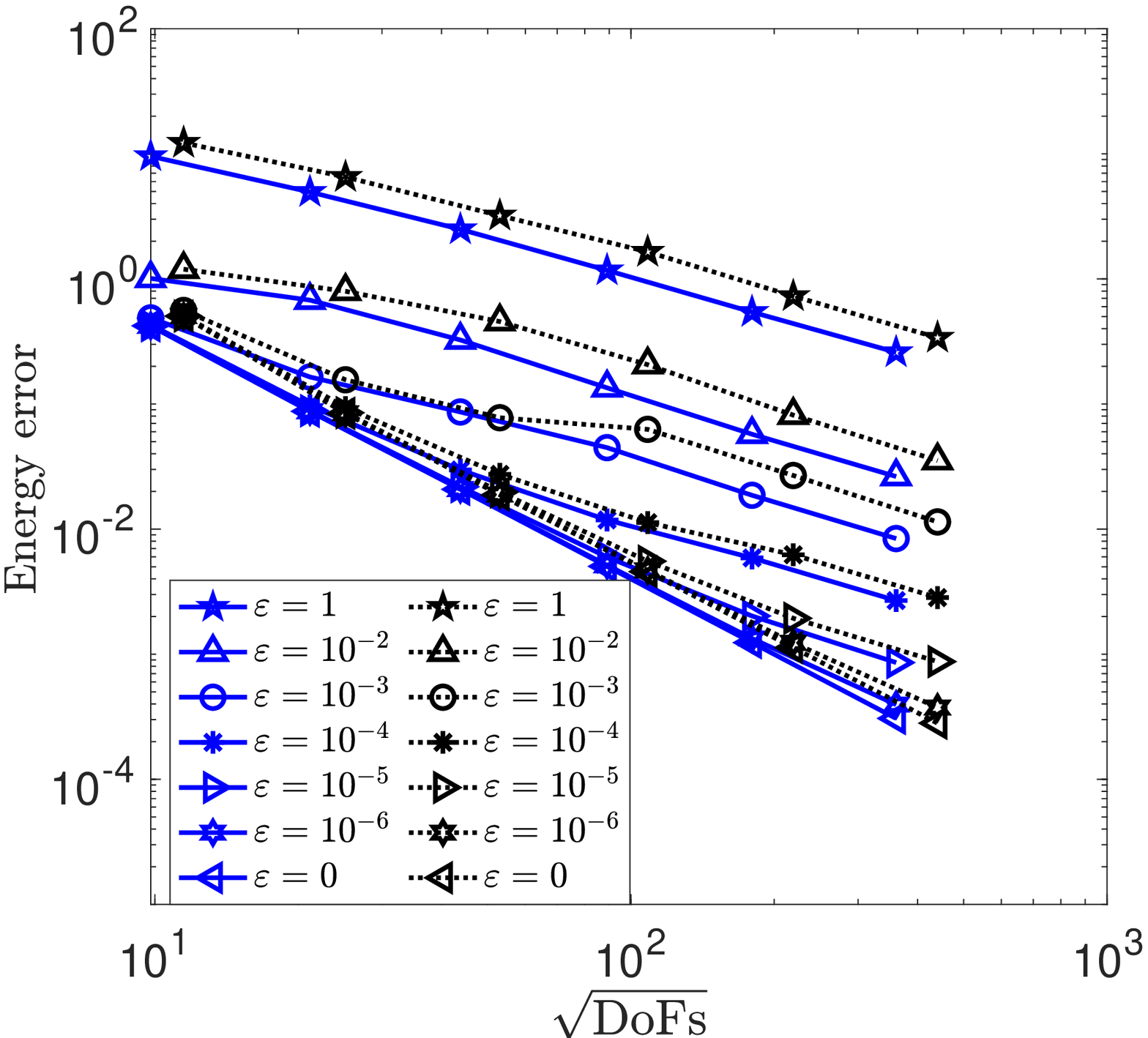} &
            \includegraphics[scale=0.4]{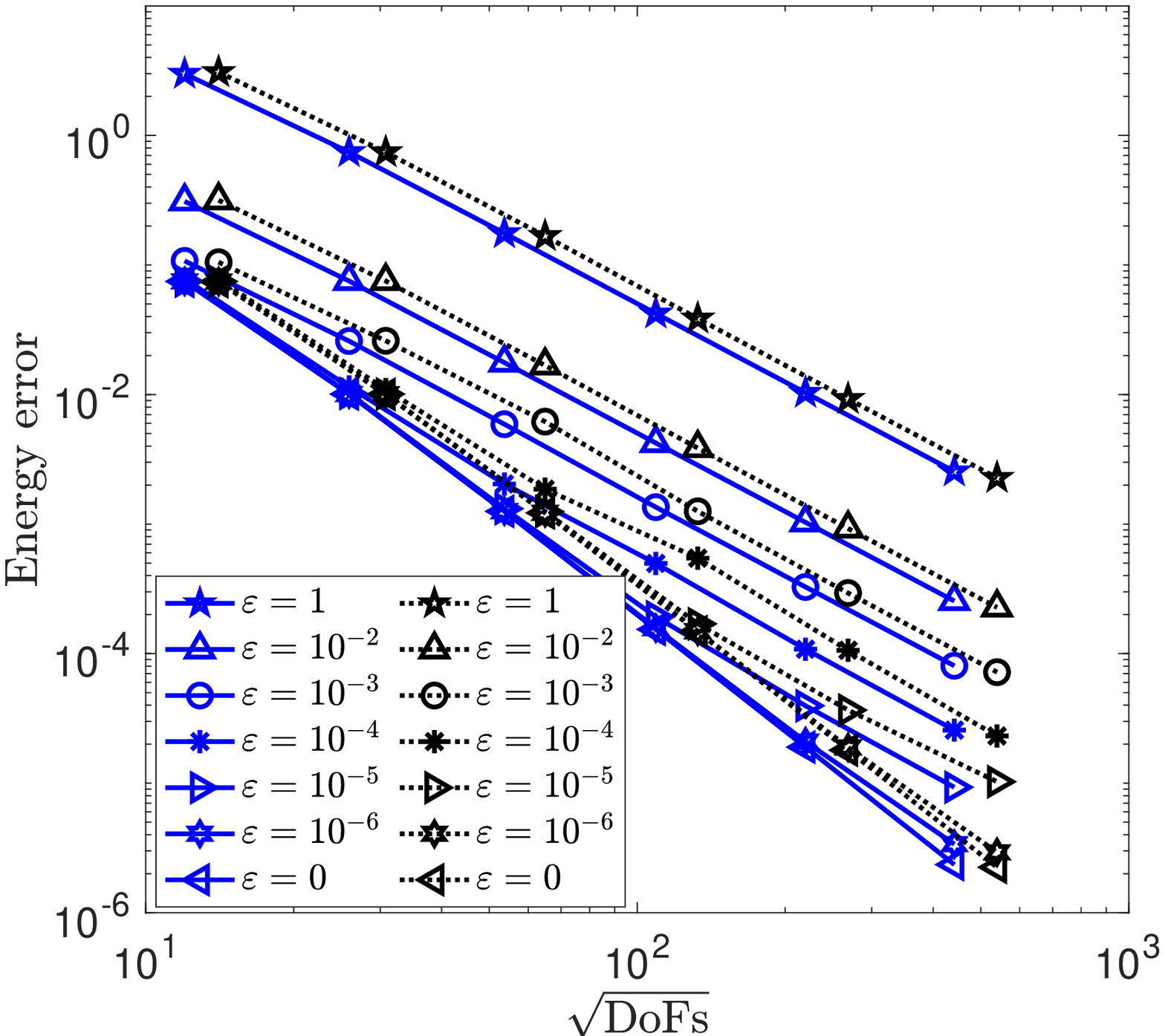} \\
            \includegraphics[scale=0.4]{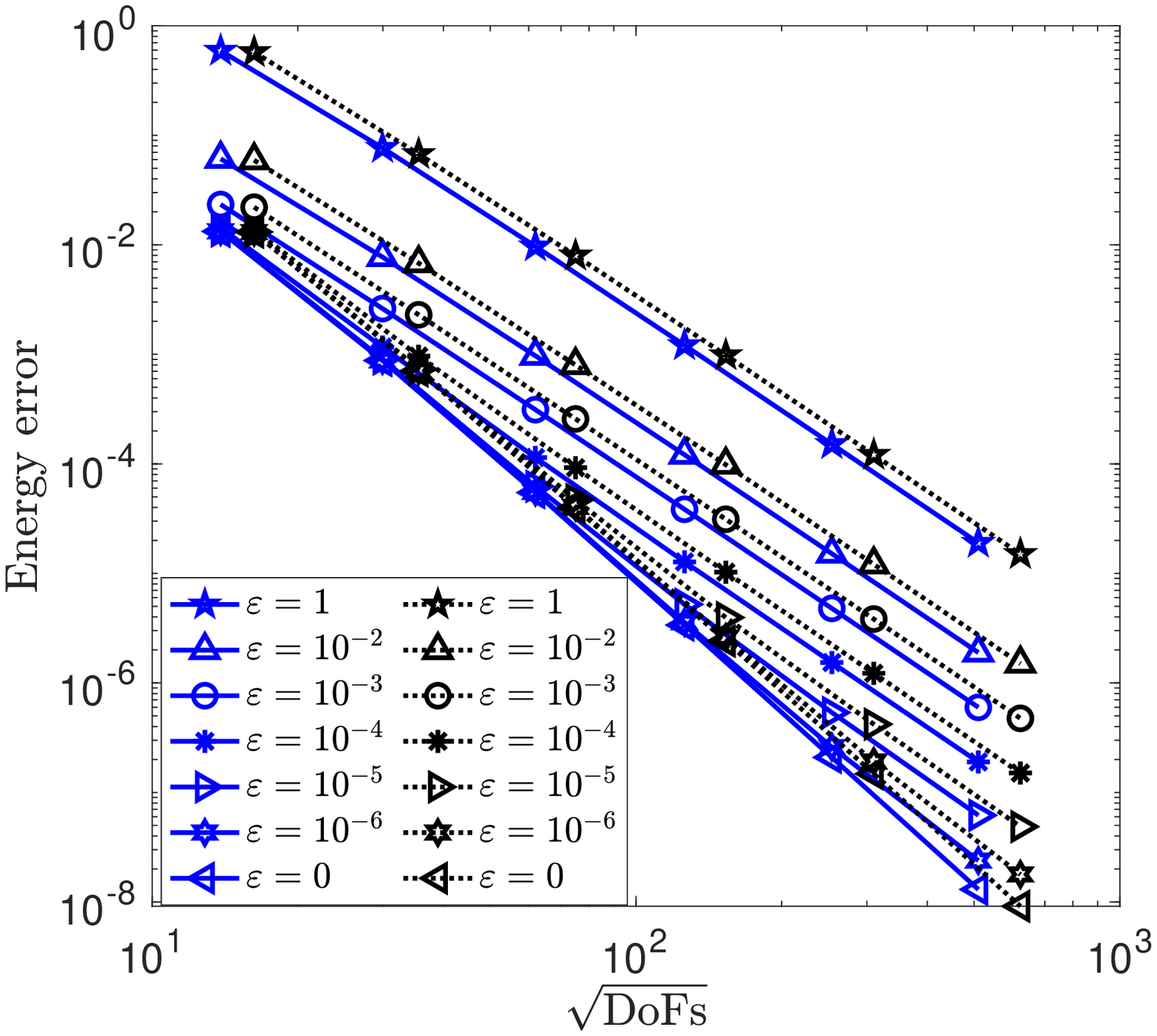} &
            \includegraphics[scale=0.4]{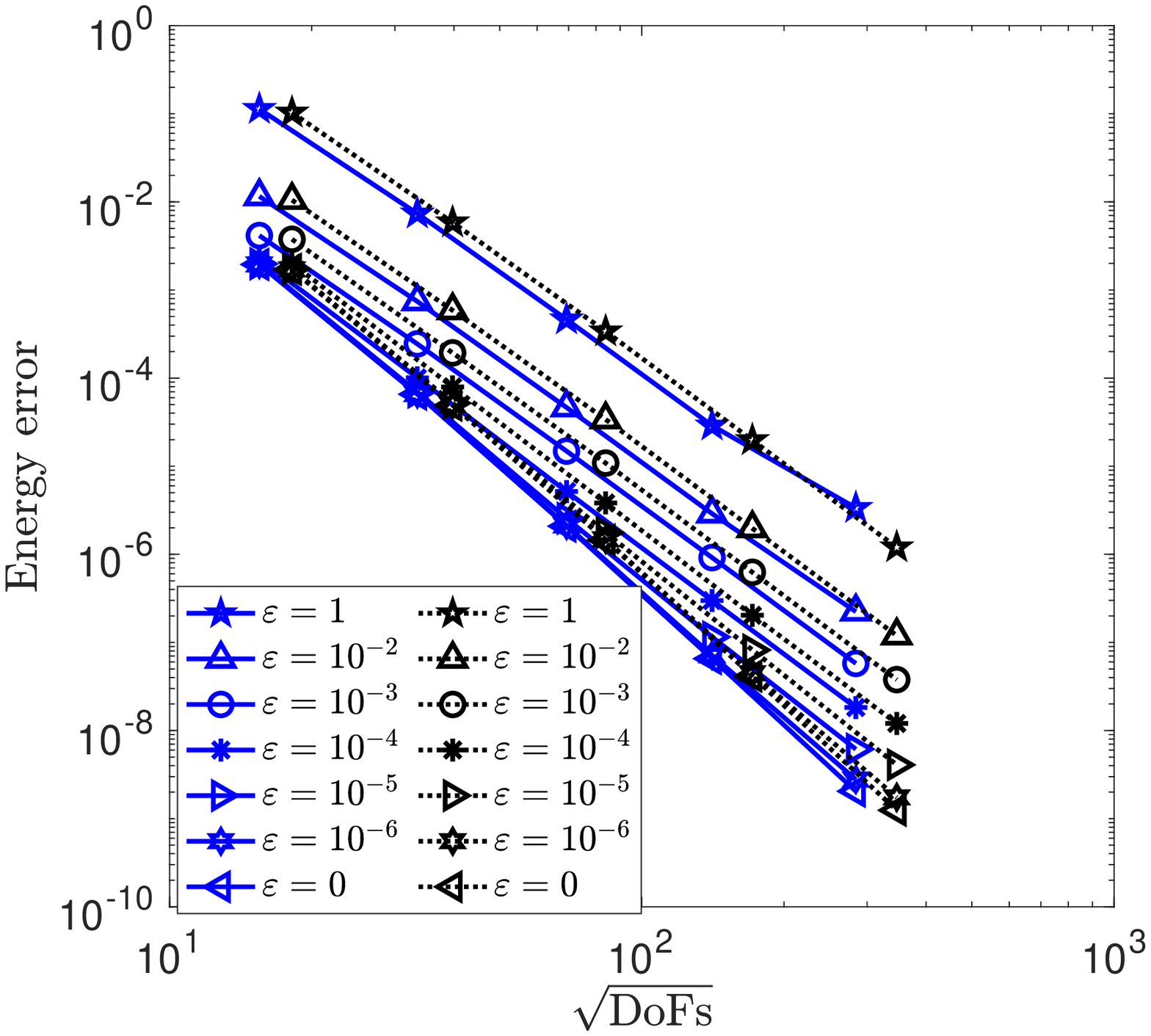} \\
        \end{tabular}
    \end{center}
    \caption{Polygonal domain and smooth solution. Convergence rates in the energy seminorm on rectangular (solid line) and polygonal (dotted line) meshes for different values of $\varepsilon$. The polynomial degree is $k=0$ (top left), $k=1$  (top right), $k=2$  (bottom left) and $k=3$  (bottom right).}\label{ex1:h-refine}
\end{figure}

We select $f$ and the boundary conditions on $\Omega:=(0,1)^2$ so that the exact solution to \eqref{pde} is
$$
u(x,y) = \sin(\pi x)^2  \sin(\pi y)^2+e^{-(x-0.5)^2-(y-0.5)^2}.
$$
We employ polynomial degrees $k\in\{0,\ldots,3\}$ and meshes consisting of $\{16,64,256,1024,4096,16384\}$ elements. The meshes can be either composed of rectangular cells or of
polygonal (Voronoi-like) cells (generated through the PolyMesher \texttt{Matlab}
library \cite{polymesher}). Two examples of rectangular and polygonal meshes, both composed of 256 cells, are shown in Figure \ref{ex1:mesh_figure}.
Despite an $hp$-error analysis falls beyond
the present scope, we implement the stabilization terms in \eqref{def: stabilisation}
and \eqref{def: stabilisation Boundary} with $h_K^{-1}$ replaced by $(k+1)^2 h_K^{-1}$
for all $K\in\mesh$.

 \begin{table}[!htb]
    \begin{center}
        \begin{tabular}{||c|c|c|c|c|c|c|c||}
            \hline
            $\#$ \text{Cells}  & $\varepsilon=1$ & $\varepsilon=10^{-2}$ & $\varepsilon=10^{-3}$ & $\varepsilon=10^{-4}$ & $\varepsilon=10^{-5}$ & $\varepsilon=10^{-6}$ & $\varepsilon=0$\\
            \hline   \hline
            &
            &
            \multicolumn{4}{c}{ Rectangular meshes with $k=0$ }
            &
            &
            \\
            \hline
            64 & 0.85 & 0.52 & 1.41 & 1.93 & 2.03 & 2.04 & 2.04\\
            \hline
            256 & 0.94 & 1.01 & 0.89 & 1.62 & 1.92 & 1.97 & 1.97\\
            \hline
            1024 & 1.06 & 1.24 & 0.93 & 1.30 & 1.82 & 1.98 & 2.00\\
            \hline
            4096 & 1.10 & 1.22 & 1.25 & 1.00 & 1.55 & 1.92 & 2.00\\
            \hline
            16384 & 1.06 & 1.12 & 1.14 & 1.13 & 1.23 & 1.75 & 2.00\\
            \hline
            \hline
            &
            &
            \multicolumn{4}{c}{ Rectangular meshes with $k=1$ }
            &
            &\\
            \hline
            64 & 1.81 & 1.84 & 1.85 & 2.49 & 2.59 & 2.60 & 2.60\\
            \hline
            256 & 1.98 & 1.99 & 2.03 & 2.37 & 2.81 & 2.86 & 2.86\\
            \hline
            1024 & 2.01 & 2.01 & 2.08 & 1.99 & 2.71 & 2.93 & 2.96\\
            \hline
            4096 & 2.01 & 2.01 & 2.02 & 2.17 & 2.26 & 2.88 & 2.99\\
            \hline
            16384 & 2.01 & 2.01 & 2.01 & 2.06 & 2.08 & 2.61 & 2.99\\
            \hline
            \hline
            &
            &
            \multicolumn{4}{c}{ Rectangular meshes with $k=2$ }
            &
            &\\
            \hline
            64 & 2.65 & 2.68 & 2.84 & 3.26 & 3.49 & 3.52 & 3.53\\
            \hline
            256 & 2.85 & 2.85 & 2.92 & 3.19 & 3.68 & 3.80 & 3.82\\
            \hline
            1024 & 2.93 & 2.93 & 2.95 & 3.09 & 3.50 & 3.86 & 3.93\\
            \hline
            4096 & 2.97 & 2.97 & 2.97 & 3.02 & 3.24 & 3.74 & 3.97\\
            \hline
            16384 & 2.98 & 2.98 & 2.98 & 3.00 & 3.10 & 3.42 & 3.98\\
            \hline
            \hline
            &
            &
            \multicolumn{4}{c}{ Rectangular meshes with $k=3$ }
            &
            &\\
            \hline
            64 & 3.55 & 3.56 & 3.68 & 4.04 & 4.35 & 4.40 & 4.40\\
            \hline
            256 & 3.80 & 3.81 & 3.85 & 4.08 & 4.54 & 4.72 & 4.74\\
            \hline
            1024 & 3.91 & 3.91 & 3.92 & 4.02 & 4.37 & 4.79 & 4.88\\
            \hline
            4096 & 3.97 & 3.95 & 3.94 & 3.99 & 4.18 & 4.63 & 4.94\\
            \hline
        \end{tabular}
    \end{center}
    \caption{Polygonal domain and smooth solution. Convergence rates in the energy seminorm on rectangular meshes for different values of $k$ and $\varepsilon$.} \label{table:h-convergence}
\end{table}

\begin{table}[!htb]
    \begin{center}
        \begin{tabular}{||c|c|c|c|c|c|c|c||}
            \hline
            $\#$ \text{Cells}  & $\varepsilon=1$ & $\varepsilon=10^{-2}$ & $\varepsilon=10^{-3}$ & $\varepsilon=10^{-4}$ & $\varepsilon=10^{-5}$ & $\varepsilon=10^{-6}$ & $\varepsilon=0$\\
            \hline   \hline
            &
            &
            \multicolumn{4}{c}{ Rectangular meshes with $k=0$ }
            &
            &
            \\
            \hline
            64 & 1.40 & 0.17 & 2.36 & 2.92 & 2.94 & 2.94 & 2.94\\
            \hline
            256 & 1.76 & 1.36 & 0.52 & 2.81 & 2.95 & 2.96 & 2.96\\
            \hline
            1024 & 1.83 & 1.66 & 0.81 & 2.36 & 3.01 & 3.03 & 3.03\\
            \hline
            4096 & 1.90 & 1.84 & 1.75 & 0.74 & 2.97 & 3.02 & 3.03\\
            \hline
            16384 & 1.98 & 1.97 & 1.94 & 1.52 & 2.57 & 3.01 & 3.01\\
            \hline
            \hline
            &
            &
            \multicolumn{4}{c}{ Rectangular meshes with $k=1$ }
            &
            &\\
            \hline
            64 & 3.35 & 3.24 & 1.71 & 3.49 & 3.46 & 3.45 & 3.45\\
            \hline
            256 & 3.68 & 3.80 & 2.96 & 3.10 & 3.77 & 3.75 & 3.74\\
            \hline
            1024 & 3.84 & 3.94 & 3.78 & 1.98 & 3.90 & 3.89 & 3.87\\
            \hline
            4096 & 3.92 & 3.97 & 3.96 & 3.42 & 2.54 & 3.98 & 3.94\\
            \hline
            16384 & 4.05 & 4.02 & 4.01 & 3.93 & 2.69 & 3.92 & 3.97\\
            \hline
            \hline
            &
            &
            \multicolumn{4}{c}{ Rectangular meshes with $k=2$ }
            &
            &\\
            \hline
            64 & 4.51 & 4.60 & 4.32 & 4.42& 4.37 & 4.36 & 4.36\\
            \hline
            256 & 4.77 & 4.81 &4.72 & 4.57 & 4.74 & 4.70 & 4.69\\
            \hline
            1024 & 4.98 & 4.93 & 4.82 & 4.57 & 4.99 & 4.88 & 4.85\\
            \hline
            4096 & 2.49 & 4.17 & 4.92 & 4.83 & 4.69 & 5.02& 4.97\\
            \hline
            \hline
            &
            &
            \multicolumn{4}{c}{ Rectangular meshes with $k=3$ }
            &
            &\\
            \hline
            64 & 5.29 & 5.26 & 5.05 & 5.12 & 5.21 & 5.20 & 5.19\\
            \hline
            256 & 5.87 & 5.87 & 5.59 & 5.37 & 5.66 & 5.62 & 5.61\\
            \hline
            1024 & 5.06 & 5.28 & 5.89 & 5.87 & 5.86 & 5.85 & 5.88\\
            \hline
        \end{tabular}
    \end{center}
    \caption{{Polygonal domain and smooth solution. Convergence rates in the $L^2$-norm on rectangular meshes for different values of $k$ and $\varepsilon$.}} \label{table:h-convergence L2 norm}
\end{table}

\begin{table}[!htb]
    \begin{center}
        \begin{tabular}{||c|c|c|c|c|c|c||}
            \hline
            $\#$ \text{Cells}  & $\varepsilon=1$ &  $\varepsilon=10^{-4}$ & $\varepsilon=10^{-5}$ & $\varepsilon=10^{-6}$ & $\varepsilon=0$ & 2nd-order  \\
            \hline   \hline
            &
            &
            \multicolumn{3}{c}{ Rectangular meshes with $k=0$ }
            &
            &
            \\
            \hline
            1024 & 2.10e+06 & 8.44e+04 & 2.49e+05 & 3.13e+05 & 3.18e+05 & 2.47e+03\\
            \hline
            4096 & 3.38e+07 & 5.99e+05 & 6.15e+05 & 1.20e+06 & 1.33e+06 & 9.72e+03\\
            \hline
            16384 & 5.42e+08 & 6.84e+06 & 1.64e+06 & 3.84e+06 & 4.91e+06 & 3.85e+04\\
            \hline
            \hline
            &
            &
            \multicolumn{3}{c}{ Rectangular meshes with $k=1$ }
            &
            &\\
            \hline
            1024 & 2.52e+07 & 3.72e+05 & 4.34e+05 & 7.98e+05 & 8.82e+05 & 4.89e+03\\
            \hline
            4096 & 3.97e+08 & 4.86e+06 & 1.20e+06 & 2.44e+06 & 3.63e+06 & 1.92e+04\\
            \hline
            16384 & 6.29e+09 & 7.76e+07 & 9.84e+06 & 6.25e+06 & 1.49e+07 & 7.60e+04\\
            \hline
            \hline
            &
            &
            \multicolumn{3}{c}{ Rectangular meshes with $k=2$ }
            &
            &\\
            \hline
            1024 & 1.45e+08 & 2.04e+06 & 6.45e+05 & 1.75e+06 & 2.27e+06 & 9.29e+03\\
            \hline
            4096 & 2.28e+09 & 2.86e+07 & 1.20e+06 & 2.44e+06 & 9.30e+06 & 3.64e+04\\
            \hline
            16384 & 3.62e+10 & 4.52e+08 & 5.96e+07 & 1.02e+07 & 3.77e+07 & 1.44e+05\\
            \hline
            \hline
            &
            &
            \multicolumn{3}{c}{ Rectangular meshes with $k=3$ }
            &
            &\\
            \hline
            1024 &  4.77e+08 & 7.37e+06 & 1.40e+06 & 2.81e+06 & 4.88e+06 & 1.28e+04\\
              \hline
            4096 &  7.51e+09 & 9.59e+07 & 1.61e+07 & 5.77e+06 & 1.98e+07 & 5.04e+04\\
            \hline
            16384 & 1.19e+11 & 1.49e+09 & 2.02e+08 & 3.17e+07 & 8.01e+07 & 1.99e+05\\
            \hline
            \hline
        \end{tabular}
    \end{center}
    \caption{{The condition number of condensed linear system on rectangular meshes for different values of $k$ and $\varepsilon$ for the proposed HHO method and the mixed order HHO method.}} \label{table:condition number}
\end{table}

Let us first verify the convergence rates. We measure relative errors in the (broken) energy seminorm used on the left-hand side of our main error estimate \eqref{abstract error bound}. The errors are reported as a function of
$\sqrt{\mathrm{DoFs}}$, where $\mathrm{DoFs}$ denotes the total number of
globally coupled
discrete unknowns (that is, the face unknowns). The results are reported
in Figure \ref{ex1:h-refine}.
The first observation is that there is almost no difference in the convergence rates obtained on rectangular and polygonal meshes, and that these rates match the prediction of Theorem~\ref{Theorem: main} for all the polynomial degrees. Next, we observe that the errors obtained with $\varepsilon =1$ and $\varepsilon =0$  converge, respectively, at the optimal rates $O(h^{k+1})$ and $O(h^{k+2})$, in agreement with Remark~\ref{rem:err2}. For values of
$\varepsilon$ between these two extreme values, a transition between the above
two regimes is observed for $\varepsilon\approx h^2$. Whenever $\varepsilon \ge h^2$, the convergence rate is $O(h^{k+1})$ as expected for a fourth-order differential operator. Instead, whenever $\varepsilon \le h^2$, the convergence rate is closer to the value $O(h^{k+2})$ expected for a second-order differential operator.
In Table \ref{table:h-convergence}, we list the convergence rates on rectangular meshes (those on polygonal meshes lead to the same conclusions).  Reading the table horizontally, the transition between the two regimes is clearly visible.
{For completeness, we also report the relative errors in the $L^2$-norm in Table \ref{table:h-convergence L2 norm}. As for the energy seminorm error, there is almost no difference in the convergence rates obtained on rectangular and polygonal meshes, so that we focus on rectangular meshes. For $\varepsilon=0$, the convergence rate is
always $O(h^{k+3})$ which is optimal (recall that polynomials of order $(k+2)$ are employed to approximate the traces on the mesh faces). For $\varepsilon=1$, the convergence rate is suboptimal for $k=0$, i.e., only $O(h^2)$, as also observed with other nonconforming finite element methods applied to fourth-order PDEs. Instead, the optimal rate $O(h^{k+3})$ is recovered for $k\ge1$. For $k=0$, a transition between second- and third-order convergence is observed as $\varepsilon\to0$.}

{Finally, we present in Table \ref{table:condition number} the (Euclidean) condition number of the linear system after static condensation for all the considered values of $\varepsilon$. To compare the value obtained for $\varepsilon=0$, we also report the condition number for the linear system discretized by a genuine HHO method for the second-order PDE, employing polynomials of order $(k+2)$ for the cell unknowns and $(k+1)$ for the face unknowns. The first observation is that the condition number for $\varepsilon=1$ scales as $\mathcal{O}(h^{-4})$ (as expected) and is larger than the condition number for all the other values of $\varepsilon$ by two orders of magnitude for all $k\ge 0$. Instead, for $\varepsilon=0$, the condition number scales as $\mathcal{O}(h^{-2})$ (again, as expected). In addition, the condition number for the genuinely second-order HHO method has the same quadratic scaling for the condition number, with values that are two orders of magnitude smaller than those reported in the column $\varepsilon=0$. This is reasonable since the proposed HHO method is not designed for genuinely second-order operators, but it still gives the optimal quadratic scaling for the condition number in the limit case $\varepsilon=0$. This encouraging observation indicates that the linear systems obtained with the present HHO method remain relatively well-behaved as $\varepsilon\to0$. Further studies are, however, needed, including, e.g., preconditioned iterative methods (as, for instance, the one conducted in \cite{huang2021morley}). Finally, for decreasing $\varepsilon$ between the two extreme values $1$ and $0$, the condition number for fixed $k$ and $h$ decreases first and then increases. This nonmonotone behavior is probably related to the two different asymptotic regimes associated with $\varepsilon=1$ and $\varepsilon=0$.}

\begin{figure}[htb!]
    \centering
    \includegraphics[scale=0.3]{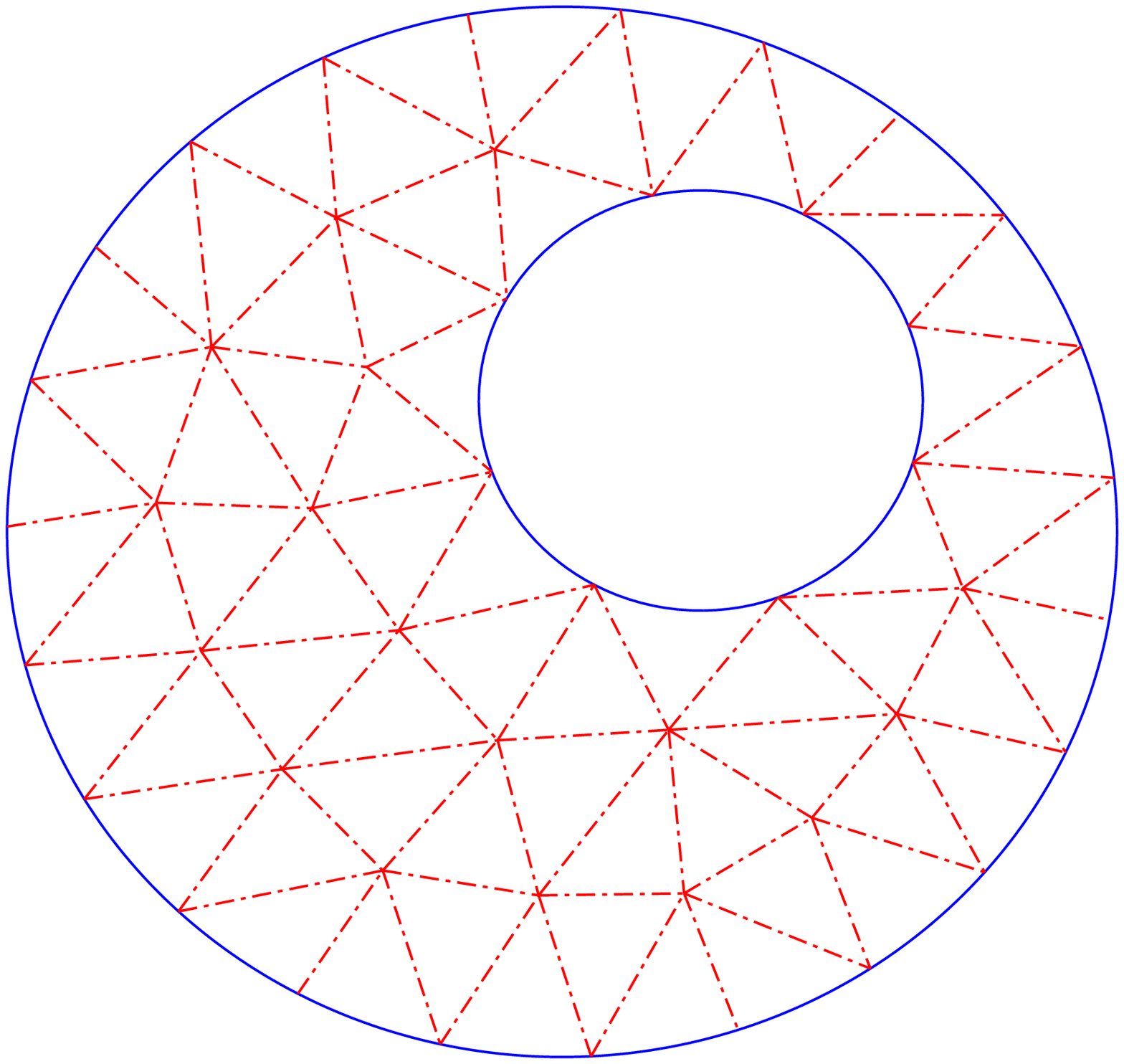}
    \hspace{1cm}
    \includegraphics[scale=0.3]{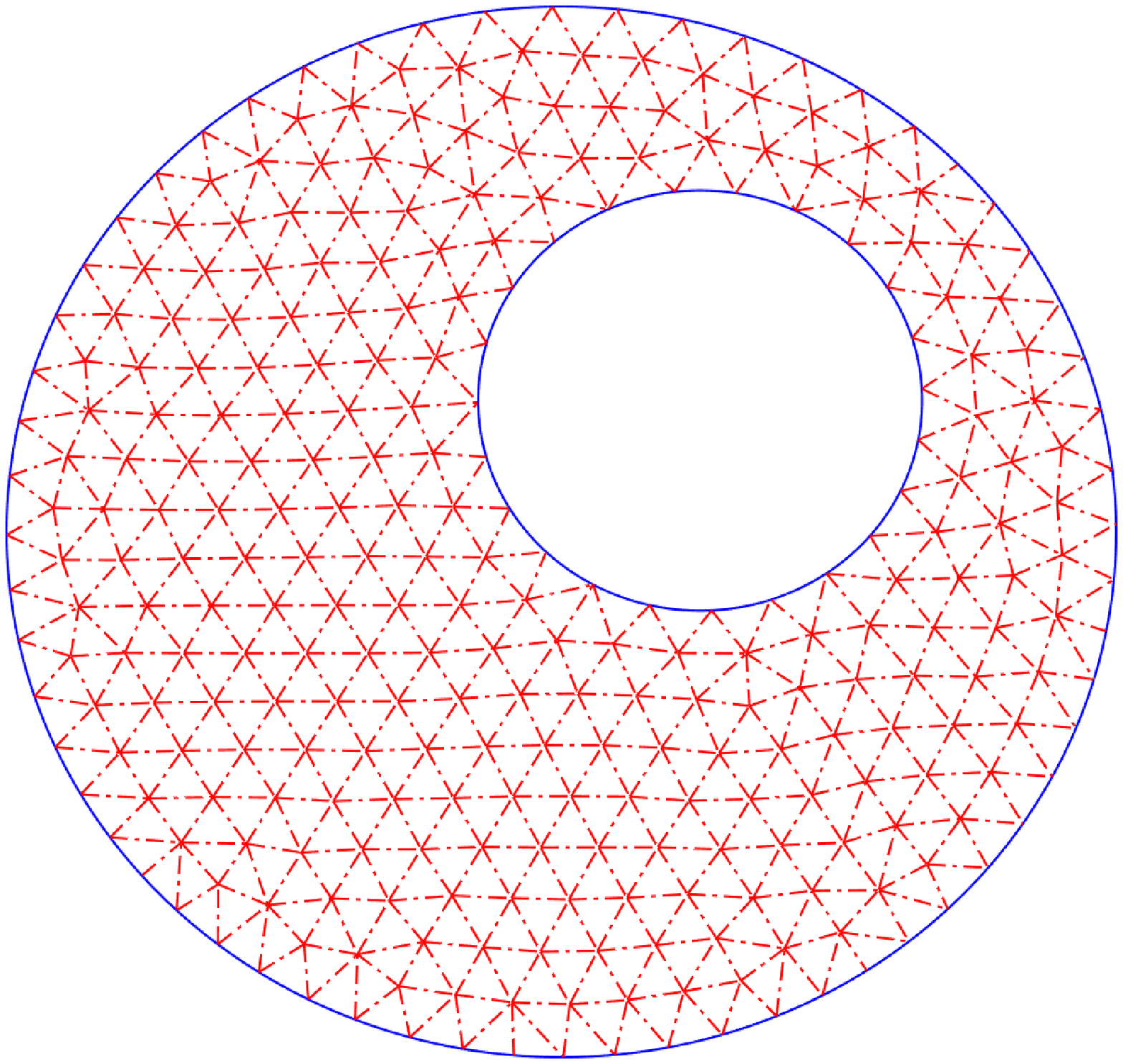}
    \caption{Two examples of curved triangular meshes with $65$ (left) and $527$ (right) cells fitting exactly the annular domain with a hole.}\label{ex2_curved_domain}
\end{figure}

\begin{figure}[htb!]
    \centering
    \includegraphics[scale=0.38]{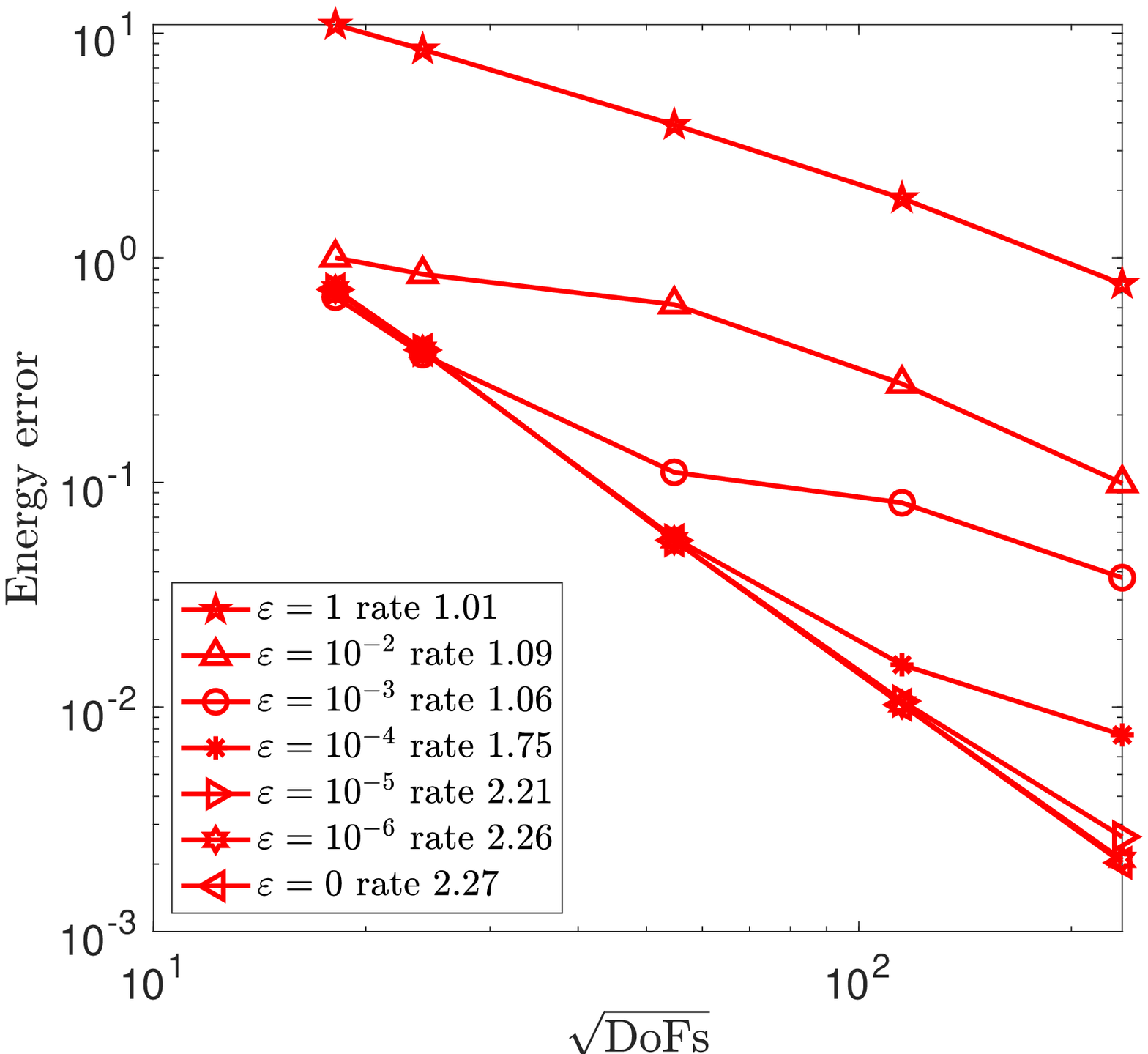}
    \includegraphics[scale=0.38]{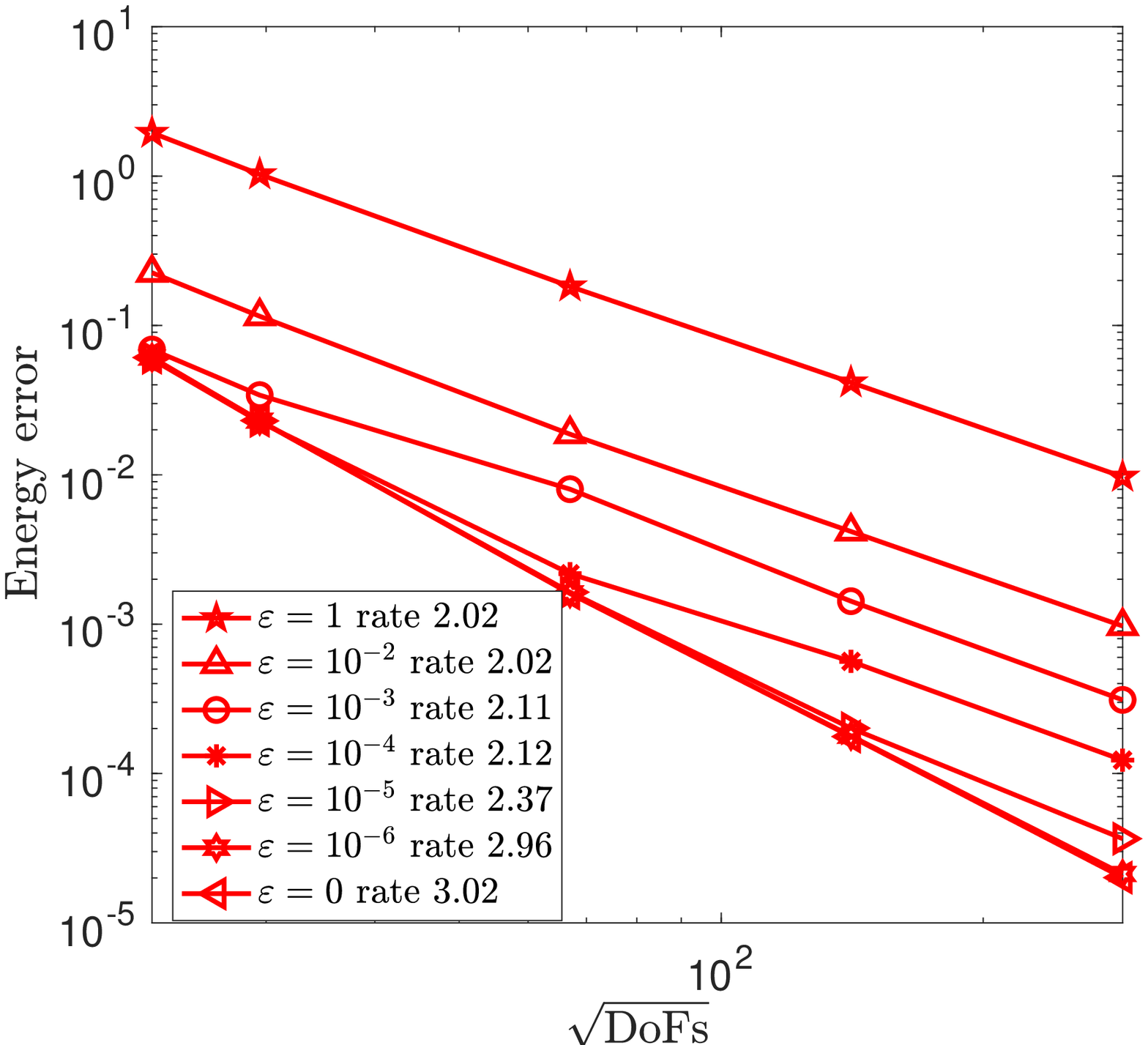} \\
    \includegraphics[scale=0.38]{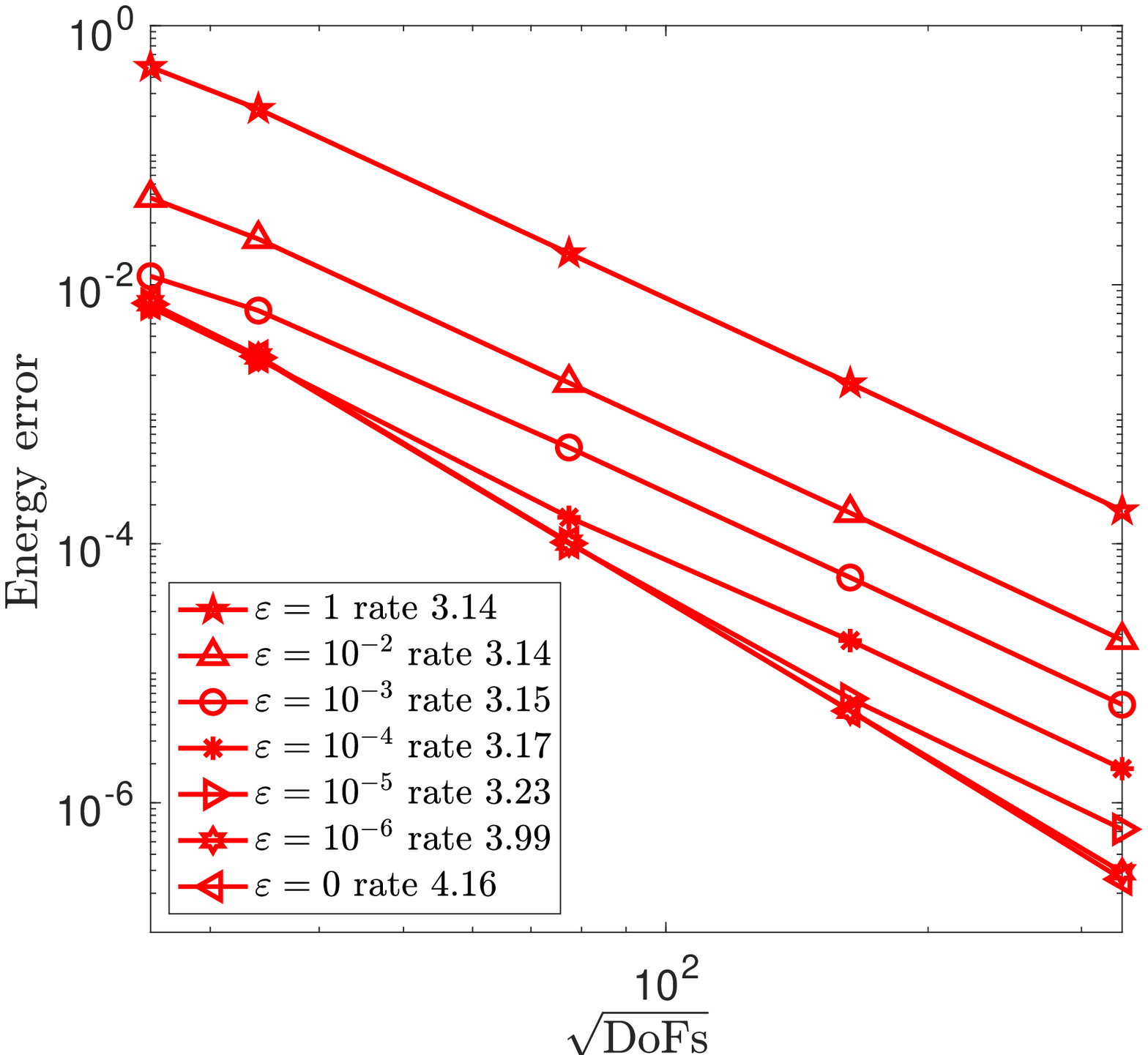}
    \includegraphics[scale=0.38]{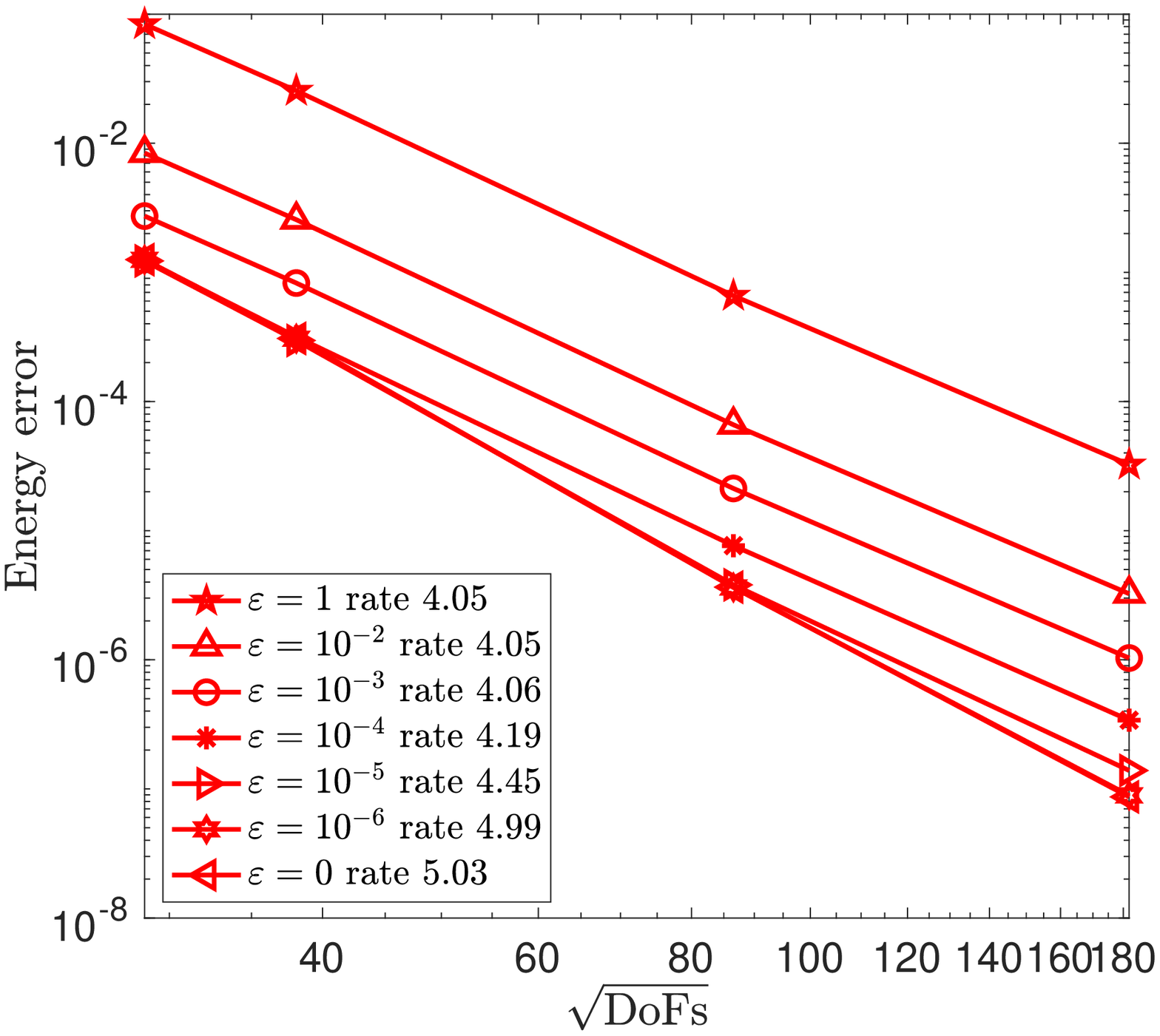}
    \caption{
        \label{ex2:h-refine, k = 0,1,2,3}
        Annular domain and smooth solution. Convergence rates in the energy seminorm on curved triangular meshes for different values of $\varepsilon$. The polynomial degree $k=0$ (top left), $k=1$  (top right), $k=2$  (bottom left) and $k=3$  (bottom right).}
\end{figure}

\subsection{Convergence rates and robustness in a domain with curved boundary}
\label{sec:res_curved}

In this second example, we consider an annular domain constructed as the unit disc centered at the  origin, with a circular hole centered at $(0.25,0.25)$ and with radius $0.4$; see Figure \ref{ex2_curved_domain}. We select $f$ and the boundary conditions so that the exact solution to \eqref{pde} is
$$
u(x,y)= (1+\sin(\pi (x^2+y^2-1)))e^{(- x^2 - y^2)}.
$$
We consider a quasi-uniform sequence of triangular meshes composed of 65, 109, 527, 2266, and 9411 triangular elements. All the meshes fit the domain $\Omega$ exactly, and for every mesh in the sequence, each interior cell has only straight edges, whereas each boundary cell has one curved edge that exactly fits the boundary of $\Omega$.

We perform the same numerical experiment as in the previous section and report the results in Figure \ref{ex2:h-refine, k = 0,1,2,3}. The conclusions are the same as in the previous test case. The transition from the $O(h^{k+2})$ to the $O(h^{k+1})$ regimes is clearly visible in all cases.

\subsection{Test case with boundary layer}
\label{sec:layer}

\begin{figure}[htb!]
    \begin{center}
        \begin{tabular}{cc}
            \includegraphics[scale=0.3]{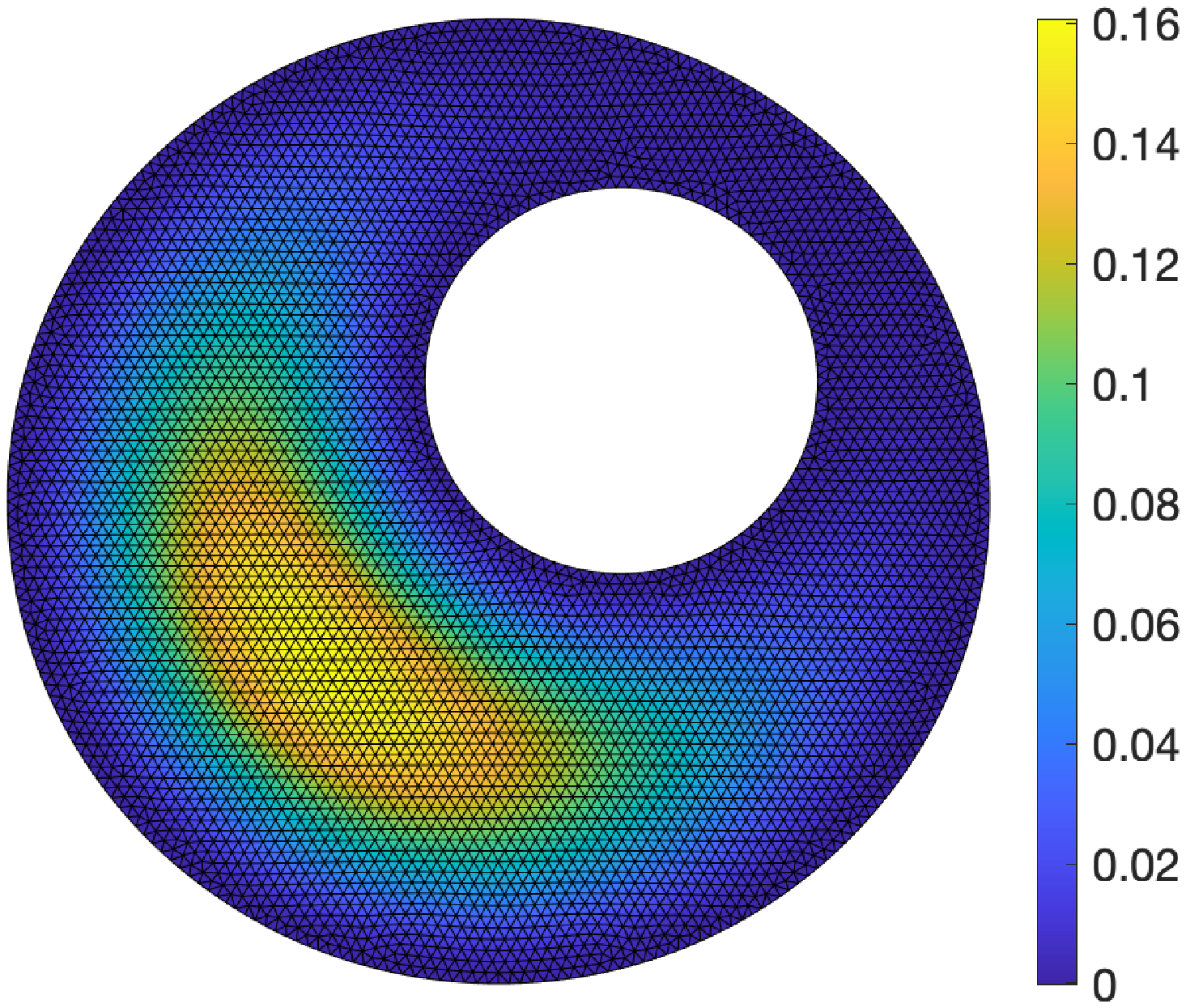}
            \includegraphics[scale=0.3]{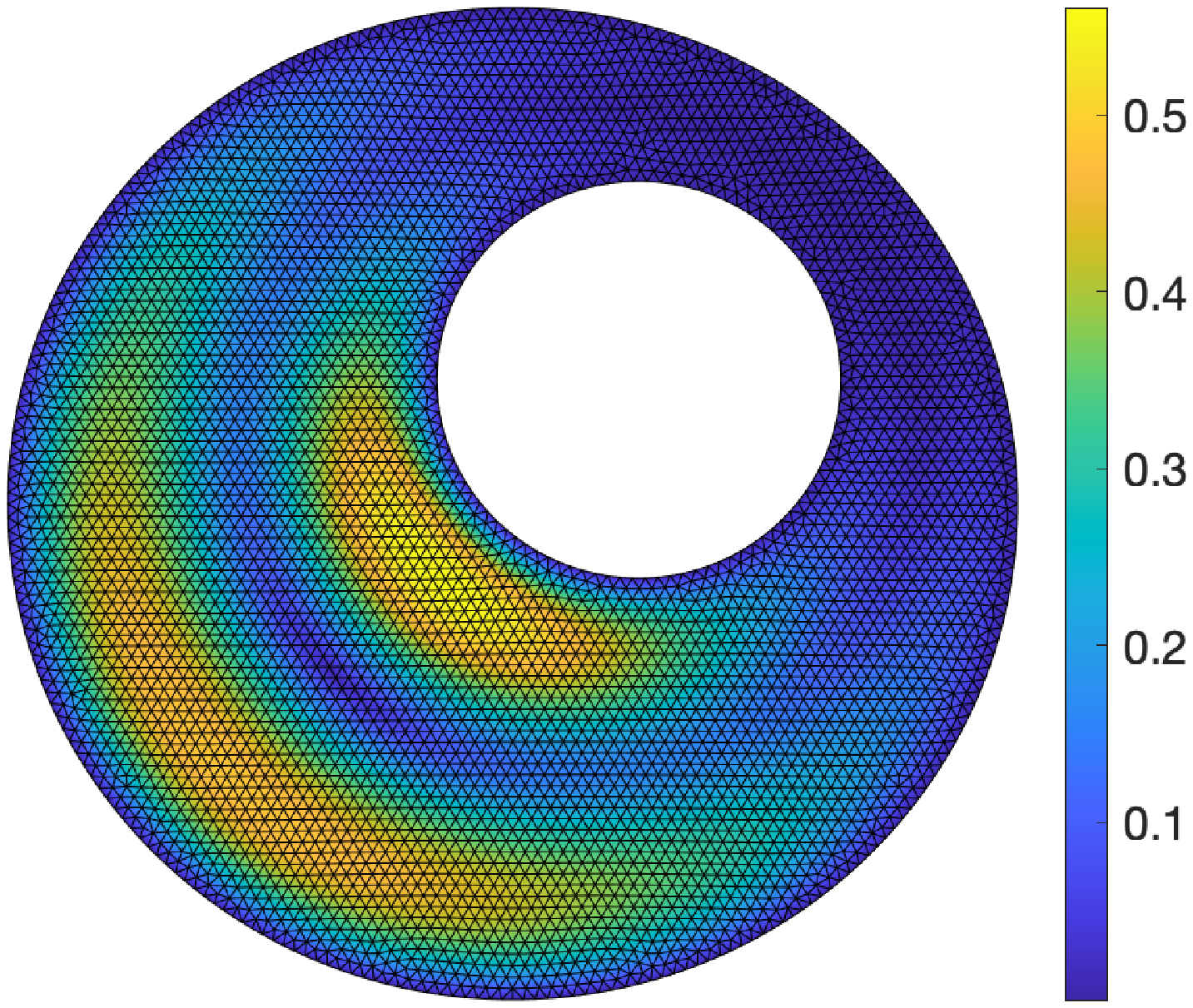}
            \includegraphics[scale=0.3]{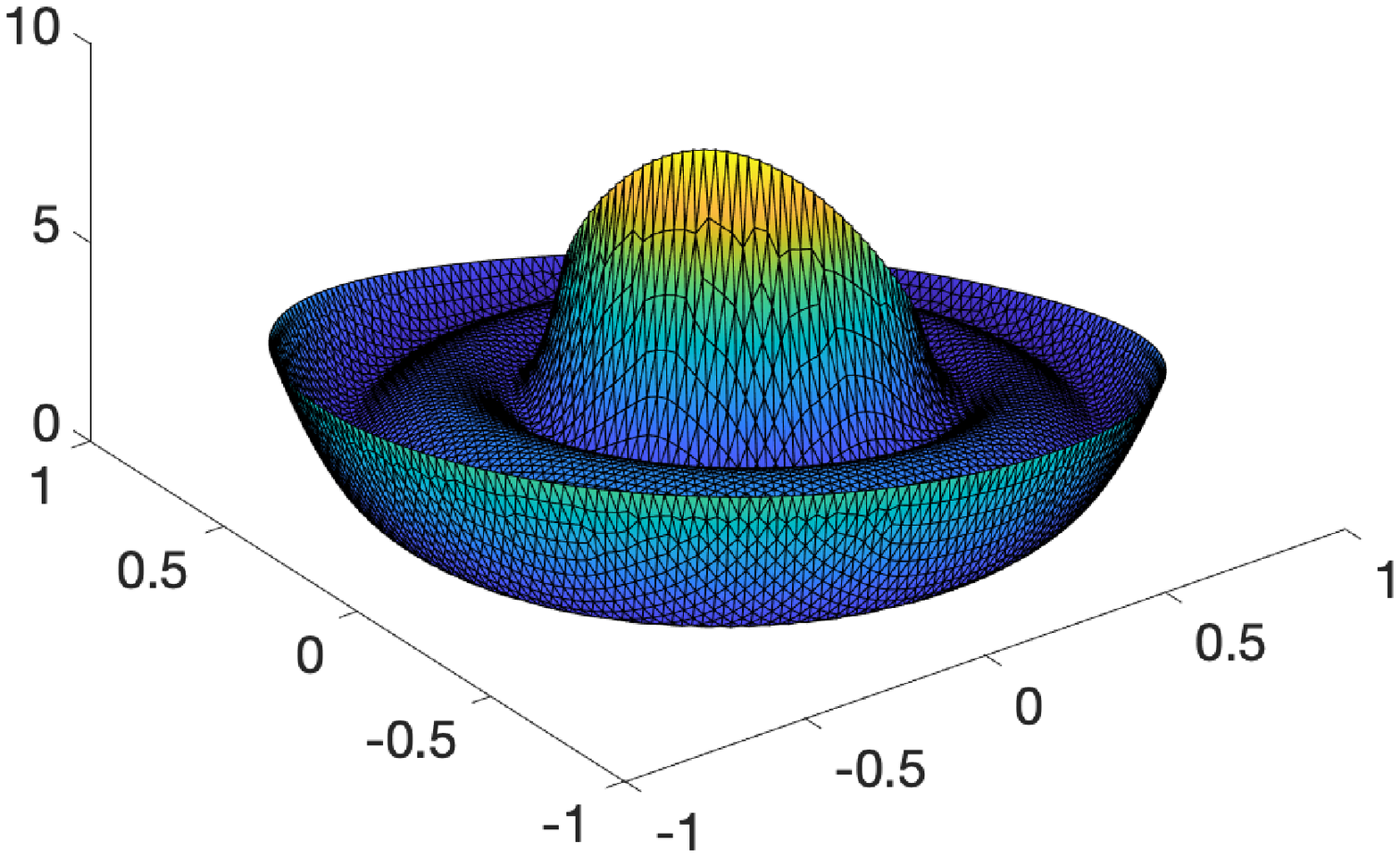} \\
            \includegraphics[scale=0.3]{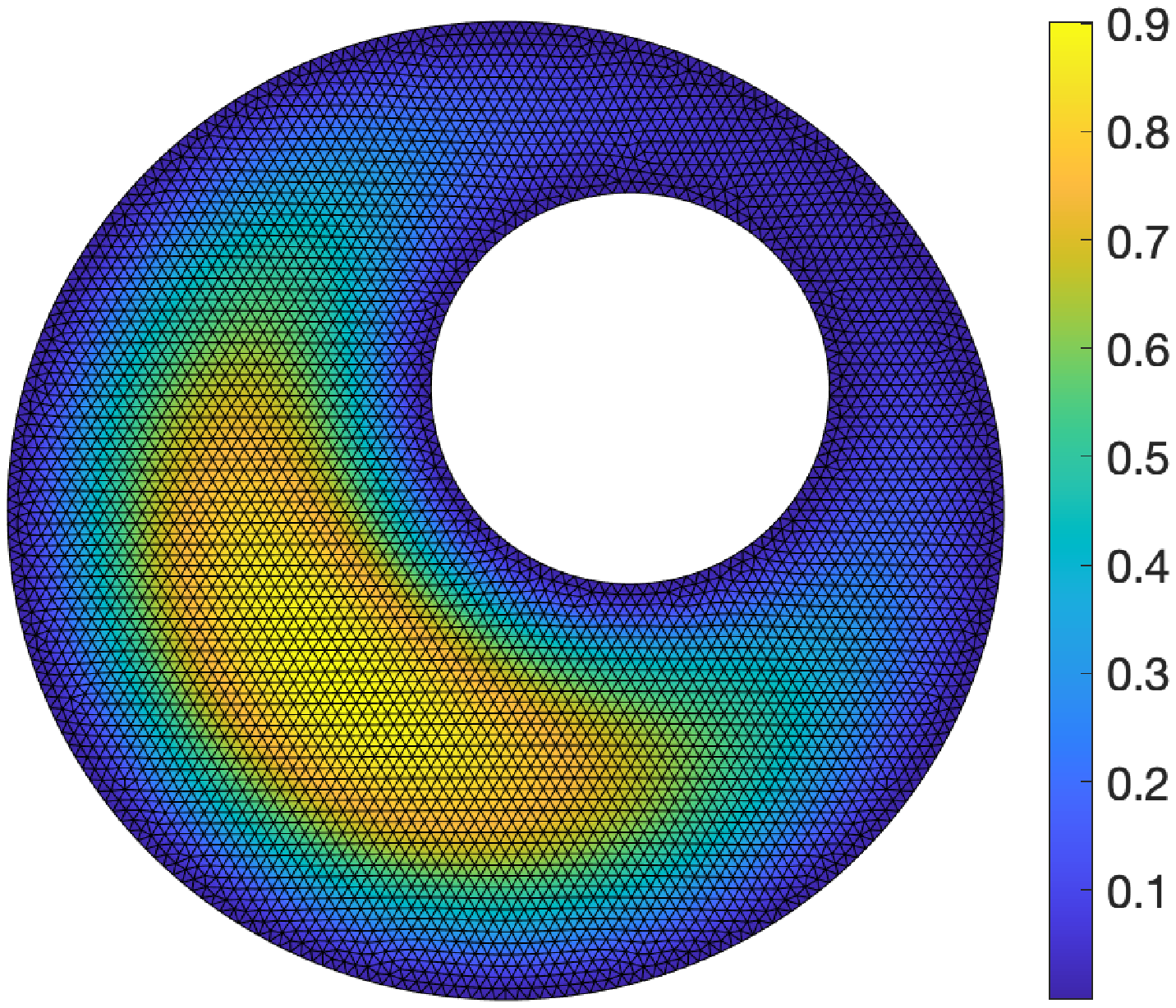}
            \includegraphics[scale=0.3]{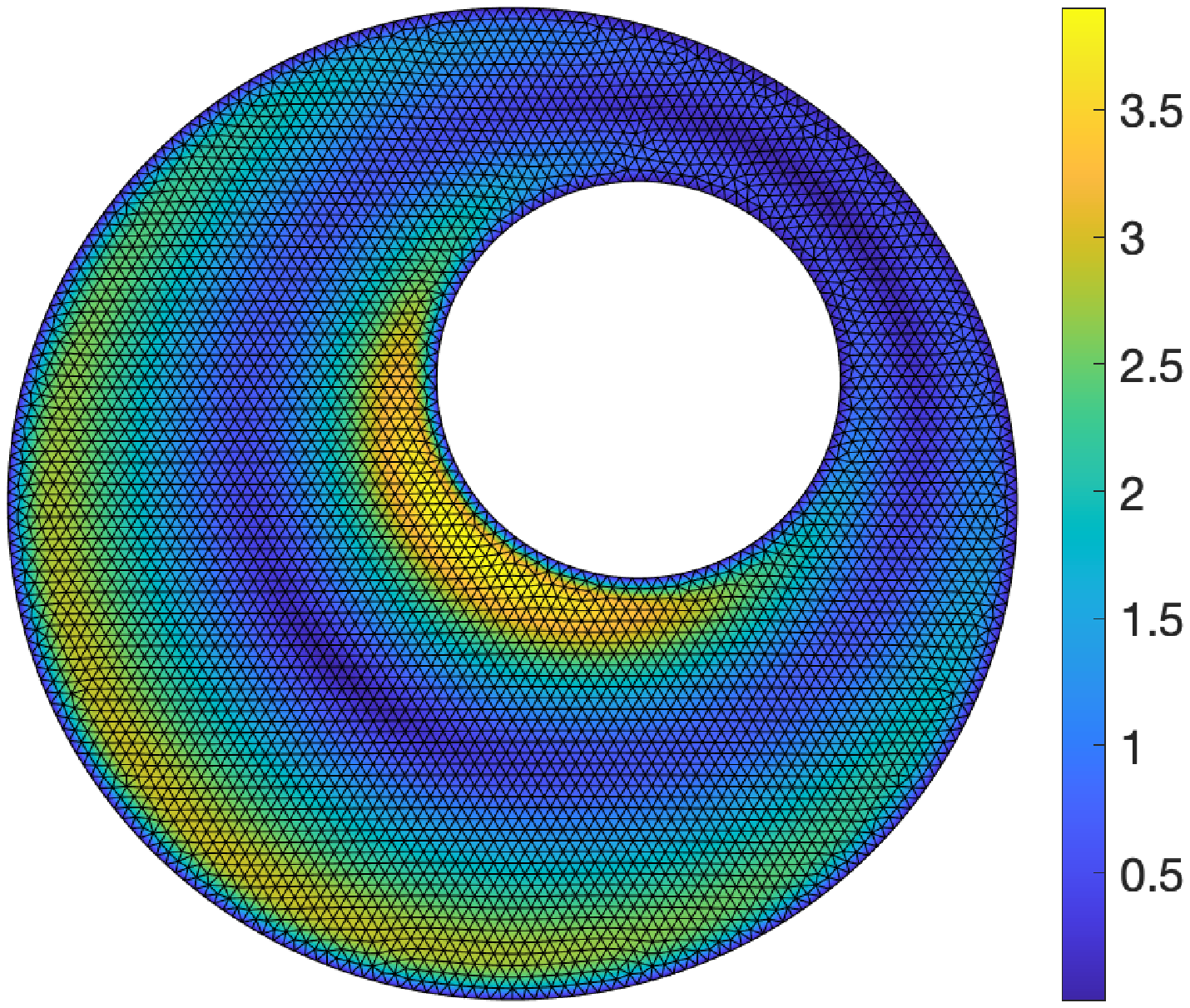}
            \includegraphics[scale=0.3]{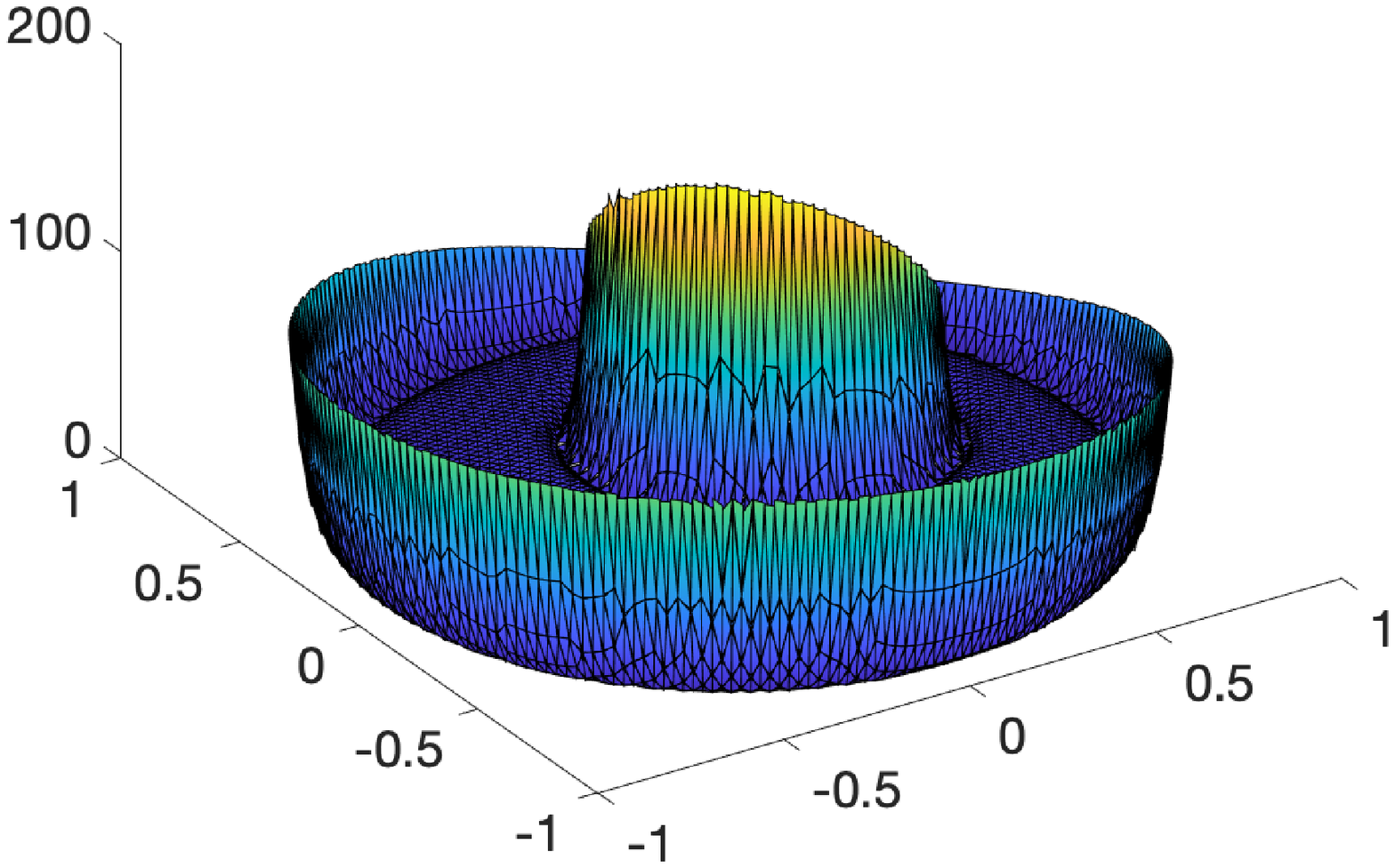}
        \end{tabular}
    \end{center}
    \caption{Test case with boundary layer on the mesh composed of $9411$ cells, $\varepsilon=10^{-1}$ (top row) and $\varepsilon=10^{-3}$ (bottom row). Left column: reconstructed solution; middle column: piecewise gradient of reconstructed solution (Euclidean norm); right column: piecewise Hessian of reconstructed solution (Frobenius norm).}\label{ex3_eps_1_2_3_fine_mesh}
\end{figure}

\begin{figure}[htb!]
    \begin{center}
        \begin{tabular}{cc}
            \includegraphics[scale=0.3]{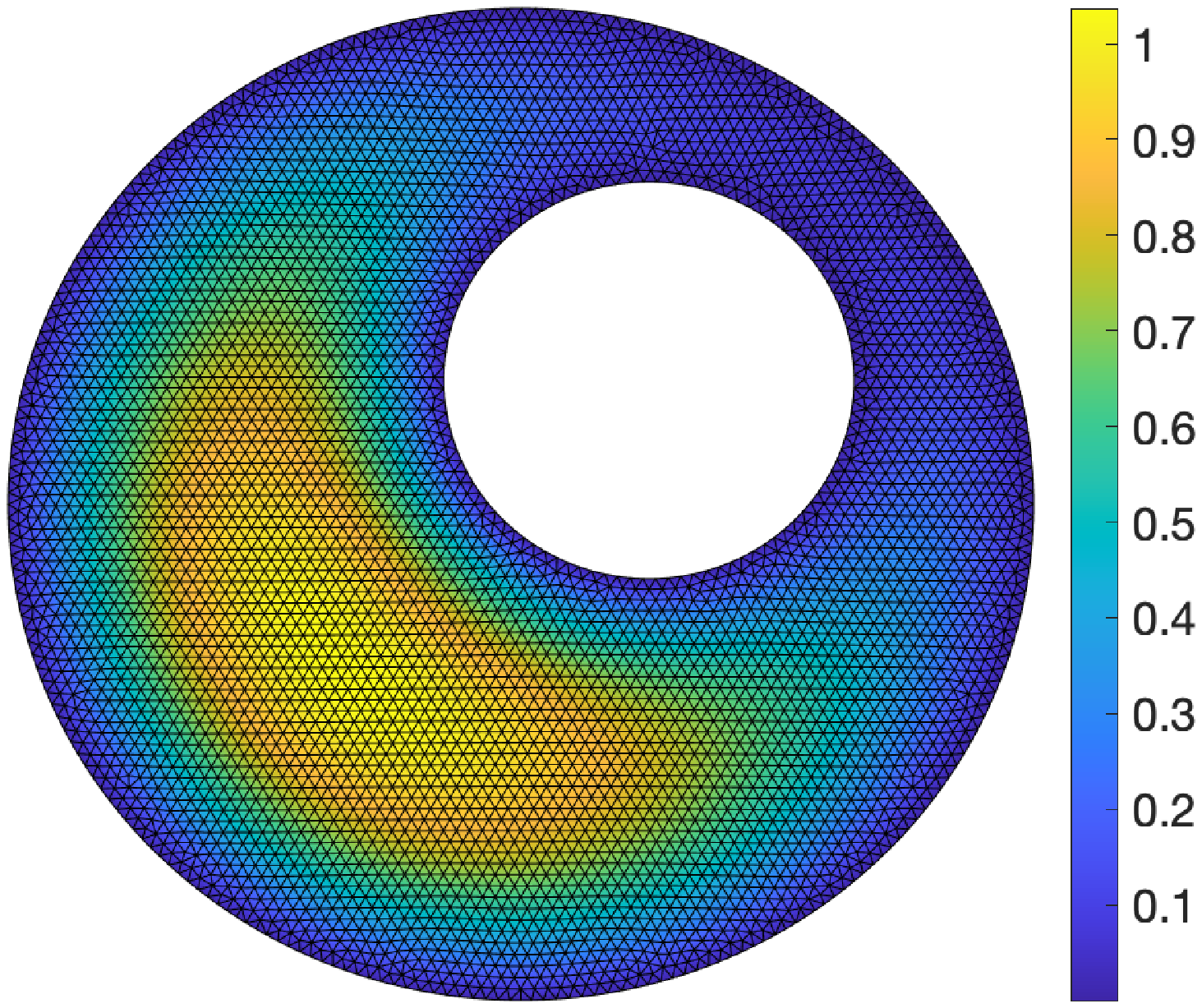}
            \includegraphics[scale=0.3]{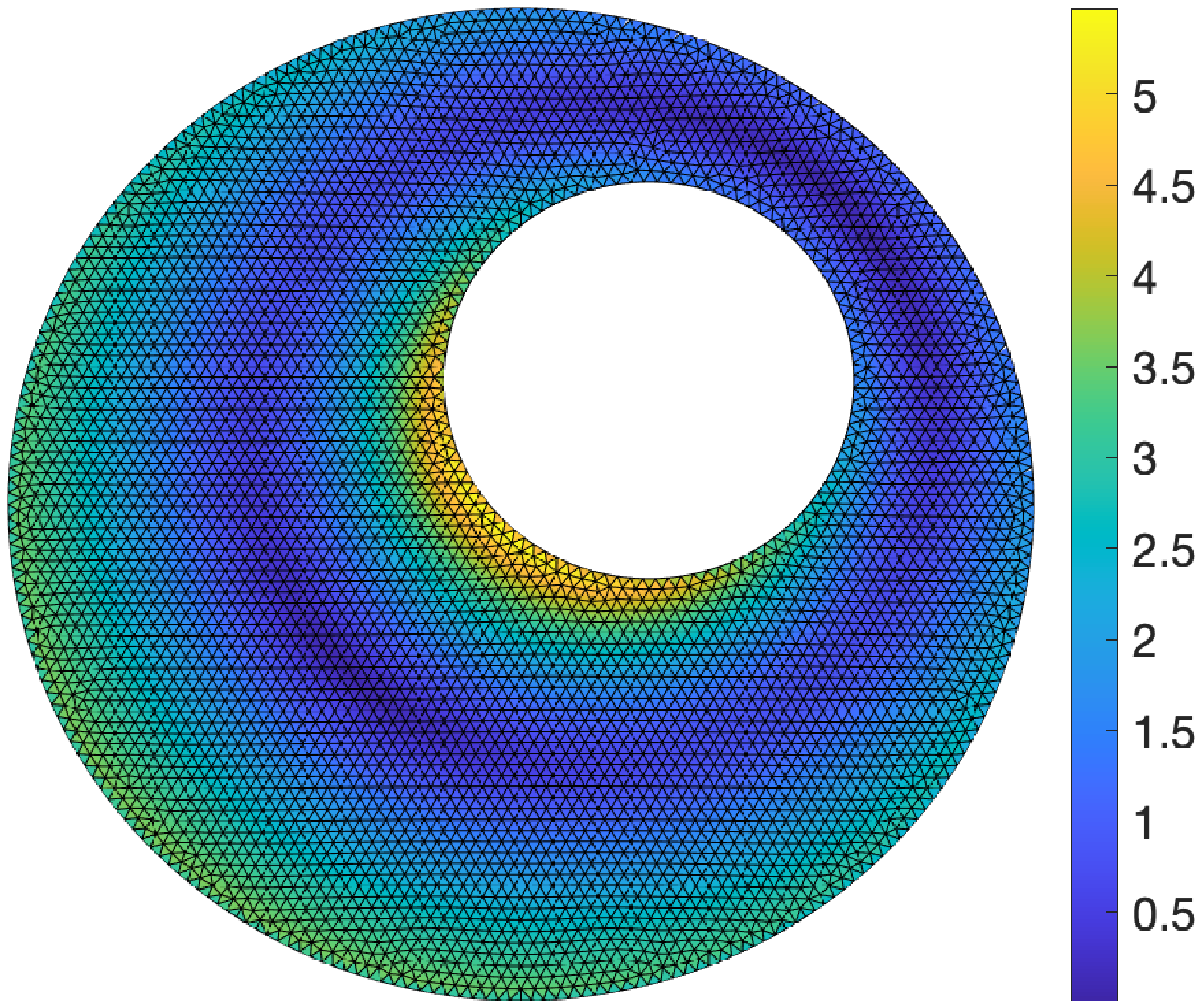} \includegraphics[scale=0.3]{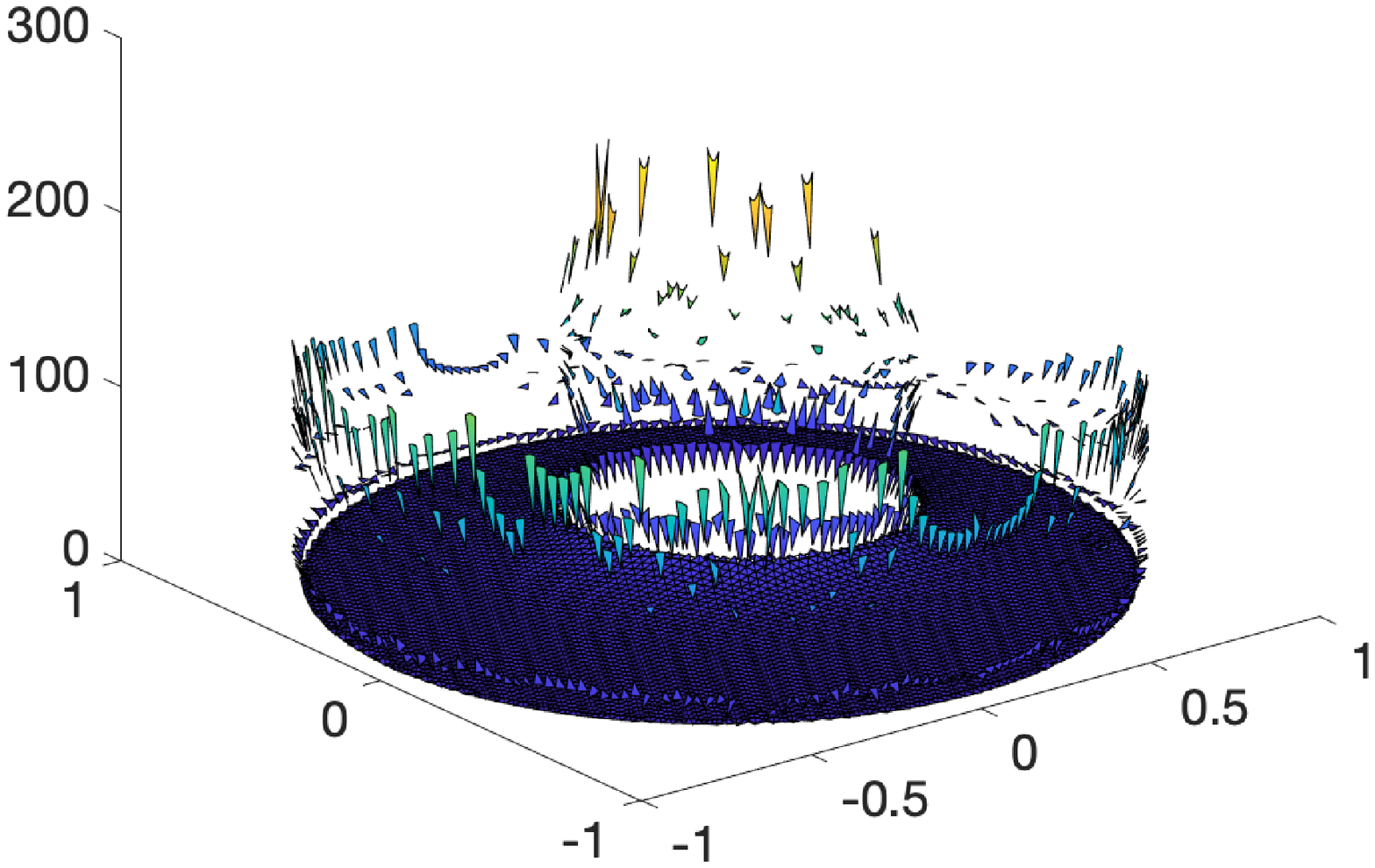}\\
            \includegraphics[scale=0.3]{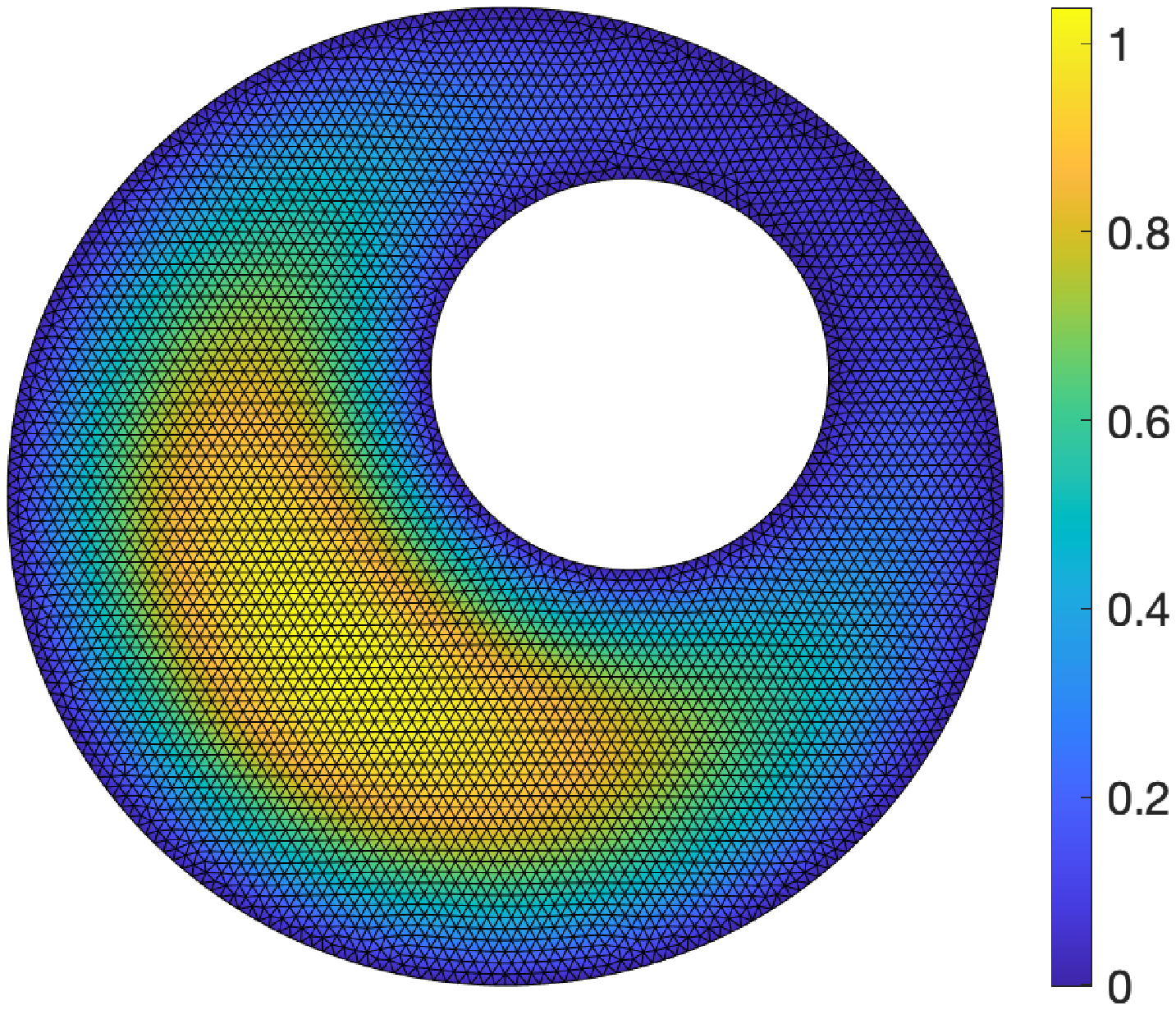}
            \includegraphics[scale=0.3]{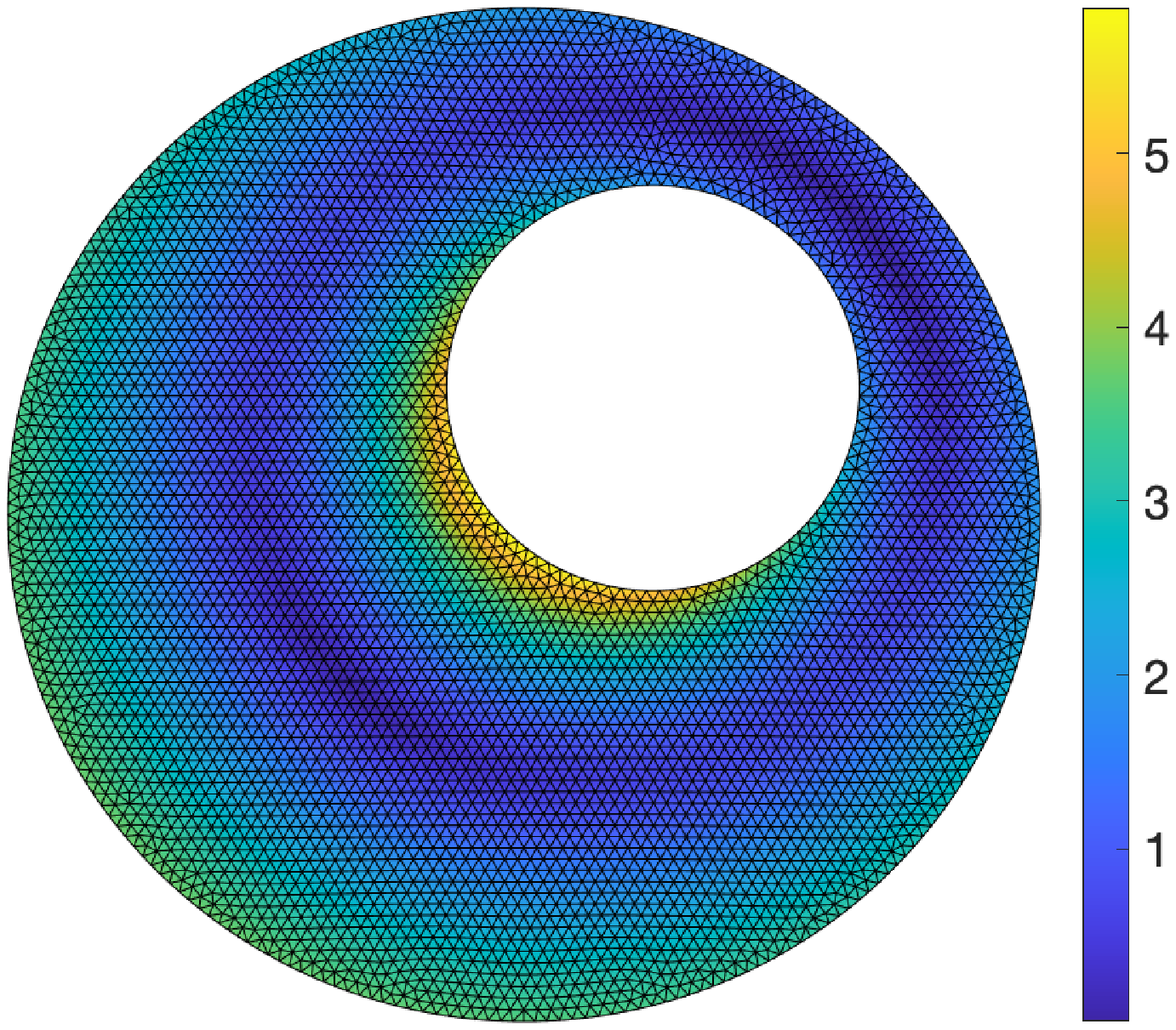} \includegraphics[scale=0.3]{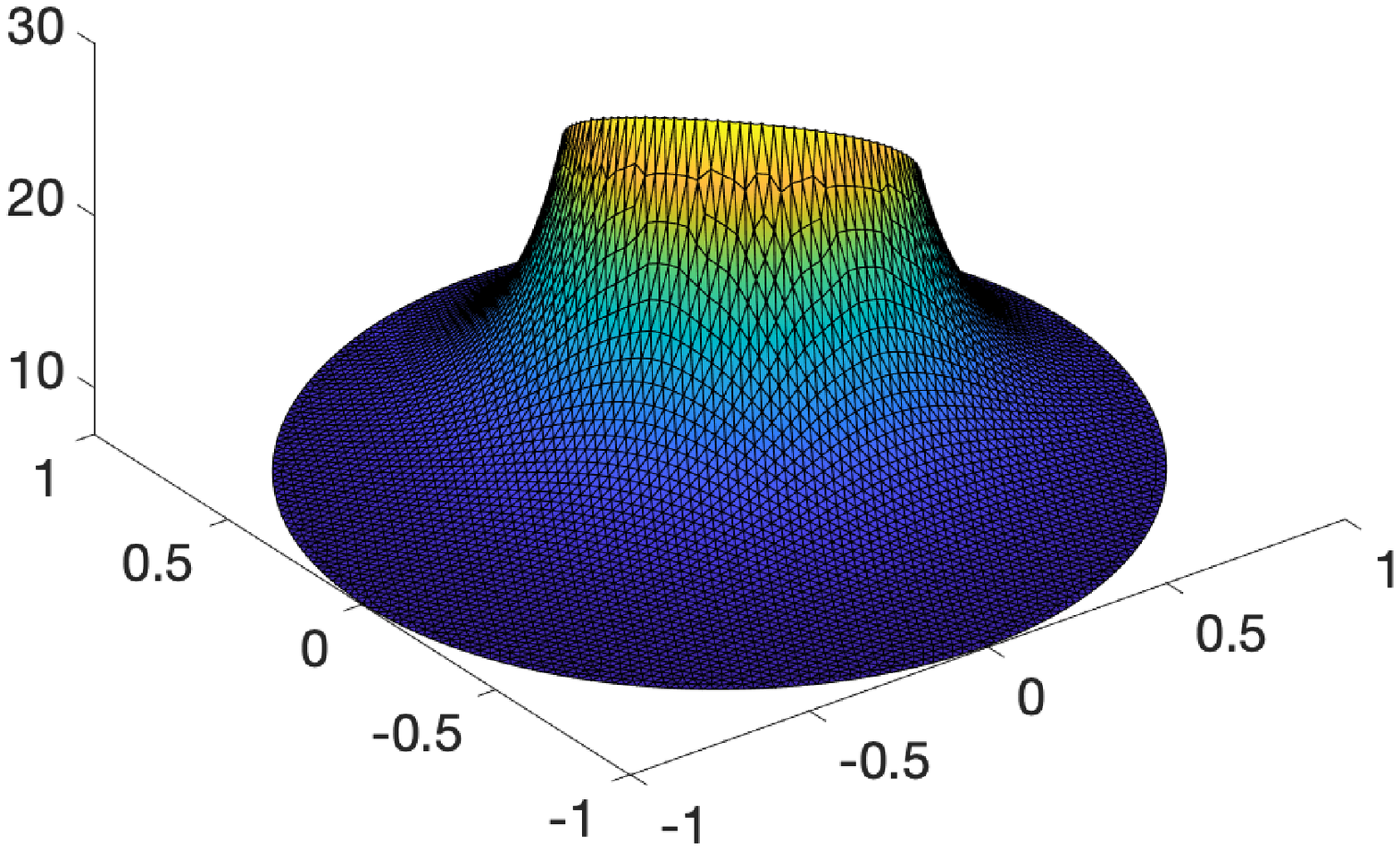}
        \end{tabular}
    \end{center}
    \caption{Test case with boundary layer on the mesh composed of $9411$ cells, $\varepsilon=10^{-6}$ (top row) and $\varepsilon=0$ (bottom row).  Left column: reconstructed solution; middle column: piecewise gradient of reconstructed solution (Euclidean norm); right column: piecewise Hessian of reconstructed solution (Frobenius norm).}\label{ex3_eps_0_coarse_mesh}
\end{figure}

We conclude this series of numerical experiments with a somewhat more challenging test case featuring a boundary layer. We consider the same annular domain as in the previous section,
we set the source term to $f := 10$ and we enforce homogeneous boundary conditions. The considered values for the singular perturbation parameter are $\varepsilon\in\{10^{-1},10^{-2},10^{-3},10^{-6},0\}$. In all cases, the analytical solution is unknown. Numerical solutions are obtained on the curved triangular meshes considered in the previous section using the polynomial degree $k=1$. We report in Figure~\ref{ex3_eps_1_2_3_fine_mesh} the reconstructed solution $R_h(\widehat{u}_h)$ defined as $R_h(\widehat{u}_h)|_K:=R_K(\widehat{u}_K)$ for all $K\in\mesh$ (with
$R_K(\widehat{u}_K)$ defined just above Theorem~\ref{Theorem: main}), its piecewise gradient (Euclidean norm), and its piecewise Hessian (Frobenius norm) for $\varepsilon=10^{-1}$ and $\varepsilon=10^{-3}$ on the mesh composed of $9411$ curved triangular cells. Since $h=0.0344$ for this mesh, the boundary layer is well resolved for $\varepsilon=10^{-1}$ and barely resolved for $\varepsilon=10^{-3}$.
Notice that $R_h(\widehat{u}_h)$ is a piecewise cubic polynomial since we are using here $k=1$. We observe in Figure~\ref{ex3_eps_1_2_3_fine_mesh} that the presence of the boundary layer is reflected by larger values of the Hessian near the boundary, whereas the reconstructed solution and its piecewise gradient take moderate values. To illustrate that the HHO method remains stable even if the boundary layer is not resolved, we present in Figure \ref{ex3_eps_0_coarse_mesh} the same quantities as in Figure~\ref{ex3_eps_1_2_3_fine_mesh} obtained on the same mesh, but this time with $\varepsilon=10^{-6}$ and $\varepsilon=0$. We notice in particular that the larger values of the Hessian remain localized close to the boundary for $\varepsilon=10^{-6}$, whereas the solution to the second-order PDE is recovered for $\varepsilon=0$.

\begin{figure}[htb!]
    \begin{center}
        \begin{tabular}{cc}
             \hspace{-0.9cm}
            \includegraphics[scale=0.33]{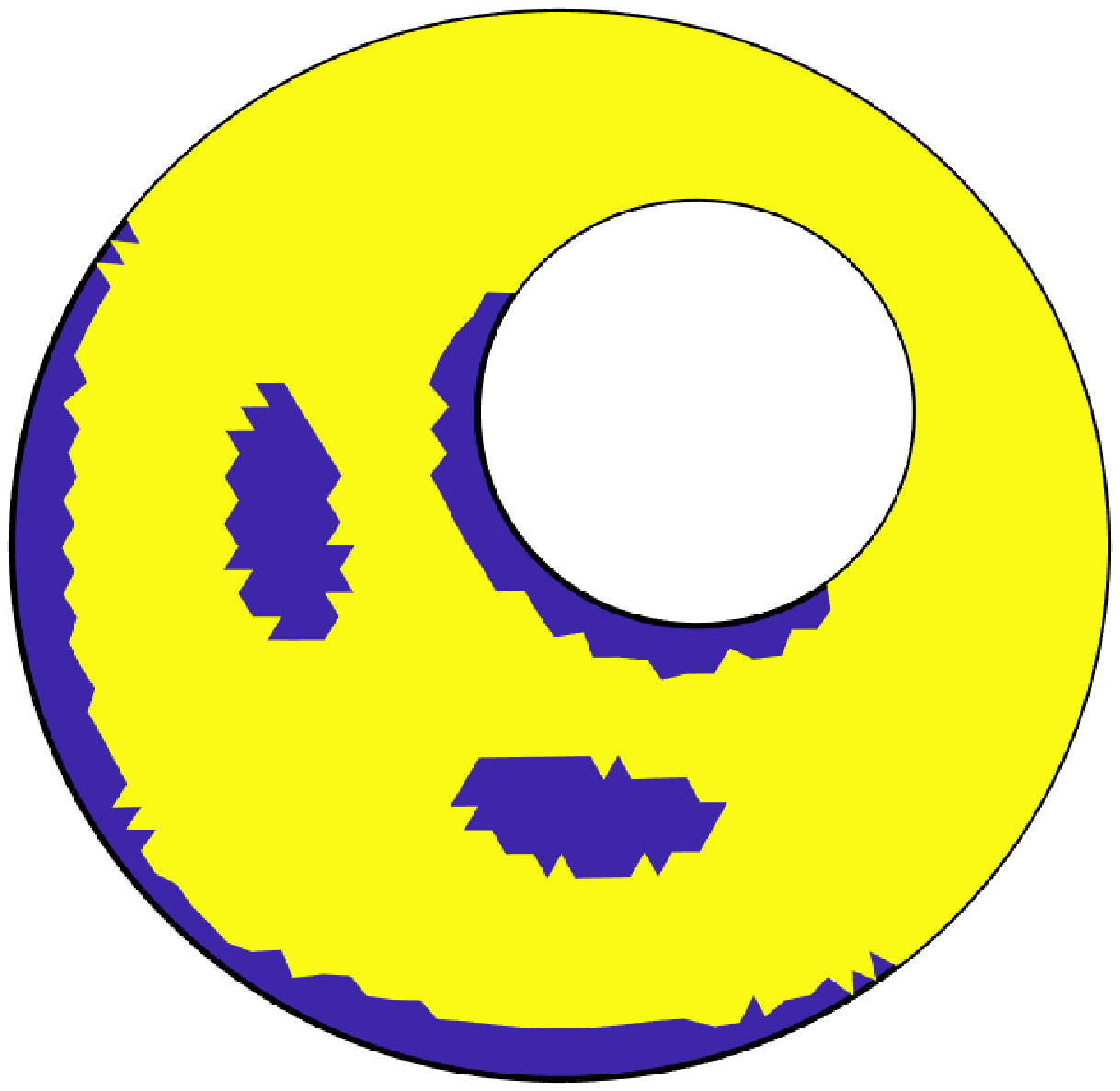}
           \includegraphics[scale=0.33]{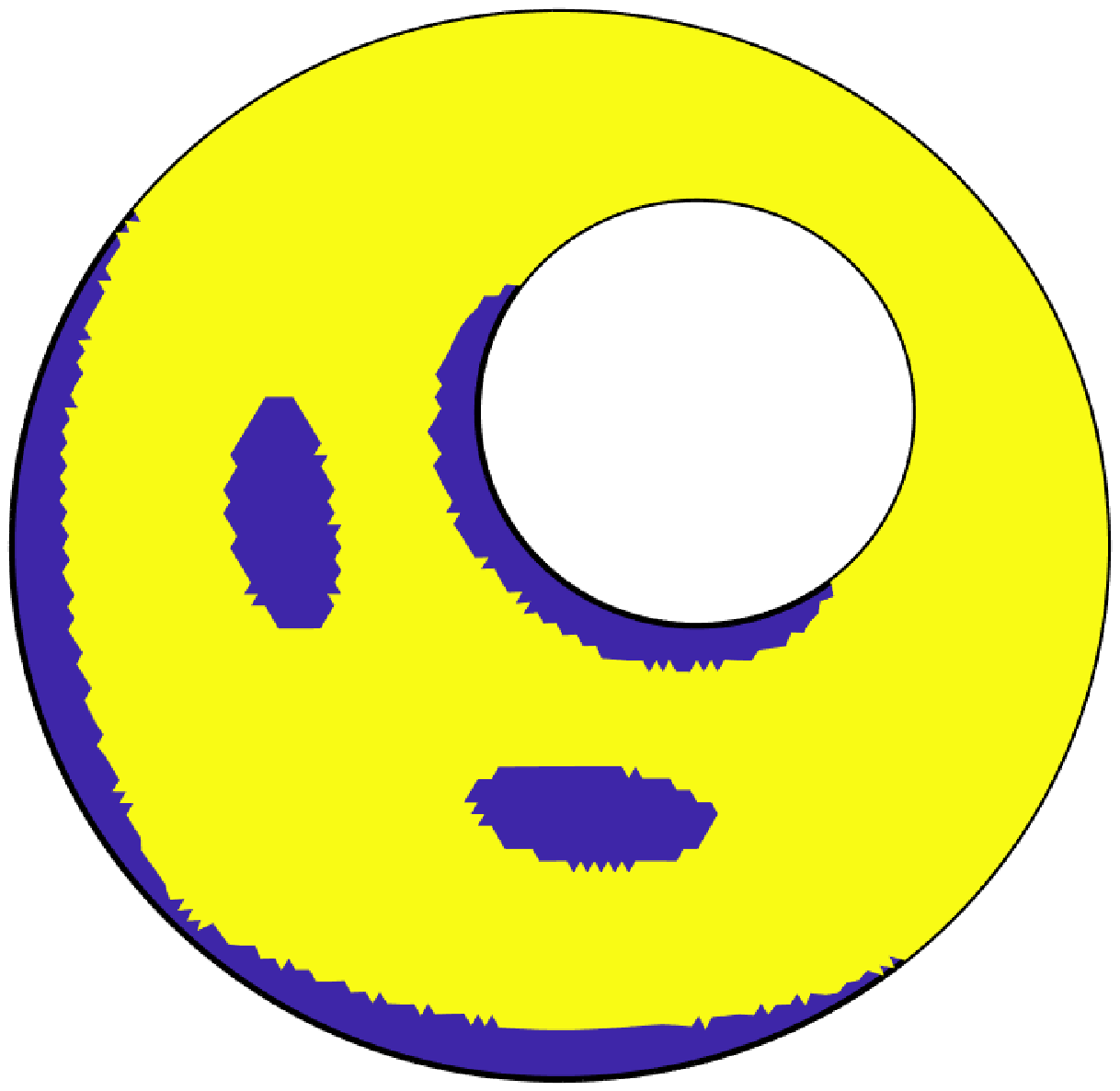}
              \includegraphics[scale=0.33]{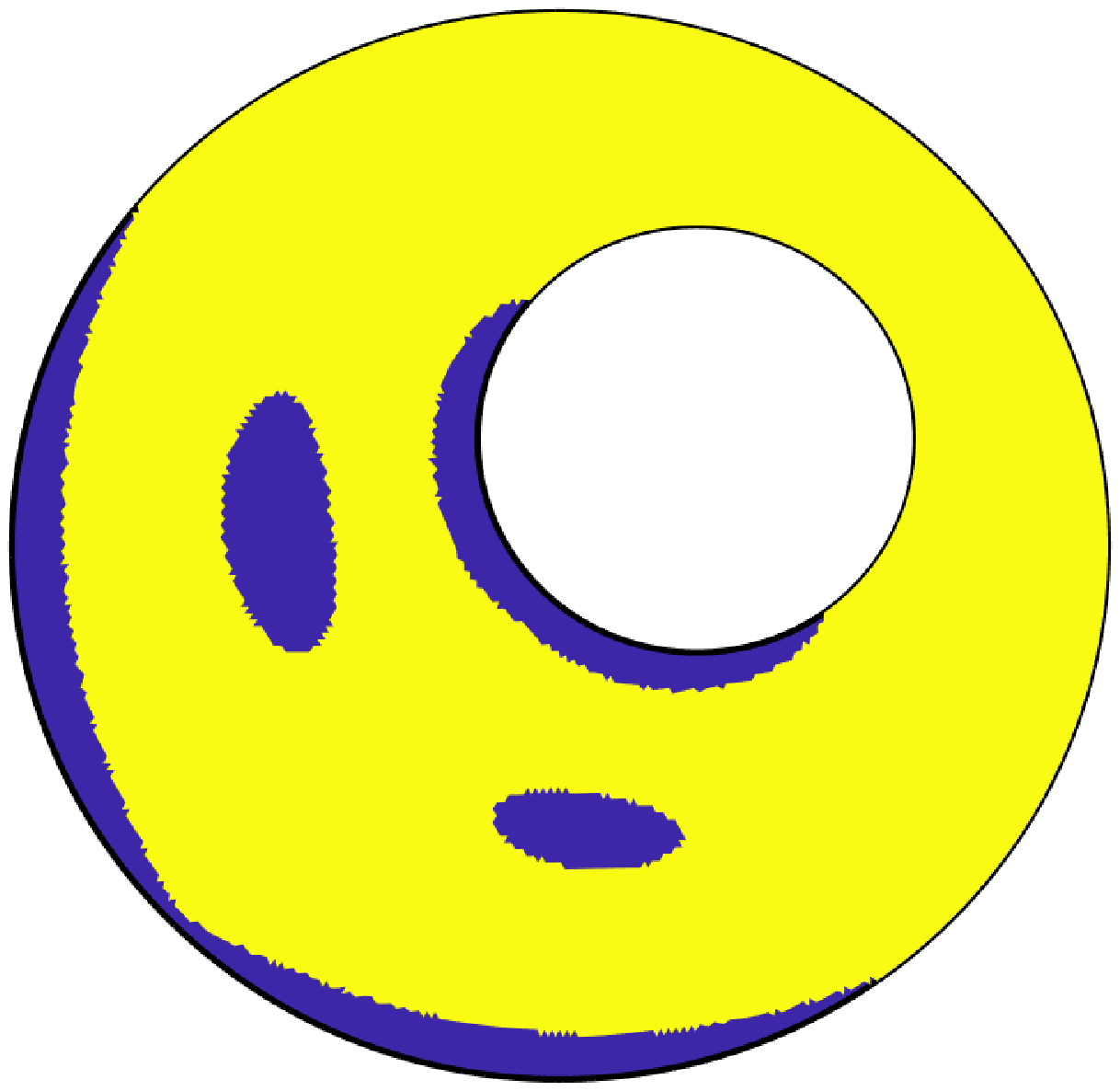}
            \\
             \hspace{-0.9cm}
            \includegraphics[scale=0.33]{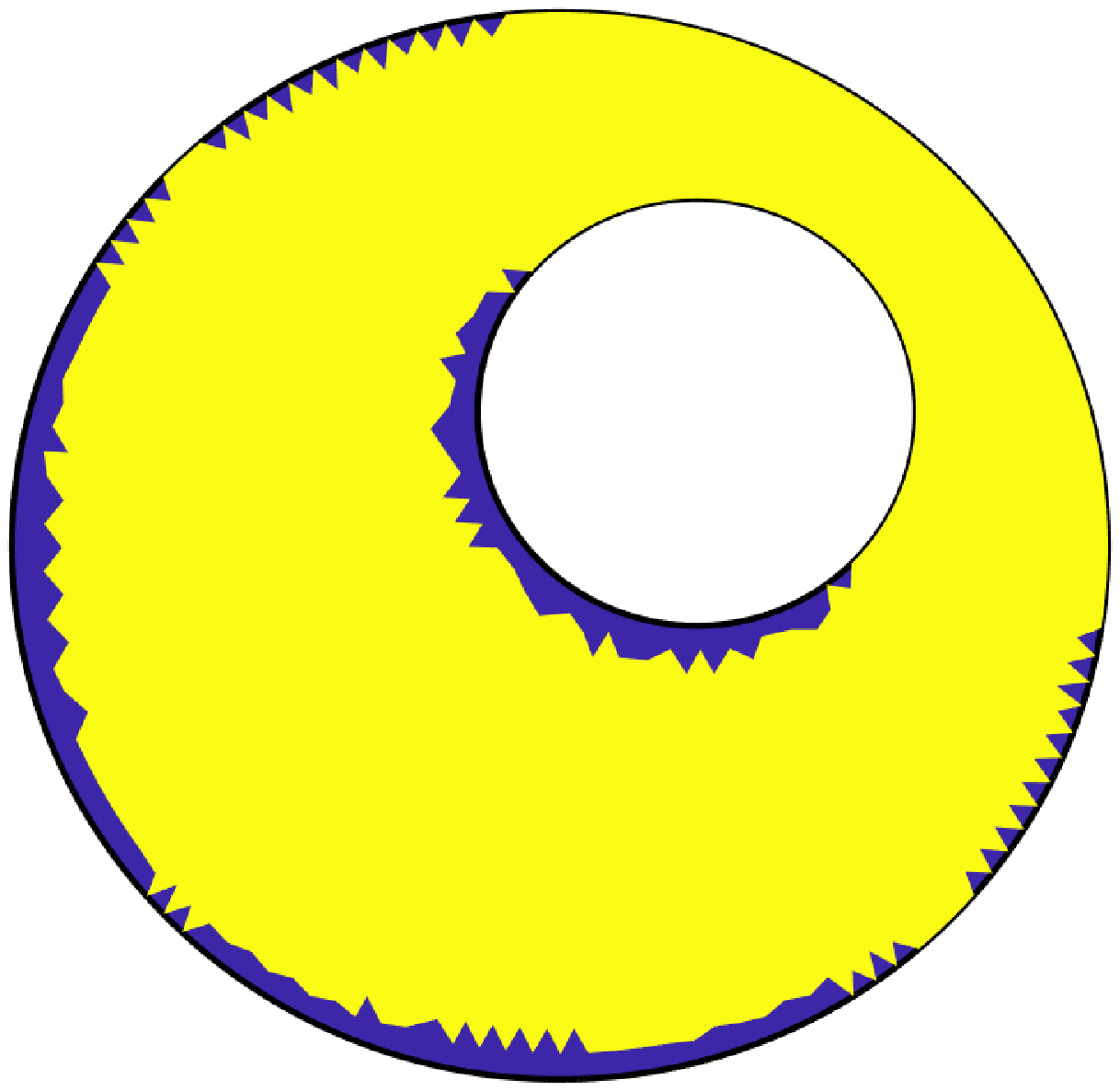}
            \includegraphics[scale=0.33]{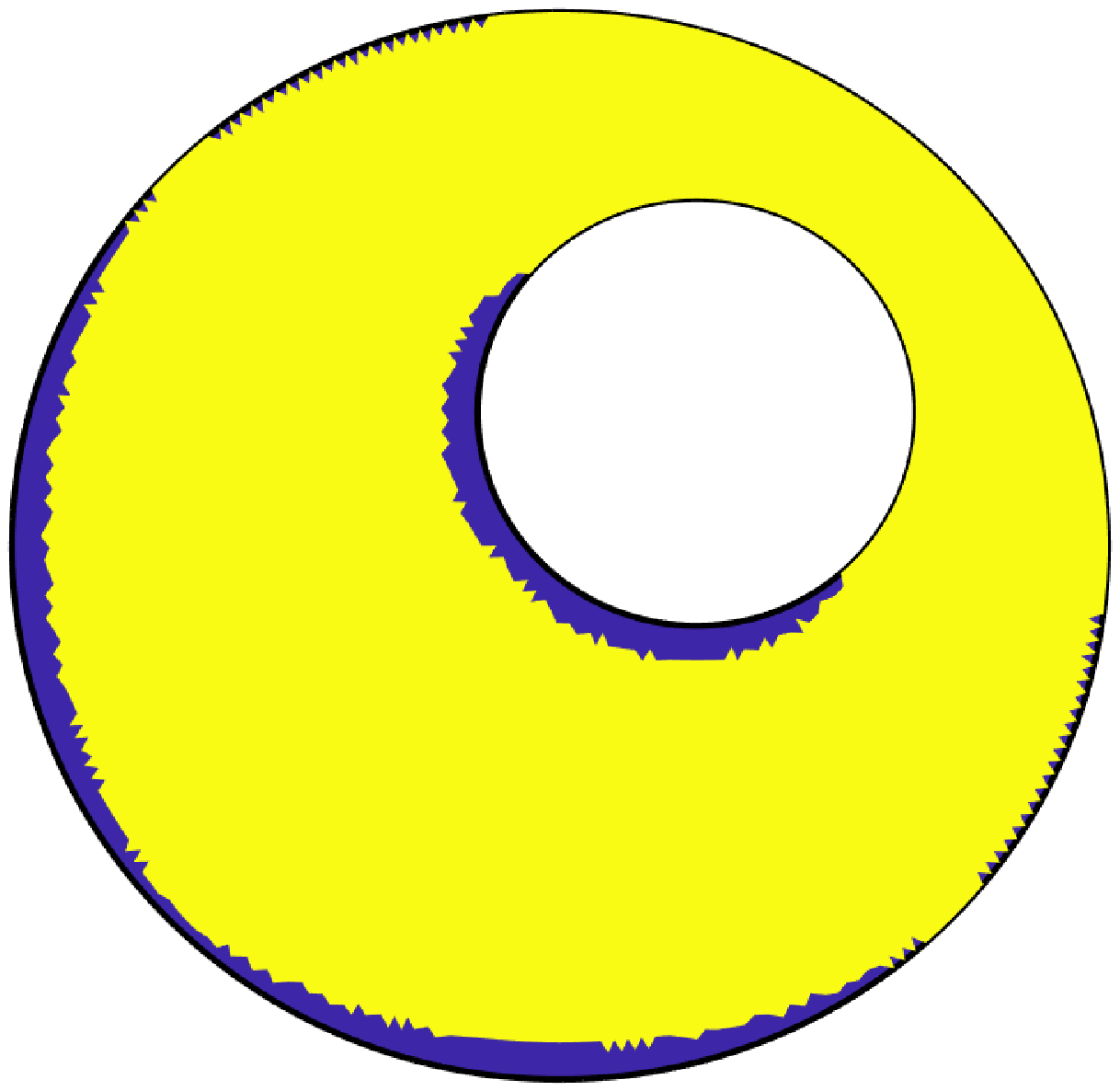}
              \includegraphics[scale=0.33]{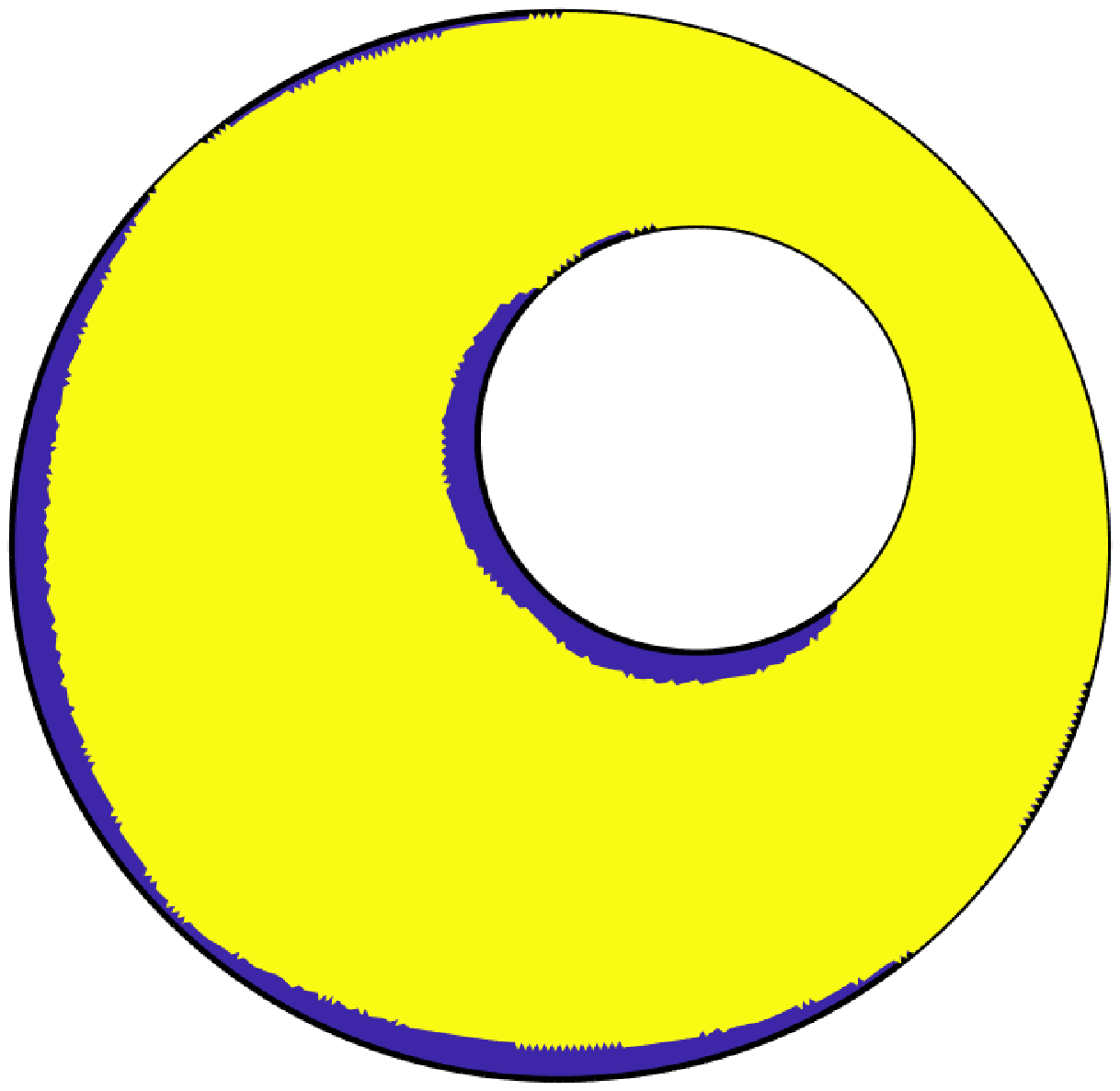}
            \\
             \hspace{-0.9cm}
            \includegraphics[scale=0.33]{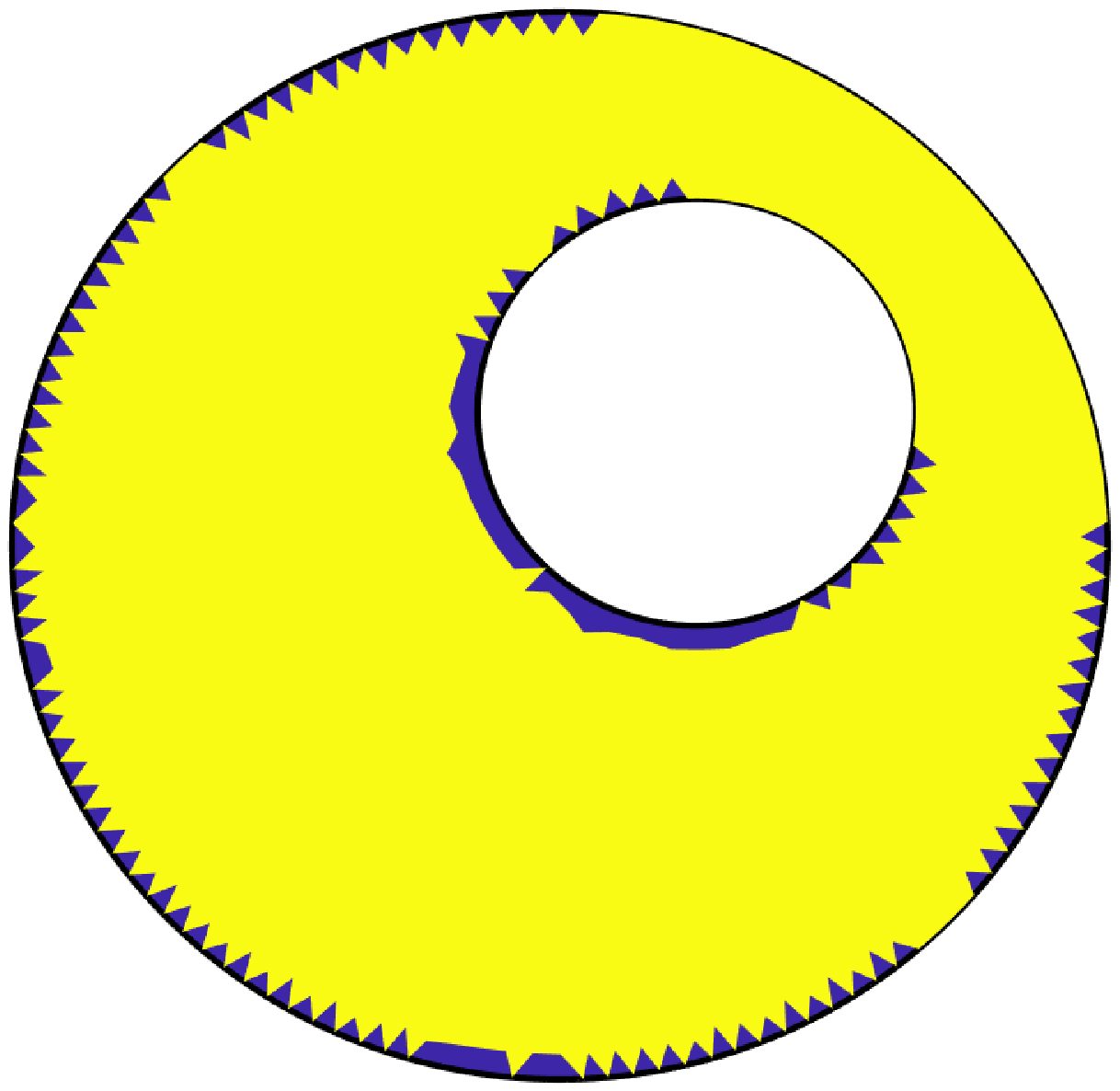}
            \includegraphics[scale=0.33]{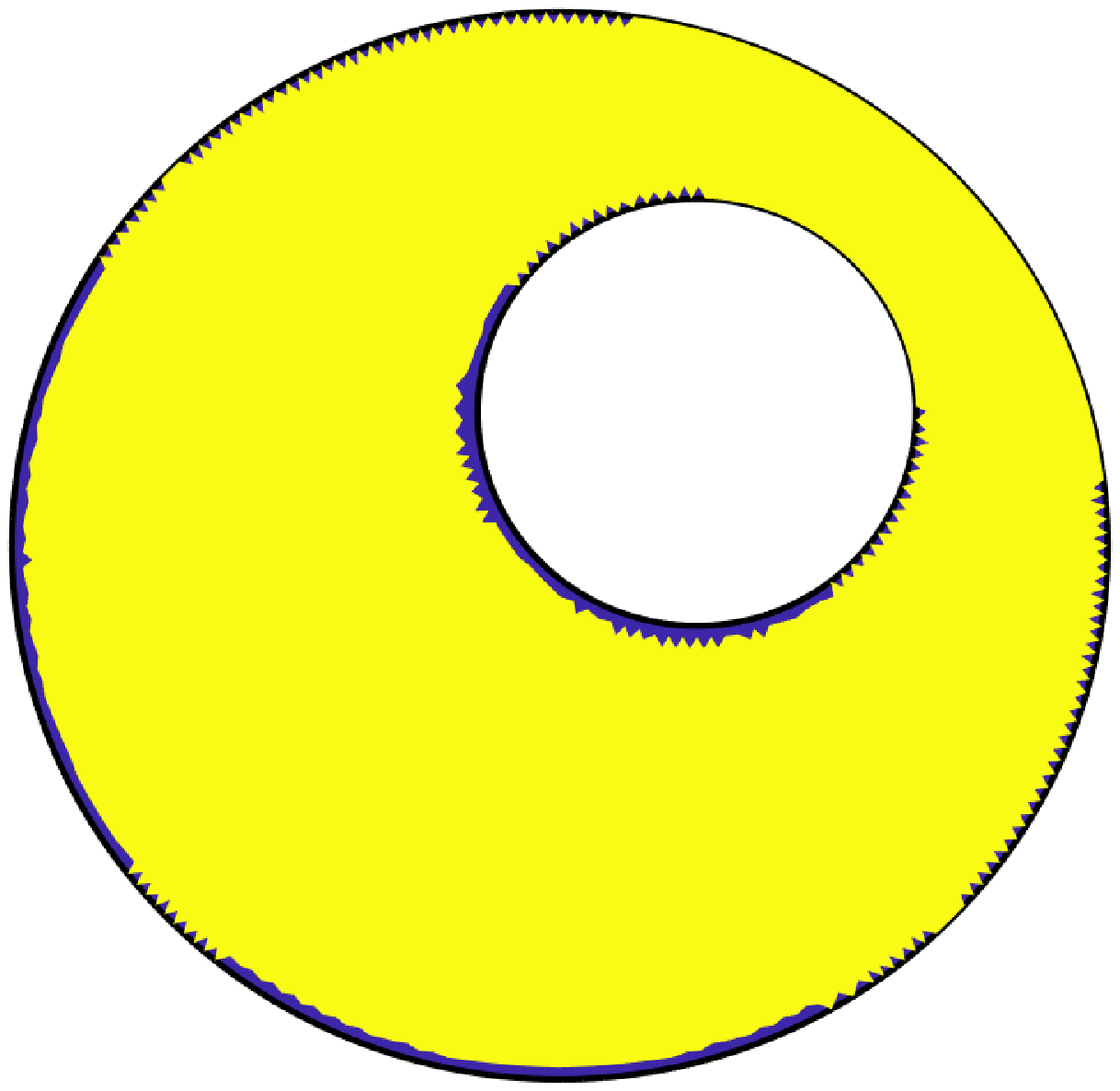}
             \includegraphics[scale=0.33]{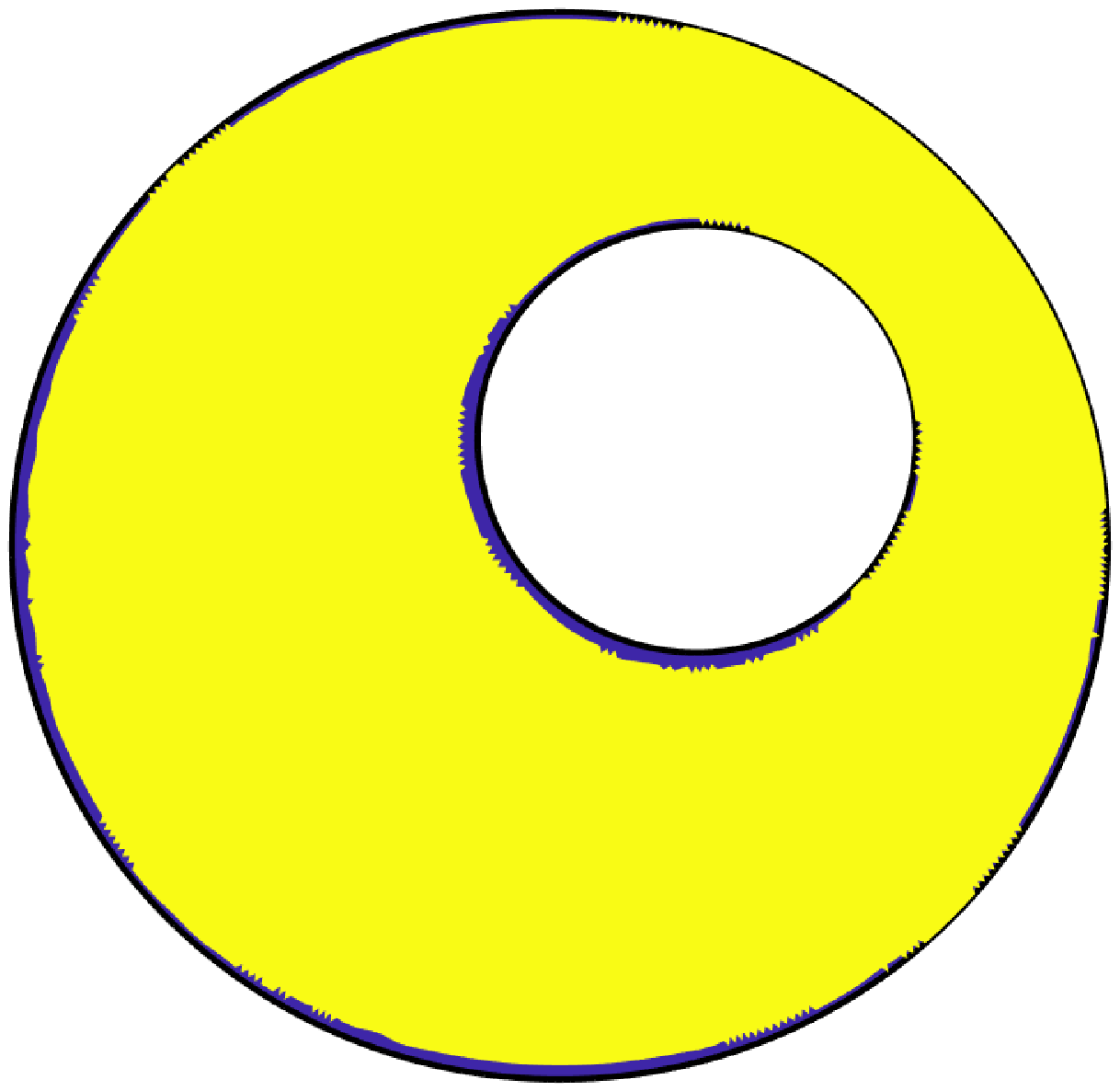}
        \end{tabular}
    \end{center}
    \caption{Region (blue color) covered by the cells in the set $\mesh^*$ for a mesh composed of $2266$ (left column), $9411$ (middle column), or 29496 cells, $\varepsilon=10^{-1}$ (top row), $\varepsilon=10^{-2}$ (middle row), or $\varepsilon=10^{-3}$ (bottom row).}\label{ex3_eps_theta_0.3_boundary_layer_mesh}
\end{figure}

Finally, to give some insight on the resolution of the boundary layer for $\varepsilon\in\{10^{-1},10^{-2},10^{-3}\}$, we flag the mesh cells as belonging to the boundary layer by means of the following criterion:
$$
\mesh^*:= \Big\{ K\in\mesh \;\big|\; \|\nabla^2 R_K(\hat{u}_K)\|_{L^\infty(K)} \geq \theta  \max_{\tilde{K} \in \mesh}\|\nabla^2 R_{\tilde{K}}(\hat{u}_{\tilde{K}} )\|_{L^\infty(\tilde{K})} \Big\},
$$
with the threshold parameter set here to $\theta:=0.3$,
and the $L^\infty$-norm estimated by computing the mean of the values taken
by the Hessian norm at the three vertices of $K$.
We report in Table~\ref{table: Maximum of Hessian} the two following quantities:
(i) the maximal value of the Hessian, $\max_{\tilde{K} \in \mesh}\|\nabla^2 R_{\tilde{K}}(\hat{u}_{\tilde{K}} )\|_{L^\infty(\tilde{K})}$;
(ii) the area of the boundary layer, $\sum_{K\in\mesh^*}|K|$.
First, we observe that the maximal value of the Hessian on the two finest meshes (which both resolve the boundary layer) are almost the same for all the values of $\varepsilon$. Moreover, the maximal value of the Hessian appears to scale as $\mathcal{O}(\varepsilon^{-\frac{3}{4}})$. This scaling is consistent with the expected scaling of the $H^2$-norm of the Hessian as $\mathcal{O}(\varepsilon^{-\frac12})$ and a boundary layer with surface scaling as $\mathcal{O}(\varepsilon^{\frac12})$. Furthermore,
the set $\mesh^*$ covers a region whose area decays a bit slower than the expected rate $\mathcal{O}(\varepsilon^{\frac12})$. This behavior indicates that somewhat finer meshes are still needed to fully resolve the geometric description of the boundary layer. This conclusion is corroborated in Figure \ref{ex3_eps_theta_0.3_boundary_layer_mesh}, where we show the region covered by the cells in $\mesh^*$ for $\varepsilon\in\{10^{-1},10^{-2},10^{-3}\}$
and the curved triangular meshes composed of 2266, 9411, or 29496 cells.
For $\varepsilon=0.1$, the set $\mesh^*$  contains not only cells close to the boundary but also cells in the interior. As $\varepsilon$ becomes smaller, the set $\mesh^*$ contains fewer and fewer cells in the interior of the domain, and for both $\varepsilon=10^{-2}$ and $10^{-3}$, the region covered by $\mesh^*$ is fully localized at the boundary. We notice, however, that even for $\varepsilon=10^{-3}$, there are boundary cells that are not flagged as members of $\mesh^*$; those cells are located in the part of the domain where the boundary of the inner disk is close to the boundary of the outer disk.


\begin{table}[!htb]
    \begin{center}
        \begin{tabular}{||c|c|c|c|c||}
            \hline
            mesh & 527 cells& 2266 cells& 9411 cells & 29496 cells \\
            \hline
            $h$&0.1416&0.0704&0.0344&0.0205\\
            \hline
            &\multicolumn{4}{c||}{Max Hessian} \\
            \hline
            $\varepsilon=10^{-1}$ & 6.94& 7.09&7.11 & 7.11\\
            \hline
            $\varepsilon=10^{-2}$  &38.13&40.48&40.88& 41.22\\
            \hline
            $\varepsilon=10^{-3}$ &88.63&128.33&147.96& 151.32\\
            \hline
             &\multicolumn{4}{c||}{Area covered by $\mesh^*$} \\
             \hline
             $\varepsilon=10^{-1}$ &5.50e-01 & 4.51e-01& 4.16e-01 &  3.86e-01  \\
             \hline
             $\varepsilon=10^{-2}$  &3.78e-01 &2.47e-01&2.16e-01 & 2.13e-01 \\
             \hline
             $\varepsilon=10^{-3}$ &4.38e-01 & 1.65e-02&1.12e-01& 1.22e-01\\
             \hline
        \end{tabular}
    \end{center}
    \caption{Test case with boundary layer.
    Maximal Hessian value 
    and area covered by $\mesh^*$ for $\varepsilon\in\{10^{-1},10^{-2},10^{-3}\}$
    and the curved triangular meshes composed of 527, 2266, 9411, and 29496 cells.}
    \label{table: Maximum of Hessian}
\end{table}


\section{Proof of main results} \label{sec: proof of main resutls}

In this section, we present the proofs of Lemma \ref{lem: stability and boundedness}, Lemma \ref{lem: approxmation}, Lemma \ref{lemma: consistency}, and Theorem~\ref{Theorem: main}.

\subsection{Proof of Lemma \ref{lem: stability and boundedness}}
\label{sec:proof_stability}

(1) We start with the lower bound in \eqref{local equivalent}. Choosing the test function $w = v_K$ in \eqref{eq: reconstruction} gives
\begin{align*}
\|\nabla  {v}_K\|^2_{K,\varepsilon} = {}&	(\nabla R_K\upi (v_K),\nabla v_K)_{K,\varepsilon} + (v_K -v_{\partial K} , \partial_n  v_K)_{\dKi}
+ (v_K  , \partial_n  v_K)_{\dKb} \\
&- \varepsilon \Big\{
(v_K -v_{\partial K} , \partial_n \Delta  v_K)_{\dKi}  - (\partial_n v_K - \gamma_{\partial K},  \partial_{nn}  v_K)_{\dKi}  - (\partial_t (v_K - v_{\partial K}),  \partial_{nt}  v_K)_{\dKi}\\
& - (v_K  , \partial_n \Delta  v_K)_{\dKb} + (\nabla  v_K ,  \nabla \partial_{n} v_K)_{\dKb}
\Big\}.
\end{align*}
Using the Cauchy--Schwarz inequality, the discrete inverse inequalities \eqref{H1 inverse inequality},  \eqref{trace inverse inequality}, \eqref{inverse inequality tangential}, and that $\partial_{nn}v_K \in \mathbb{P}^{k}(\dKi)$ to introduce the projection $\Pi_{\dKi}^k$, we infer that
\begin{multline*}
\|\nabla  {v}_{K}\|_{K,\varepsilon} \leq \|\nabla  R_K\upi (\hat{v}_K)\|_{K,\varepsilon}
+ C\Big(  h_K^{-\frac12} \| v_K -v_{\partial K}\|_{\dKi} +  h_K^{-\frac12} \| v_K \|_{\dKb}\\
+ \varepsilon^{\frac12} \Big\{ h_K^{-\frac32} \| v_K -v_{\partial K}\|_{\dKi}
+ h_K^{-\frac12} \|\Pi_{\dKi}^k(\partial_n v_K - \gamma_{\partial K} )\|_{\dKi}
+h_K^{-\frac32} \|v_K\|_{\dKb}+h_K^{-\frac12} \|\nabla  v_K\|_{\dKb} \Big\} \Big).
\end{multline*}
Since $\sigma_K = \max\{1,\varepsilon h^{-2}_K\}$, this implies that
\[
\|\nabla  {v}_{K}\|_{K,\varepsilon} \le \|\nabla  R_K\upi (\hat{v}_K)\|_{K,\varepsilon}
+ C\big( S\upi_{\dK}(\hat{v}_K,\hat{v}_K)+S\upb_{\dK}(v_K,v_K)\big)^{\frac12}.
\]
It remains to bound the four boundary terms on the right-hand side of \eqref{H2_seminorm_elem}. It is clear that
\[
\sigma_K h_K^{-1}\|v_{\partial K} - v_K\|_{\dKi}^2
+ \sigma_K h_K^{-1}  \|v_{K}\|^2_{\dKb}
+ \varepsilon h_K^{-1}  \|\nabla v_{K}\|^2_{\dKb} \le
S\upi_{\dK}(\hat{v}_K,\hat{v}_K)+S\upb_{\dK}(v_K,v_K),
\]
so that it only remains to bound ${\sigma_K}{h_K} \|\gamma_{\partial K} - \partial_n v_K\|_{\dKi}^2$. To this purpose, using a Poincar\'e--Steklov inequality followed by a discrete trace inequality on $\dKi$, we observe that
\begin{align*}
\sigma_K^{\frac12} h_K^{\frac12} \| \gamma_{\dK} - \partial_n v_K\|_{\dKi} \le {}& \sigma_K^{\frac12}h_K^{\frac12} \|\Pi_{\dKi}^{k}(\gamma_{\dK} - \partial_n v_K)\|_{\dKi}
+ \sigma_K^{\frac12}h_K^{\frac12}\|\partial_nv_K-\Pi_{\dKi}^{k}(\partial_n v_K)\|_{\dKi} \\
\le {}& S\upi_{\dK}(\hat{v}_K,\hat{v}_K)^{\frac12} + C \sigma_K^{\frac12}h_K^{\frac32} \|\partial_{nt}v_K\|_{\dK} \\
\le {}& S\upi_{\dK}(\hat{v}_K,\hat{v}_K)^{\frac12} + C\sigma_K^{\frac12}h_K\|\nabla^2 v_K\|_K.
\end{align*}
The definition of $\sigma_K$ implies that $\sigma_K^{\frac12}\le 1+ \varepsilon^{\frac12}h_K^{-1}$, so that
\[
\sigma_K^{\frac12}h_K\|\nabla^2 v_K\|_K\le h_K\|\nabla^2 v_K\|_K + \varepsilon^{\frac12}\|\nabla^2 v_K\|_K \le C\|\nabla v_K\|_K + \varepsilon^{\frac12}\|\nabla^2 v_K\|_K \le C\|\nabla v_K\|_{K,\varepsilon},
\]
where we used the discrete inverse inequality \eqref{H1 inverse inequality}.
Owing to the above bound on $\|\nabla v_K\|_{K,\varepsilon}$, we infer that
\[
\sigma_K^{\frac12} h_K^{\frac12} \| \gamma_{\dK} - \partial_n v_K\|_{\dKi} \le
C\Big(\|\nabla  R_K\upi (\hat{v}_K)\|_{K,\varepsilon}
+ \big( S\upi_{\dK}(\hat{v}_K,\hat{v}_K)+S\upb_{\dK}(v_K,v_K)\big)^{\frac12}\Big).
\]
Combining the above bounds, we conclude that the lower bound in \eqref{local equivalent} holds true.

(2) Let us now establish the upper bound in \eqref{local equivalent}. This time we choose $w = R_K\upi (\hat{v}_K)$ in \eqref{eq: reconstruction} and proceeding as above yields
\[
\|\nabla  R_K\upi (\hat{v}_K)\|_{K,\varepsilon} \le \|\nabla  {v}_{K}\|_{K,\varepsilon}
+ C\big( S\upi_{\dK}(\hat{v}_K,\hat{v}_K)+S\upb_{\dK}(v_K,v_K)\big)^{\frac12}.
\]
Moreover, it is clear that
\begin{multline*}
S\upi_{\dK}(\hat{v}_K,\hat{v}_K)+S\upb_{\dK}(v_K,v_K) \le
\sigma_K h_K^{-1}\|v_{\partial K} - v_K\|_{\dKi}^2
+ \sigma_K h_K \| \gamma_{\dK} - \partial_n v_K\|_{\dKi}^2 \\
+ \sigma_K h_K^{-1}  \|v_{K}\|^2_{\dKb}
+ \varepsilon h_K^{-1}  \|\nabla v_{K}\|^2_{\dKb},
\end{multline*}
since $\|\Pi^{k}_{\dKi} ( \partial_n  v_K - \gamma_{\partial K} ) \|_{\dKi}
\leq \sigma_K h^{\frac12}_K  \|\partial_{n} v_{K} - \gamma_{\partial K}  \|_{\dKi}$.
This completes the proof.

\subsection{Proof of Lemma \ref{lem: approxmation}}
\label{sec:proof_approximation}

(1) Let us first bound $\|v - \mathcal{E}_K (v)\|_{\sharp, K}$. The triangle inequality
followed by the discrete inverse inequalities \eqref{trace inverse inequality} and
\eqref{H1 inverse inequality} implies that
\begin{align*}
\|v - \mathcal{E}_K (v)\|_{\sharp, K} &\leq \|v - \Pi_K^{k+2}(v)\|_{\sharp, K} + \|\mathcal{E}_K(v) - \Pi_K^{k+2}(v)\|_{\sharp, K} \\
&\leq \|v - \Pi_K^{k+2}(v)\|_{\sharp, K} + C\|\nabla (\mathcal{E}_K(v) - \Pi_K^{k+2}(v))\|_{K,\varepsilon},
\end{align*}
so that we only need to bound the last term on the right-hand side.
Straightforward algebra shows that, for all $\xi \in \mathbb{P}^{k+2}(K)$,
\begin{align*}
(\nabla \mathcal{E}_K(v)-\nabla \Pi_{K}^{k+2}(v), \nabla \xi)_{K,\varepsilon}
= {}& -( \Pi_{K}^{k+2}(v) -v , \partial_n   \xi)_{\dK}\\
& + \varepsilon	\Big\{ (\Pi_{K}^{k+2}(v) -v , \partial_n \Delta  \xi)_{\dK} - (\partial_n (\Pi_{K}^{k+2}(v)  - v),  \partial_{nn}  \xi)_{\dK} \\
&- (\partial_t (\Pi_{K}^{k+2}(v) - \Pi^{k+2}_{\dKi} (v)),  \partial_{nt}  \xi)_{\dKi}
- (\partial_t (\Pi_{K}^{k+2}(v) - v),  \partial_{nt}  \xi)_{\dKb}\Big\}.
\end{align*}
Choosing $\xi =  \mathcal{E}_K(v) - \Pi_{K}^{k+2}(v)$ and
using the discrete inverse inequalities \eqref{trace inverse inequality}, \eqref{H1 inverse inequality}, \eqref{inverse inequality tangential} gives
\begin{align*}
\|\nabla (\mathcal{E}_K(v)  - \Pi_{K}^{k+2}(v))\|_{K,\varepsilon} \leq {}&
C\Big( \sigma_K^{\frac12} h_K^{-\frac12}\|v-\Pi_{K}^{k+2}(v)\|_{\dK}
+ \varepsilon^{\frac12} h_K^{-\frac12} \|\nabla(v-\Pi_{K}^{k+2}(v))\|_{\dK}\\
&+\varepsilon^{\frac12} h_K^{-\frac12}\|\partial_t (\Pi_{K}^{k+2}(v) - \Pi^{k+2}_{\dKi} (v))\|_{\dKi}\Big).
\end{align*}
The first two terms on the right-hand side are bounded using \eqref{eq:tr_PK2} and \eqref{eq:tr_sigma} leading to
\[
\sigma_K^{\frac12} h_K^{-\frac12}\|v-\Pi_{K}^{k+2}(v)\|_{\dK}
+ \varepsilon^{\frac12} h_K^{-\frac12} \|\nabla(v-\Pi_{K}^{k+2}(v))\|_{\dK}
\le C \|\nabla (v-\Pi_{K}^{k+2}(v))\|_{K,\varepsilon}.
\]
Moreover, for the last term on the right-hand side, proceeding as in \cite{DongErn2021biharmonic}, we invoke the discrete trace inequality \eqref{inverse inequality tangential} and observe that $\Pi_{K}^{k+2}(v) - \Pi^{k+2}_{\dKi} (v)=\Pi^{k+2}_{\dKi}(\Pi_{K}^{k+2}(v)-v)$. Since $\Pi^{k+2}_{\dKi}$ is $L^2$-stable, we conclude that
\[
\varepsilon^{\frac12} h_K^{-\frac12}\|\partial_t (\Pi_{K}^{k+2}(v) - \Pi^{k+2}_{\dKi} (v))\|_{\dKi}
\leq C\varepsilon^{\frac12} h_K^{-\frac32}\|v-\Pi_{K}^{k+2}(v)\|_{\dKi}.
\]
Invoking \eqref{eq:tr_PK2}, this yields
\[
\varepsilon^{\frac12} h_K^{-\frac12}\|\partial_t (\Pi_{K}^{k+2}(v) - \Pi^{k+2}_{\dKi} (v))\|_{\dKi}
\leq C\|\nabla (v-\Pi_{K}^{k+2}(v))\|_{K,\varepsilon}.
\]
We have thus shown that
\[
\|\nabla (\mathcal{E}_K(v)  - \Pi_{K}^{k+2}(v))\|_{K,\varepsilon} \le C\|\nabla (v-\Pi_{K}^{k+2}(v))\|_{K,\varepsilon}.
\]
Putting the above bounds together shows that $\|v - \mathcal{E}_K (v)\|_{\sharp, K} \leq C
\|v - \Pi_K^{k+2}(v)\|_{\sharp, K}$.

(2) Let us now bound $S\upi_{\dK}(\mathcal{\hat{I}}^k_K(v),\mathcal{\hat{I}}^k_K(v))$.
We have
\begin{align*}
S\upi_{\dK}(\mathcal{\hat{I}}^k_K(v),\mathcal{\hat{I}}^k_K(v))
={}& \sigma_K h_K^{-1} \| \Pi^{k+2}_{\partial K} (v)- \Pi^{k+2}_K (v) \|^2_{\dKi} + \sigma_K h_K \| \Pi^k_{\dKi}(\Pi^k_{\dKi} (\partial_{n}v)- \partial_{n} (\Pi^{k+2}_K (v))) \|^2_{\dKi} \\
\leq {}& \sigma_K h_K^{-1} \| v- \Pi^{k+2}_K (v) \|^2_{\dKi} + \sigma_K h_K \| \partial_{n}(v- \Pi^{k+2}_K (v)) \|^2_{\dKi},
\end{align*}
since $\Pi^{k+2}_{\dKi}(\Pi^{k+2}_K (v))=\Pi^{k+2}_K (v)_{| \dKi}$, $\Pi^{k+2}_{\partial K}$ is $L^2$-stable, and $\Pi^k_{\dKi}\circ \Pi^k_{\dKi}=\Pi^k_{\dKi}$.
The first term on the right-hand side is bounded by means of \eqref{eq:tr_sigma}, yielding
\begin{align*}
\sigma_K h_K^{-1} \| v- \Pi^{k+2}_K (v) \|^2_{\dKi} &\leq
C \|\nabla (v-\Pi_{K}^{k+2}(v))\|_{K,\varepsilon}^2.
\end{align*}
Moreover, we have
\begin{align*}
\sigma_K h_K \| \partial_{n}(v-\Pi^{k+2}_K (v)) \|^2_{\dKi}
&\leq \varepsilon h_K^{-1}\|  \partial_{n}(v-  \Pi_K^{k+2}(v)) \|^2_{\dKi}+ h_K \|  \partial_{n}(v-  \Pi_K^{k+2}(v)) \|^2_{\dKi} \\
&\leq C \|\nabla(v-  \Pi_K^{k+2}(v))\|_{K,\varepsilon}^2 + h_K \|  \partial_{n}(v-  \Pi_K^{k+2}(v)) \|^2_{\dKi} \\
&\leq C\|v- \Pi^{k+2}_K (v)\|_{\sharp,K}^2,
\end{align*}
where we used \eqref{eq:tr_PK2} and the definition of the $\|{\cdot}\|_{\sharp,K}$-norm.
Putting the above bounds together, we infer that $S\upi_{\dK}(\mathcal{\hat{I}}^k_K(v),\mathcal{\hat{I}}^k_K(v))\leq C\|v- \Pi^{k+2}_K (v)\|_{\sharp,K}^2$.

(3) Finally, let us bound $S\upb_{\dK}(v-\Pi_K^{k+2}(v),v-\Pi_K^{k+2}(v))$. We have
\begin{align*}
S\upb_{\dK}(v-\Pi_K^{k+2}(v),v-\Pi_K^{k+2}(v)) &=
\sigma_Kh_K^{-1}\|v-\Pi_K^{k+2}(v)\|_{\dKb}^2
+ \varepsilon h_K^{-1} \|\nabla(v-\Pi_K^{k+2}(v))\|_{\dKb}^2,
\end{align*}
and the two terms on the right-hand side can be bounded by invoking the same arguments as in Step~(2). This concludes the proof.

\subsection{Proof of Lemma \ref{lemma: consistency}}
\label{sec:proof_consistency}

Recalling the definition \eqref{def: lifting} of the lifting operator and the definition \eqref{def: stabilisation Boundary} of the boundary stabilization operator, we observe that
\[
\ell_h(\hat{w}_h) = \su \Big\{
(f,w_K)_K - (\nabla \mathcal{L}_K(u), \nabla R_K\upi(\hat{w}_K))_{K,\varepsilon}
+ S_{\dK}\upb(u,w_K) \Big\}.
\]
This implies that
\begin{equation} \label{eq:consist}
\langle \delta_h,\hat{w}_h\rangle = \su \Big\{
(f,w_K)_K - (\nabla  \mathcal{E}_K(u), \nabla R_K\upi(\hat{w}_K))_{K,\varepsilon}
- S_{\dK}\upi(\mathcal{\hat{I}}^k_K(u),w_K)- S_{\dK}\upb(\Pi_K^{k+2}(u)-u,w_K) \Big\}.
\end{equation}
We bound the four terms on the right-hand side of \eqref{eq:consist}.
The first two terms are combined together
by using that $f=\varepsilon \Delta^2 u - \Delta u$, integration by parts, the regularity
of the exact solution, and the definition of $R_K\upi(\hat{w}_K)$.
Let us set $\eta|_K:= u|_K - \mathcal{E}_K (u|_K)$ for all $K\in\mesh$.
Proceeding as in \cite{DongErn2021biharmonic} for the fourth-order operator and as in
\cite{DiPEL:14,BurErCiDe:21} for the second-order operator, we obtain
\begin{align*}
&\su \Big\{
(f,w_K)_K - (\nabla \mathcal{E}_K(u|_K), \nabla R_K\upi(\hat{w}_K))_{K,\varepsilon} \Big\} \\
&= \su \bigg\{(\nabla \eta, \nabla {w}_K)_{K,\varepsilon} - \big(\partial_n\eta , w_K - w_{\partial K} \big)_{\dKi}
-\big(\partial_n \eta , w_{K} \big)_{\dKb}\\
& \quad + \varepsilon \Big( \big( \partial_n \Delta\eta , w_K - w_{\partial K} \big)_{\dKi}
- \big(\partial_{nn}  \eta , \partial_{n} w_K  - \chi_{\partial K}\big)_{\dKi}
- \big( \partial_{nt}  \eta , \partial_{t} (w_K - w_{\partial K}) \big)_{\dKi} \\
& \quad + \big(\partial_n \Delta \eta , w_{K} \big)_{\dKb}- \big(\nabla \partial_{n} \eta, \nabla w_{K} \big)_{\dKb} \Big) \bigg\}.
\end{align*}
Invoking the Cauchy--Schwarz inequality and the discrete inverse inequalities from Lemma~\ref{lemma: Inverse inequality}, we infer that
\[
\su \Big\{
(f,w_K)_K - (\nabla\mathcal{E}_K(u |_K), \nabla R_K\upi(\hat{w}_K))_{K,\varepsilon} \Big\}
\leq C
\left( \su \| \eta \|^2_{\sharp,K}\right)^{\frac12} \|\hat{w}_h\|_{\fes}.
\]
Moreover, the third and fourth terms on the right-hand side of \eqref{eq:consist} are estimated by means of Lemma~\ref{lem: approxmation}. Putting everything together concludes the proof.

\subsection{Proof of Theorem~\ref{Theorem: main}}
\label{sec:proof_main}

(1) Let us set $\hat{e}_h^k : = \mathcal{\hat{I}}_h^k(u)-\hat{u}_h \in \fes$, so that $\hat{e}_K^k : = \mathcal{\hat{I}}_K^k(u|_K) - \hat{u}_K \in \fesE$ for all $K\in \mesh$. The property \eqref{coercivity} and the identity $a_h(\hat{e}_h^k,\hat{e}_h^k)= - \langle \delta_h,\hat{e}_h^k \rangle$ imply that
\begin{equation*}
\alpha \su \|\nabla (R_K\upi(\hat{e}_K^k))\|_{K,\varepsilon}^2
\le \alpha \|\hat{e}^k_h\|_{\fes}^2 \leq
a_h(\hat{e}_h^k,\hat{e}_h^k) = \langle \delta_h,\hat{e}_h^k \rangle.
\end{equation*}
Owing to the Lemma \ref{lemma: consistency}, we infer that
$$
\su \|\nabla (R_K\upi(\hat{e}_K^k))\|_{K,\varepsilon}^2  \leq C \su \|u-\Pi_K^{k+2} (u)\|^2_{\sharp,K}.
$$
By adding and subtracting $R_K\upi({\mathcal{\hat{I}}_K^k(u|_K)})$, we have
$$
u_{|K} -  R_K\upi(\hat{u}_K) - \mathcal{L}_K(u|_K)
=(u_{|K} - \mathcal{E}_K(u|_K))+ R_K\upi(\hat{e}_K^k).
$$
Using the  triangle inequality and Lemma~\ref{lem: approxmation}, the error estimate~\eqref{abstract error bound} is derived.

(2) Let us assume that $u|_K\in H^{k+3}(K)$ for all $K\in\mesh$ if $k\ge1$
and that $u|_K\in H^{4}(K)$ for all $K\in\mesh$ if $k=0$. Combining
the multiplicative trace inequality \eqref{trace inequality} with
the polynomial approximation result \eqref{Polynomial approximation} shows that
$\|u-\Pi_K^{k+2}(u)\|_{\sharp,K} \le C\sigma_K^{\frac12}h_K^{k+2}|u|_{H^{k+3}(K)}$
if $k\ge1$ and
$\|u-\Pi_K^{k+2}(u)\|_{\sharp,K} \le C\sigma_K^{\frac12}h_K^{k+2}(|u|_{H^{k+3}(K)}+h_K|u|_{H^4(K)})$
if $k=0$. Summing these estimates over the mesh cells and using the error estimate derived in Step (1) proves \eqref{eq:err2}.

\bibliographystyle{siam}
\bibliography{HHO_N}

\end{document}